\documentclass[11pt]{article}
\usepackage{mathbbold}
\usepackage{bbm}

\usepackage{mathrsfs}
\usepackage{amsfonts}
\usepackage{amssymb}
\usepackage{amsfonts,amssymb,mathrsfs,amsmath,amssymb,theorem,float,xypic}
\usepackage{color}

\newtheorem{theorem}{Theorem}[section]

\newtheorem{cor}[theorem]{Corollary}
\newtheorem{lemma}[theorem]{Lemma}
\newtheorem{definition}[theorem]{Definition}
\newtheorem{remark}[theorem]{Remark}
\newtheorem{example}[theorem]{Example}

\floatstyle{ruled}
\newfloat{algorithm}{t}{loa}
\floatname{algorithm}{Algorithm}
\newcommand{\Inp}[1]
  {\noindent\begin{tabular}{@{}p{1.8cm}@{}p{13.2cm}@{}}
   {\bf Input: }&#1 \end{tabular}}
\newcommand{\Outp}[1]
  {\noindent\begin{tabular}{@{}p{1.8cm}@{}p{13.2cm}@{}}
   {\bf Output: }&#1 \end{tabular}}

\def\imply{\Rightarrow}
\def\qed{{\vrule height5pt width2pt depth2pt}}

\def\and{\cap}

\def\bref#1{(\ref{#1})}
\def\proof{{\noindent\em Proof:} }

\newcommand{\SPC}{\hspace*{15pt}}
\newcommand{\qedd}{\hspace*{\fill}$\Box$\medskip}
\newcommand{\bm}[1]{\mbox{\boldmath{$#1$}}}

\floatstyle{ruled}
\newfloat{algorithm}{t}{loa}
\floatname{algorithm}{Algorithm}

\def\st#1{{\em{#1}}}

\def\X{{\mathbb{X}}}
\def\Y{{\mathbb{Y}}}

\def\V{{\mathbb{V}}}
\def\L{{\mathbb{L}}}
\def\J{{\mathcal {J}}}
\def\T{{\mathbb{T}}}

\def\FB{{\mathbb{F}}}
\def\RB{{\mathbb{R}}}

\def\XB{{\mathbf{X}}}

\def\AH{{\mathbb{A}}}
\def\IH{{\mathbb{I}}}
\def\CH{{\mathbb{C}}}

\def\A{{\mathcal A}}
\def\B{{\mathcal B}}
\def\C{{\mathcal C}}
\def\D{{\mathcal D}}

\def\MI{{\mathcal{M}}}

\def\D{{\mathbb D}}
\def\G{{\mathbb G}}

\def\P{{\mathbb P}}
\def\Q{{\mathbb Q}}

\def\bu{{\mathbf{u}}}

\def\balpha{{\bm\alpha}}
\def\bbeta{{\bm\beta}}
\def\bgamma{{\bm\gamma}}
\def\beps{{\bm\epsilon}}

\def\falpha{{\bm\bbalpha}}
\def\fbeta{{\bm\bbbeta}}

\def\prem{\hbox{\rm prem}}

\def\grem{\hbox{\rm grem}}

\def\sat{\hbox{\rm{sat}}}

\def\max{\hbox{\rm{max}}}

\def\Im{\hbox{\rm{Im}}}

\def\span{\hbox{\rm{Span}}}
\def\lcm{\hbox{\rm{lcm}}}
\def\gcd{\hbox{\rm{gcd}}}

\def\gcd{\hbox{\rm{gcd}}}
\def\deg{\hbox{\rm{deg}}}

\def\init{\hbox{\rm{I}}}
\def\ord{\hbox{\rm{ord}}}

\def\H{\hbox{\rm{H}}}
\def\lead{\hbox{\rm{ld}}}

\def\dim{\hbox{\rm{dim}}}

\def\lv{\hbox{\rm{lvar}}}
\def\mod{\hbox{\rm{mod}}}

\def\ord{\hbox{\rm{ord}}}
\def\rk{\hbox{\rm{rk}}}

\def\coeff{\hbox{\rm{coeff}}}

\def\Spec{\hbox{\rm{Spec}}}

\def\I{\mathcal{I}}
\def\CI{{\mathcal{I}}}

\def\LT{{\bf LT}}
\def\c{{\bf c}}
\def\b{{\bf b}}
\def\f{{\bf f}}
\def\g{{\bf g}}
\def\h{{\bf h}}
\def\e{{\bf e}}
\def\s{{\bf s}}
\def\a{{\bf a}}

\def\Z{{\mathbb{Z}}}
\def\Q{{\mathbb{Q}}}
\def\N{{\mathbb{N}}}

\def\F{{\mathcal{F}}}
\def\E{{\mathcal{E}}}

\def\dtrdeg{\hbox{$\triangle$\rm{tr.deg}}}
\def\Jac{\hbox{\rm{Jac}}}

\def\ix{{\rm{i}}}

\def\K{{\mathbb{K}}}

\def\grem{\hbox{\rm grem}}
\def\and{\cap}

\newcounter{bean}
\def\bl{\begin{list}{Step \arabic{bean}}{\usecounter{bean}}\labelwidth=34pt}
\def\el{\end{list}}

\def\deg{{\rm deg}}
\def\lvar{{\rm lvar}}

\def\init{{\rm I}}

\def\normalization1{{\rm normalization1}}
\def\normalization{{\rm normalization}}

\def\irrfactor1{{\rm irrfactor1}}
\def\irrfactor{{\rm irrfactor}}

\def\sat{{\rm sat}}

\def\gb{Gr\"obner基}

\def\grem{\hbox{\rm{grem}}}

\def\m{\mathfrak{m}}
\def\p{\mathfrak{p}}

\def\st{\hbox{ {\rm{s.t.}} }}

\def\Zx{\Z[x]}
\def\Zxn{\Z[x]^n}

\def\fb{{\mathbbm{f}}}
\def\gb{{\mathbbm{g}}}
\def\hb{{\mathbbm{h}}}

\def\mb{{\mathbbm{m}}}

\def\lr#1{{\langle{#1}\rangle}}

\begin{document}
\title{ Binomial Difference Ideal and Toric Difference Variety\thanks{Partially
       supported by a National Key Basic Research Project of China (2011CB302400) and  by grants from NSFC
       (60821002,11101411).}}
\author{Xiao-Shan Gao, Zhang Huang, Chun-Ming Yuan\\
 KLMM,  Academy of Mathematics and Systems Science\\
 Chinese Academy of Sciences, Beijing 100190, China\\
 \{xgao,cmyuan\}@mmrc.iss.ac.cn}
\date{}
\maketitle

\begin{abstract}\noindent
In this paper, the concepts of binomial difference ideals and toric
difference varieties are defined and their properties are proved.
Two canonical representations for Laurent binomial difference ideals
are given using the reduced Gr\"obner basis of $\Zx$-lattices and
regular and coherent difference ascending  chains, respectively.
Criteria for a Laurent binomial difference ideal to be reflexive,
prime, well-mixed, perfect, and toric are given in terms of their
support lattices which are $\Zx$-lattices. The reflexive,
well-mixed, and perfect closures of a Laurent binomial difference
ideal are shown to be binomial. Four equivalent definitions for
toric difference varieties are presented. Finally, algorithms are
given to check whether a given Laurent binomial difference ideal
$\I$ is reflexive, prime, well-mixed,  perfect, or toric, and in the
negative case, to compute the reflexive, well-mixed, and perfect
closures of $\I$. An algorithm is given to decompose a finitely
generated perfect binomial difference ideal as the intersection of
reflexive prime binomial difference ideals.


\end{abstract}

\tableofcontents

\section{Introduction}
\label{sec-intro}

The theory of toric varieties has been extensively studied since its
foundation in the early 1970s by Demazure \cite{dema},  Miyake-Oda
\cite{oda}, Mumford et al. \cite{mum}, and Satake \cite{sata}, due
to its deep connections with polytopes, combinatorics, symplectic
geometry, topology, and its applications in physics, coding
theory, algebraic statistics, and hypergeometric functions
\cite{cox-2010,fulton,gelfand,oda-book}.
Toric varieties are often used as an effective testing ground for general theories of algebraic geometry.
In \cite{es-bi}, Eisendbud and Sturmfels initiated the study of
binomial ideals which enriched the algebraic aspects of toric
varieties and leaded to more applications \cite{miller1,pachter1}.

%
In this paper, we initiate the study of binomial difference ideals and toric
difference varieties and hope that they will play similar roles in
difference algebraic geometry to their algebraic counterparts in
algebraic geometry.
Difference algebra and difference algebraic geometry were founded by
Ritt \cite{ritt-dd1} and Cohn \cite{cohn}, who aimed to study
algebraic difference equations as algebraic geometry to polynomial
equations. Besides the early achievements, some of the major recent
advances in this field include the difference Galois theory and its
applications \cite{put-dd}, the establishment of the theory of
difference scheme and the proof of the Jacobi bound for difference
equations  \cite{Hrushovski1}, and the explicit
description of invariant varieties under coordinatewise difference
operators defined with univariate polynomials  \cite{med1}.

We now describe the main results of this paper.
In Section \ref{sec-zxm}, we prove basic properties of $\Z[x]$-lattices. By a $\Z[x]$-lattice, we mean a $\Z[x]$-module in
$\Z[x]^n$. $\Z[x]$-lattices play the same role as $\Z$-lattices in
the study of binomial ideals and toric varieties. Here, $x$ is used
to denote the difference operator $\sigma$. For instance,
$a^3\sigma(a)^2$ is denoted as $a^{2x+3}$.  Many properties of
binomial difference ideals can be described or proved with the help
of $\Zx$-lattices.
Since $\Zx$ is not a PID, the Hermite normal form for a matrix with
entries in $\Zx$ does not exist. In this section, we introduce the
concept of generalized Hermite normal form and show that a matrix is
a generalized Hermite normal form if and only if its columns form a
reduced Gr\"obner basis for a $\Zx$-lattice.
In general, a $\Zx$-lattice is not a free $\Zx$-module. We prove
that the kernel of a matrix with entries in $\Zx$ is a free $\Zx$
module, which plays a key role in the study of toric difference
varieties.

In Section \ref{sec-lbi}, we prove basic properties of Laurent
binomial difference ideals.
Gr\"obner bases play an important role in the study of binomial
ideals \cite{es-bi}. In general, a binomial difference ideal is not
finitely generated and does not have a finite Gr\"obner basis.
Instead, the theory of characteristic set for difference polynomial
systems \cite{gao-dcs} is used for similar purposes.
In Sections 4.2 and 4.3, we give canonical forms for Laurent
binomial difference ideals in terms of $\Zx$-lattices and characteristic sets.
Let $\F$ be a difference field with a difference operator $\sigma$,
$\F^*=\F\setminus\{0\}$, $\Y=\{y_1,\ldots,y_n\}$ a set of difference
indeterminates, and $\F\{\Y^{\pm}\}$ the ring of Laurent difference
polynomials in $\Y$. Then
\begin{theorem}\label{th-m1}
 $\I\subset\F\{\Y^{\pm}\}$ is a proper Laurent binomial
difference ideal if and only if
\begin{description}
\item[(1)]$\I=[\A]$, where $\A=\: \Y^{\f_{1}}-c_1, \ldots, \Y^{\f_{s}}-c_s$,
$\f_i\in\Zxn$, $c_i\in\F^*$, $\fb=\{\f_{1}, \ldots, \f_{s}\}$ is a
reduced Gr\"obner basis of a $\Zx$-lattice, and $[\A]\ne[1]$.
\item[(2)] $\I=[\A]$, where $\A=\: \Y^{\f_{1}}-c_1, \ldots, \Y^{\f_{s}}-c_s$,
$\f_i\in\Zxn$, $c_i\in\F^*$, and $\A$ is a regular and coherent
difference ascending chain.
\item[(3)]
$\I=\I(\rho)=[\Y^{\f}-\rho(\f)\,|\, \f\in L_{\rho}]$, where $\rho$
is a homomorphism from a $\Z[x]$-lattice $L_{\rho}$ generated by
$\fb$ to the multiplicative group $\F^{\ast}$ satisfying
$\rho(\f_i)=c_i$ and $\rho(\sigma(\f))=\sigma(\rho(\f))$ for $\f\in
L_\rho$.
\end{description}
In (1) and (2), $\A$ is a characteristic set of the difference ideal
$\I$.
\end{theorem}

In Sections \ref{sec-pr1} and \ref{sec-pr3}, we give criteria for a
Laurent binomial difference ideal to be prime, reflexive,
well-mixed, and perfect in terms of its lattice support.
\begin{theorem}\label{th-m2}
Let $\I$ be a proper Laurent binomial difference ideal and
$L=\{\f\,|\,\Y^\f-c_\f\in\I\}$ the {\em support lattice} of $\I$. If
$\F$ is algebraically closed and inversive, then
\begin{description}
\item[(1)] $\I$ is prime if and only if $L$ is $\Z$-saturated, that
is, $k\f\in L$ implies $\f\in L$ for $k\in\N$ and $\f\in\Zxn$.

\item[(2)] $\I$ is reflexive if and only if $L$ is $x$-saturated,
that is, $\f^x\in L$ implies $\f\in L$ for $\f\in\Zxn$.

\item[(3)] If the well-mixed closure of $\I$ is not
$[1]$, then $\I$ is well-mixed if and only if $L$ is
{\rm{M}}-saturated, that is, $k\f\in L$ implies $(x-o_k)\f\in L$ for
$k\in\N$ and $\f\in\Zxn$, where $o_k\in\N$ is a number determined by
the  difference field $\F$.

\item[(4)] If the perfect closure $\{\I\}$ of $\I$ is not
$[1]$, then  $\I$ is perfect if and only if $L$ is
{\rm{M}}-saturated and $x$-saturated.
\end{description}
\end{theorem}
The criterion for prime ideals is similar to the algebraic case, but
the criteria for reflexive, well-mixed, and perfect difference ideals
are unique to difference algebra and are first proposed in this
paper. In particular, the criterion for well-mixed difference ideals is quite intriguing.

Based on the above theorem, it is shown that the reflexive,
well-mixed, and perfect closures of a Laurent binomial difference
ideal $I$ with support lattice $L$ are still binomial, whose support
lattices are the $x$-, $M$-, and the $x$-$M$-saturation lattice
of $L$, respectively.
It is further shown that any perfect Laurent binomial difference
ideal $\I$ can be written as the intersection of Laurent reflexive
prime binomial difference ideals whose support lattices are the 
$x$-$\Z$- saturation of the support lattice of $\I$.
Since difference polynomial rings are not Notherian, binomial
difference ideals do not have an analog of the primary decomposition
theorem in the algebraic case.
%
%

In Section \ref{sec-bi}, binomial difference ideals are studied. It
is shown that a large portion of the properties for binomial ideals
proved in \cite{es-bi} can be easily extended to the difference case.
We also identify a class of normal binomial difference ideals which
are in a one to one correspondence with Laurent difference binomial
ideals.
%
%
With the help of this correspondence, properties proved for Laurent
binomial difference ideals can be extended to the non-Laurent case.

In Section \ref{sec-tv}, four equivalent definitions for difference
toric varieties are given. A difference variety is called toric if
it is the Cohn closure of the values of a set of Laurent difference
monomials. It is proved that
\begin{theorem}\label{th-m3} A difference variety $X$ is toric if and only
if one of the following properties is valid.
\begin{description}
\item[(1)] $X\cong\Spec^{\sigma}(\F\{M\})$,  where $M$ is a finitely generated affine $\N[x]$-module in $\Zx^m$ and $\F\{M\}=$ $\{\sum_{\balpha \in M}$
$a_{\balpha}\T^{\balpha}\,|\, a_{\balpha}\in \F, a_{\balpha}\ne0
\hbox{ for finitely many } \balpha\}$ for a set of difference
indeterminates $\T=\{t_1,\ldots,t_m\}$.

\item[(2)] The defining ideal of $X$ is a  toric difference ideal, that is,
 $\IH(X) = [\Y^{\f^+}-\Y^{\f^-}\,|\,\f\in L]$, where $L\subset\Zxn$ is a
 $\Zx$\hbox{-saturated} $\Zx$-lattice and $\f^+,\f^-$ are
 the positive and negative parts of $\f$, respectively.

\item[(3)] $X$ contains a difference torus $T^{*}$ as a Cohn open subset
and with a group action of $T^{*}$ on $X$ extending the natural
group action of $T^{*}$ on itself.

\end{description}
\end{theorem}

In the algebraic case, any prime binomial ideal of the form
$(\Y^{\f^+}-\Y^{\f^-}\,|\,\f\in L)$ has a monomial parametrization
and hence is toric \cite{es-bi}. In the difference case, this is not
valid (see Example \ref{ex-ig11}).
Also different from the algebraic case, the difference torus is not
necessarily isomorphic to $(\AH^*)^n$ (see Example \ref{ex-ext1}), 
and this makes the definition of difference torus more complicated.

It is shown that the difference sparse resultant can be defined as the difference Chow form of a difference toric variety, and a Jacobi style order bound for a difference toric variety is derived from this connection.

In Section \ref{sec-alg}, algorithms are given to check whether a
$\Zx$-lattice is $\Z$-, $x$-, M-, P-, or $\Zx$-saturated, or
equivalently, whether a Laurent binomial difference ideal  is prime,
reflexive, well-mixed, perfect, or toric. If the answer is negative, we can also
compute the $\Z$-, $x$-, M-, {\rm{P}}-, or $\Zx$-saturation of $L$.
Based on the above algorithms, we give algorithms to compute the
reflexive, well-mixed, and perfect closures of a Laurent binomial difference
ideal and an algorithm to decompose a perfect binomial difference
ideal as the intersection of reflexive prime difference ideals. More precisely, we have
\begin{theorem}\label{th-m4}
Let $\F$ be an algebraically closed and inversive difference field
and $F=\{f_1,\ldots,$ $f_s\}$ $\subset\F\{\Y\}$ a set of difference
binomials. Then we can compute difference regular and coherent
ascending chains $\A_1,\ldots,\A_t$ such that
$$\{F\}=\cap_{i=1}^t \sat(\A_i)$$
where for each $i$, $\A_i$ consists of either $y_c\in\Y$ or
difference binomials, and $\sat(\A_i)$ is a reflexive prime
difference ideal. If $t=0$, we mean $\{F\}=[1]$.
\end{theorem}
The above result is stronger than the general decomposition algorithm  
given in \cite{gao-dcs} in that a finitely generated perfect difference ideal 
is decomposed as the intersection of reflexive prime difference ideals. 
For general difference polynomials, this is still an open problem, 
because we do not know how to check whether $\sat(\A)$ is a reflexive 
prime difference ideal for a difference ascending chain $\A$.

Finally, we make a comparison with differential algebra.
The study of binomial differential ideals is more difficult, 
because the differentiation of a binomial is generally not a binomial anymore.
Differential varieties were defined in \cite{d-sres} and were used to connect the differential Chow form \cite{gao} and differential sparse resultant.
But, contrary to the difference case, the defining ideal for a differential toric variety is generally not binomial and further study of differential toric varieties is not carried out.

\section{Preliminaries about difference algebra}
\label{sec-pre}
In this section, some basic notations and preliminary results about
difference algebra and characteristic set for difference polynomial
systems will be given. For more details about difference algebra,
please refer to \cite{cohn, Hrushovski1, levin, wibmer}. For more
details about characteristic set for difference  polynomial systems,
please refer to \cite{gao-dcs}.

\subsection{Difference polynomial and Laurent difference polynomial}

An ordinary difference field, or simply a $\sigma$-field, is a field
$\F$ with a third unitary operation $\sigma$ satisfying that for any
$a, b\in\F$, $\sigma(a+b)=\sigma(a)+\sigma(b)$,
$\sigma(ab)=\sigma(a)\sigma(b)$, and $\sigma(a)=0$ if and only if
$a=0$.
We call $\sigma$ the {\em transforming operator} of $\F$. If
$a\in\F$,  $\sigma(a)$ is called the transform of $a$ and is denoted
by $a^{(1)}$. And for $n\in\Z_{>0}$,
$\sigma^n(a)=\sigma^{n-1}(\sigma(a))$ is called the $n$-th transform
of $a$ and denoted by $a^{(n)}$,  with the usual assumption
$a^{(0)}=a$.
If $\sigma^{-1}(a)$ is defined for each $a\in\F$,  $\F$ is called
{\em inversive}. Every difference field has an inversive closure
\cite{cohn}. A typical example of inversive difference field is
$\Q(\lambda)$ with $\sigma(f(\lambda))=f(\lambda+1)$.

In this paper, $\F$ is assumed to be inversive and of characteristic
zero. Furthermore, we use $\sigma$- as the abbreviation for
difference or transformally.

We introduce the following useful notation.
Let $x$ be an algebraic indeterminate and $p=\sum_{i=0}^s c_i x^i
\in\Z[x]$. For $a$ in any $\sigma$-over field of $\F$, denote
 $$ a^p = \prod_{i=0}^s (\sigma^i a)^{c_i}.$$
For instance,  $a^{x^2-1} = a^{(2)}/a$. It is easy to check that for
$p, q\in\Z[x]$,  we have
 $$a^{p+q}=a^{p} a^{q}, \quad a^{pq}= (a^{p})^{q}.$$
By $a^{[n]}$ we mean the set $\{a, a^{(1)}, \ldots, a^{(n)}\}$. If
$S$ is a set of elements, we denote $S^{[n]}=\cup_{a\in S} a^{[n]}$.

Let $S$ be a subset of a  $\sigma$-field $\mathcal{G}$ which
contains $\mathcal {F}$.   We will  denote respectively by $\mathcal
{F}[S]$,  $\mathcal {F}(S)$,  $\mathcal {F}\{S\}$,  and $\mathcal
{F}\langle S\rangle$  the smallest subring,  the smallest subfield,
the smallest $\sigma$-subring,  and the smallest $\sigma$-subfield
of $\mathcal{G}$ containing $\mathcal {F}$ and $S$.  If we denote
$\Theta(S)=\{\sigma^ka|k\geq0, a\in S\}$,  then we have $\mathcal
{F}\{S\}=\mathcal    {F}[\Theta(S)]$ and $\mathcal {F}\langle
S\rangle=\mathcal    {F}(\Theta(S))$.

A subset $\mathcal {S}$ of a  $\sigma$-extension field $\mathcal
{G}$ of $\mathcal {F}$ is said to be {\em $\sigma$-dependent} over
$\mathcal {F}$ if the set $\{\sigma^{k}a\big|a\in\mathcal {S},
k\geq0\}$ is algebraically dependent over $\mathcal {F}$,  and is
said to be {\em $\sigma$-independent} over $\mathcal {F}$,  or to be
a family of {\em $\sigma$-indeterminates} over $\mathcal {F}$ in the
contrary case.
In the case $\mathcal{S}$ consists of one element $\alpha$,  we say
that $\alpha$ is {\em $\sigma$-algebraic} or {\em
$\sigma$-transcendental} over $\mathcal {F}$,  respectively. The
maximal subset $\Omega$ of $\mathcal {G}$ which are
$\sigma$-independent over $\mathcal {F}$ is said to be a
$\sigma$-transcendence basis of $\mathcal {G}$ over $\mathcal {F}$.
We use $\dtrdeg \, \mathcal {G}/\mathcal {F}$  to denote the {\em
$\sigma$-transcendence degree} of $\mathcal {G}$ over $\mathcal
{F}$,  which is the cardinal number of $\Omega$.
%

Now suppose $\Y=\{y_{1},   \ldots,  y_{n}\}$ is a set of
$\sigma$-indeterminates over $\F$.   The elements of $\mathcal
{F}\{\Y\}=\mathcal {F}[y_j^{(k)}:j=1, \ldots, n;k\in \N]$ are called
{\em $\sigma$-polynomials} over $\F$ in $\Y$,  and $\mathcal
{F}\{\Y\}$ itself is called the {\em $\sigma$-polynomial ring } over
$\F$ in $\Y$. A {\em $\sigma$-polynomial ideal}, or simply a
$\sigma$-ideal, $\mathcal {I}$ in $\mathcal {F}\{\Y\}$ is an
ordinary algebraic ideal which is closed under transforming, i.e.
$\sigma(\mathcal {I})\subset\mathcal {I}$. If $\mathcal{I}$ also has
the property that $a^{(1)}\in\mathcal{I}$ implies that
$a\in\mathcal{I}$,  it is called a {\em reflexive $\sigma$-ideal}.
And a prime $\sigma$-ideal is a $\sigma$-ideal which is prime as an
ordinary algebraic polynomial ideal. For convenience,  a prime
$\sigma$-ideal is assumed not to be the unit ideal in this paper.
A $\sigma$-ideal $\I$ is called {\em well-mixed} if $fg\in\I$
implies $fg^x\in\I$ for $f,g\in\F\{\Y\}$.
A $\sigma$-ideal $\I$ is called {\em perfect} if for any
$a\in\N[x]\setminus\{0\}$ and $p\in\F\{\Y\}$, $p^a\in\I$ implies
$p\in\I$.
If $S$ is a subset of  $\F\{\Y\}$,  we use $(S)$, $[S]$, $\langle
S\rangle$, and $\{S\}$ to denote the algebraic ideal, the
$\sigma$-ideal, the well-mixed $\sigma$-ideal, and the perfect
$\sigma$-ideal generated by $S$.

An $n$-tuple over $\F$ is an $n$-tuple of the form $\eta=(\eta_1,
\ldots, \eta_n)$ where the $\eta_i$ are selected from a
$\sigma$-overfield of $\F$. For a $\sigma$-polynomial
$f\in\F\{\Y\}$, $\eta$ is called a $\sigma$-zero of $f$ if when
substituting $y_i^{(j)}$ by $\eta_i^{(j)}$ in $f$, the result is
$0$.

An $n$-tuple $\eta$ is called a {\em generic zero} of a
$\sigma$-ideal $\CI\subset\F\{\Y\}$ if for any  $P\in\F\{\Y\}$ we
have $P(\eta)=0 \Leftrightarrow P\in\CI$. It is well known that a
$\sigma$-ideal possesses a generic zero if and only if it is a
reflexive prime $\sigma$-ideal other than the unit ideal
\cite[p.77]{cohn}.
Let $\CI$ be a reflexive prime $\sigma$-ideal and $\eta$ a generic
zero of $\CI.$ The {\em dimension} of $\CI$ is defined to be
$\dtrdeg\F\langle\eta\rangle/\F$.


We now define the concept of $\sigma$-variety. Let $\F$ be an
inversive $\sigma$-field,  following \cite{wibmer}, we denote the
category of $\sigma$-field extensions of $\F$ by $\mathscr{E}_\F$,
the category of $\E^{n}$ by $\mathscr{E}_{\F}^{n}$ where $\E \in
\mathscr{E}_{\F}$. Let $(\AH)^{n}$ be the functor from
$\mathscr{E}_{\F}$ to $\mathscr{E}_{\F}^{n}$ satisfying
$(\AH)^{n}(\E)=(\E)^{n}$ where $\E \in \mathscr{E}_{\F}$.
A {\em $\sigma$-variety} over $\F$ is a functor $\V$ from
$\mathscr{E}_\F$ to the category of sets with the form $\V(P)$ for
$P\subset\F\{\Y\}$ satisfying $\V_{\E}(P)= \{\eta\in\E^n\,|\,
\forall p\in P, p(\eta)=0\}$.
It is well known that $\sigma$-varieties are in a one to one
correspondence with perfect $\sigma$-ideals.

\vskip10pt
For $\f=(f_1, \ldots, f_n)^\tau\in \Z[x]^{n}$,  we define $\Y^\f =
\prod_{i=1}^n y_i^{f_i}$. $\Y^\f$ is called a {\em Laurent
$\sigma$-monomial} in $\Y$ and $\f$ is called its {\em support}.
A vector $\f=(f_1, \ldots, f_n)^\tau\in\Z[x]^n$ is said to be {\em
normal} if  the leading coefficient of $f_s$ is positive,  where $s$
is the largest subscript such that $f_s\ne0$.

A {\em Laurent $\sigma$-polynomial} over $\F$ in $\Y$ is an
$\F$-linear combinations of Laurent $\sigma$-monomials in $\Y$.
Clearly, the set of all Laurent $\sigma$-polynomials form a
commutative $\sigma$-ring under the obvious sum, product, and the
usual transforming operator $\sigma$, where all Laurent
$\sigma$-monomials are invertible. 
We denote the $\sigma$-ring of Laurent $\sigma$-polynomials with
coefficients in $\mathcal {F}$ by $\F\{\Y^{\pm}\}$.


Let $p$ be a Laurent $\sigma$-polynomial in $\F\{\Y^{\pm}\}$. An
$n$-tuple $(a_1,$ $\ldots,a_n)$ over $\F$ with each $a_i\neq0$ is
called a {\em nonzero $\sigma$-solution} of $p$  if
$p(a_1,\ldots,a_n)=0.$
%
%
The concept of generic point for a Laurent $\sigma$-ideal can be
defined similarly and it can be proved that a proper Laurent
$\sigma$-ideal is reflexive and prime if and only if it has a
generic point.

\subsection{Characteristic set for a difference polynomial system}
\label{sec-cs}
%

Let $f$ be a $\sigma$-polynomial in $\F\{\Y\}$.  The order of $f$
w.r.t. $y_i$ is defined to be the greatest number $k$ such that
$y_{i}^{(k)}$ appears effectively in $f$,  denoted by $\ord(f,
y_{i})$. If $y_{i}$ does not appear in $f$,  then we set $\ord(f,
y_{i})=-\infty$. The {\em order} of $f$ is defined to be $\max_{i}\,
\ord(f, y_{i})$,  that is, $\ord(f)=\max_{i}\, \ord(f, y_{i})$.

The {\em elimination ranking} $\mathscr{R}$ on $\Theta
(\Y)=\{\sigma^ky_i|1\leq i\leq n, k\geq0\}$ is used in this paper:
$\sigma^{k} y_{i}>\sigma^{l} y_{j}$ if and only if $i>j$ or $i=j$
and $k>l$, which is a total order over $\Theta (\Y)$.     By
convention, $1<\theta y_{j}$ for all $\theta y_{j}\in \Theta (\Y)$.

Let $f$ be a $\sigma$-polynomial in  $\mathcal {F}\{\Y\}$.  The
greatest $y_j^{(k)}$ w.r.t.  $\mathscr{R}$ which  appears
effectively in $f$ is called the {\em leader} of $f$,     denoted by
$\lead(f)$ and correspondingly $y_j$ is called the {\em leading
variable }of $f$, denoted by $\lv(f)=y_j$.
%
%
The leading coefficient of $f$  as a univariate polynomial in
$\lead(f)$ is called the {\em initial} of $f$ and is denoted by
$\init_{f}$.

    Let $p$ and $q$ be two $\sigma$-polynomials in $\F\{\Y\}$.
    $q$ is said to be of higher rank than $p$ if

    1) $\lead(q)>\lead(p)$,  or

    2) $\lead(q)=\lead(p)=y_j^{(k)}$ and $\deg(q, y_j^{(k)})>\deg(p,  y_j^{(k)})$.

Suppose $\lead(p)=y_j^{(k)}$.  $q$ is said to be {\em reduced}
w.r.t. $p$ if  $\deg(q, y_j^{(k+l)})<\deg(p, y_j^{(k)})$  for all
$l\in\N$.

 A finite sequence of nonzero $\sigma$-polynomials $\mathcal
{A}=A_1, \ldots, A_m$  is said to be a
    {\em difference ascending chain}, or simply a {\em $\sigma$-chain}, if

    1) $m=1$ and $A_1\neq0$ or

    2) $m>1$,  $A_j>A_i$ and $A_j$ is reduced
    w.r.t. $A_i$ for $1\leq i<j\leq m$.

A $\sigma$-chain $\mathcal{A}$ can be written as the following form
\begin{equation}\label{eq-asc}
\mathcal{A}=\left\{\begin{array}{l}A_{11}, \ldots, A_{1k_1}\\
\ldots \\ A_{p1}, \ldots, A_{pk_p} \end{array} \right.
\end{equation}
where $\lv(A_{ij})=y_{c_i}$ for $j=1, \ldots, k_i$ and $\ord(A_{ij},
y_{c_i})<\ord(A_{il}, y_{c_i})$ for $j<l$.
%

\begin{example}\label{ex-L11}
The following are three $\sigma$-chains
\[
 \begin{array}{llllll}
 \A_1& =& y_1^{x}-1, &y_1^2y_2^2-1, &y_2^{x}-1& \\
 \A_2& =& y_1^2-1,   &y_1^{x}-y_1, &y_2^2-1, &y_2^{x}-y_2\\
 \A_3& =& y_2^2-y_1^{x},& y_3^{2}-y_1, &y_3^{x}-y_2&\\
\end{array}\]
\end{example}

Let $\mathcal {A}=A_{1},A_{2},\ldots,A_{t}$ be a $\sigma$-chain with
$\init_{i}$ as the initial of $A_{i}$, and $f$ any
$\sigma$-polynomial.
   Then there exists an
   algorithm, which reduces
   $f$ w.r.t. $\mathcal {A}$ to a  polynomial $r$ that is
   reduced w.r.t. $\mathcal {A}$ and satisfies the relation
   \begin{equation}\label{eq-prem}
   \prod_{i=1}^t \init_{i} ^{e_{i}} \cdot f \equiv
   r, \mod \, [\mathcal {A}],\end{equation}
   where the $e_{i}\in\N[x]$.
The $\sigma$-polynomial $r=\prem(f,\A)$ is called the {\em
$\sigma$-remainder} of $f$ w.r.t. $\A$~\cite{gao-dcs}.

A $\sigma$-chain $\mathcal {C}$ contained in a $\sigma$-polynomial
set $\mathcal {S}$ is said to be a {\em characteristic set} of
$\mathcal {S}$, if  $\mathcal {S}$ does not contain any nonzero
element reduced w.r.t. $\mathcal {C}$. A characteristic set
$\mathcal{C}$ of a $\sigma$-ideal $\mathcal {J}$ reduces to zero all
elements of $\mathcal {J}$.

Let $\A:A_1,\ldots,A_t$ be a $\sigma$-chain, $I_i =\init(A_i)$,
$y_{l_i}^{(o_i)} = \lead(A_i)$.
$\A$ is called {\em regular} if for any $j\in \N$, $I_i^{x^j}$ is
invertible w.r.t $\A$ \cite{gao-dcs} in the sense that
$[A_1,\ldots,A_{i-1},I_i^{x^j}]$ contains a nonzero
$\sigma$-polynomial involving no $y_{l_i}^{(o_i+k)},k=0,1,\ldots$.
%
%
%
To introduce the concept of coherent $\sigma$-chain, we need to
define the {\em $\Delta$-polynomial} first. If $A_i$ and $A_j$ have
distinct leading variables, we define $\Delta(A_i,A_j)=0$. If $A_i$
and $A_j$ ($i<j$) have the same leading variable $y_l$, then
$o_i=\ord(A_i,y_l) < o_j=\ord(A_j,y_l)$.
Define
 \begin{equation}\label{eq-delta}
 \Delta(A_i,A_j) =\prem((A_i)^{x^{o_j-o_i}},A_j).\end{equation}
 Then $\A$ is called {\em
coherent} if $\prem(\Delta(A_i,A_j),\A)=0$ for all $i< j$
\cite{gao-dcs}.

Let $\mathcal {A}$ be a $\sigma$-chain. Denote $\mathbb{I}_{\mathcal
{A}}$ to be  the minimal multiplicative set containing the initials
of elements of $\mathcal{A}$ and their transforms.    The {\em
saturation ideal} of $\A$ is defined to be
 $$\sat(\A)=[\mathcal  {A}]:\mathbb{I}_{\mathcal
 {A}} = \{p\in\F\{\Y\}: \exists h\in \mathbb{I}_{\mathcal {A}},  \, {\text s. t.}\,  hp\in[A]\}.$$
%
%
 The following result is needed in this paper.
\begin{theorem}\cite[Theorem 3.3]{gao-dcs}\label{th-rp}
A $\sigma$-chain $\A$ is a characteristic set of $\sat(A)$ if and
only if $\A$ is regular and coherent.
\end{theorem}
%
%
%
%

\section{$\Zx$-lattice}
\label{sec-zxm}
In this section, we prove basic properties of $\Z[x]$-lattices,
which will play the role of lattices in the study of binomial ideals
and toric varieties.

\subsection{Gr\"{o}bner basis and generalized Hermite normal form}
\label{subsec-Gro}

For brevity, a $\Z[x]$-module in $\Z[x]^n$ is called a {\em $\Z[x]$-lattice}.
Since $\Z[x]$ is a Noetherian ring, we have
\begin{lemma}
Any $\Z[x]$-lattice is finitely generated.
\end{lemma}
As a consequence, any $\Zx$-lattice $L$ has a finite set of
generators $\fb=\{\f_1,\ldots,\f_s\}\subset\Z[x]^n$:
 $$L = \span_{\Zx}\{\f_1,\ldots,\f_s\} = (\f_1,\ldots,\f_s).$$
A {\em matrix representation} of $\f$ or $L$ is
 $$M = [\f_1,\ldots,\f_s]_{n\times s},$$
with $\f_i$ to be the $i$-th column of $M$. We also denote $L=(M)$.
\begin{definition}
The {\em rank} of a $\Zx$-lattice $L$ is defined to be the rank of
any matrix representation of  $L$.
\end{definition}
The concept of rank is clearly well defined.

A standard form to represent the submodules in $\Z[x]^n$ is the
Gr\"{o}bner basis. We list some basic concepts and properties of
Gr\"{o}bner basis of modules. For details, please refer to
\cite{cox-1998}.

Denote ${\beps}_i$ to be the $i$-th standard basis vector
$(0,\ldots,0,1,0,\ldots,0)^\tau\in\Z[x]^n$, where $1$ lies in the
$i$-th row of ${\beps}_i$. A {\em monomial} ${\bf m}$ in $\Z[x]^n$
is an element of the form $ax^k{\beps}_i\in\Z[x]^n$, where $a\in\Z$
and $k\in\N$.
The following {\em monomial order} $>$ of $\Z[x]^n$ will be used in
this paper: $ax^\alpha{\beps}_i > bx^\beta{\beps}_j$ if $i
> j$, or $i=j$ and $\alpha > \beta$, or $i=j$, $\alpha = \beta$, and
$|a|>|b|$.

 With the above order, any $\f\in\Z[x]^n$ can be written in a unique
way as a  linear combination of monomials,
$$ \f = \sum_{i=1}^s \f_i,$$
where ${\f_i}\ne 0$ and $\f_1 > \f_2>\cdots > \f_s$. The {\em
leading term} of $\f$  is defined to be $
\LT(\f) = \f_1.$
For any $\G\subset \Zxn$, we denote by $\LT(\G)$ the  set of leading
terms of $\G$.

The order $>$ can be extended to elements of $\Zxn$ as follows: for
$\f,\g\in\Zxn$, $\f <\g$ if and only if $\LT(\f) < \LT(\g)$.

Let $\G\subset\Z[x]^n$ and  $\f\in\Z[x]^n$. We say that $\f$ is {\em
G-reduced} with respect to $\G$ if any monomial of $\f$ is not a
multiple of $\LT(\g)$ by an element in $\Zx$ for any $\g\in\G$.

\begin{definition}
A finite set $\fb = \{\f_1,\ldots,\f_s \}\subset\Z[x]^n$ is called a
{\em Gr\"{o}bner basis} for the $\Zx$-lattice $L$ generated by $\fb$
if $(\LT(L)) = (\LT(\fb))$.
A Gr\"{o}bner basis $\fb$ is called {\em reduced} if for any
$\f\in\fb$, $\f$ is G-reduced with respect to $\fb\setminus\{\f\}$.
\end{definition}

Let $\fb$ be a Gr\"{o}bner basis. Then any $\f\in\Zxn$ can be
reduced to a unique normal form by $\fb$, denoted by
$\grem(\f,\fb)$, which is G-reduced with respect to $\fb$.

\begin{definition}\label{def-sv}
Let $\f,\g\in\Z[x]^n$,  $\LT(\f)=ax^k{\beps}_i$,
$\LT(\g)=bx^s{\beps}_j$, and $t = \max(k,s)$.
Then the {\em S-vector} of $\f$ and $\g$ is defined as follows: if
$i\ne j$ then $S(\f,\g)=0$; otherwise assume $|a|\ge |b|$ and let
 $$S(\f,\g) = x^{t-k}\f - cx^{t-s}\g$$
where $c=s\lfloor\frac{|a|}{|b|}\rfloor$ and $s$ is the sign of
$ab$.
\end{definition}
In the later case, let $a = bc+r$. Then $0\le |r| < |b|$ and
$\LT(S(\f,\g)) =rx^t{\beps}_i$.
Using the above notations, we have the following basic property for
Gr\"{o}bner basis \cite{Kandri}.
\begin{theorem}[Buchberger's Criterion]\label{th-buch}
For a set $\fb = \{\f_1,\ldots,\f_s \}\subset\Z[x]^n$, the following
statements are equivalent.
\begin{description}
\item[1)] $\fb$ is a Gr\"{o}bner basis.

\item[2)] $\grem(S(\f_i,\f_j),G)=0$  for all $i,j$.

\item[3)] $\f\in(\fb)$ if and only if $\grem(\f,\fb)=0$.
\end{description}

\end{theorem}

We will study the structure of a Gr\"obner basis for a $\Zx$-lattice
by introducing the concept of generalized Hermite normal form.
First, we consider the case of $n=1$.

\begin{lemma}\label{lm-21}
Let  $B = \{b_1,\ldots,b_k\}$ be a reduced Gr\"{o}bner basis of a
$\Zx$-module in $\Z[x]$, $b_1<\cdots<b_k$, and
$\LT(b_i)=c_ix^{d_i}\in\N[x]$. Then

1) $0\le d_1 < d_2 < \cdots<d_k$.

2) $c_k | \cdots | c_2 | c_1 $ and $c_i \ne c_{i+1}$ for $1\le i\le
k-1$.

3) $\frac{c_i}{c_k} | b_i$ for $1\le i< k$. Moreover if
$\widetilde{b}_1$ is the primitive part of $b_1$, then
$\widetilde{b}_1 | b_i$ for $1< i\le k$.
%

4) The S-polynomial $S(b_i,b_j)$ can be reduced to zero by $B$ for
any $i,j$.
\end{lemma}
\proof  1) and 4) are consequences of Theorem~\ref{th-buch}. To
prove 2), assume that there exists an $l$ such that $c_{l-1} |
\cdots | c_2 | c_1 $ but $c_l \not| c_{l-1}$. Let $r =
\gcd(c_l,c_{l-1}) = p_1 c_l+p_2c_{l-1}$, where $p_1,p_2\in\Z$.
Then $|r| < |c_{l-1}|$ and $|r| < |c_{l}|$. Since $c_{l-1} | \cdots
| c_2 | c_1 $, we have $|r| < |c_{i}|, i=1,\ldots,l$. Let $g =
p_1b_l + p_2x^{d_l-d_{l-1}}b_{l-1}$. Then $\LT(g) = rx^{d_l}$
which is reduced w.r.t. $B$ and $g\in (B)$, contradicting to the
definition of Gr\"{o}bner bases.

We prove 3) by induction on $k$. When $k=2$, let $b_1 = c_1x^{d_1} +
s_{11}x^{d_1-1} + \cdots + s_{1d_1}$ and $b_2 = c_2x^{d_2} +
s_{21}x^{d_2-1} + \cdots + s_{2d_2}$. Then, $c_2 | c_1$ and $d_1 <
d_2$. Let $c_1 = c_2 t$, we need to show $t | b_1$. Since the
S-polynomial $S(b_1,b_2) = tb_2-x^{d_2-d_1}b_1$ can be reduced to
zero by $b_1$, we have $tb_2-x^{d_2-d_1}b_1 = u(x)b_1$, where
$u(x)\in\Zx$ and $\deg(u(x))<d_2-d_1$. Then, $tb_2 = (x^{d_2-d_1} +
u(x))b_1$, and $t | b_1$ follows since $x^{d_2-d_1} + u(x)$ is a
primitive polynomial in $\Zx$. The claim is true. Assume that for
$k=l-1$, the claim is true, then
$\widetilde{b}_1 | b_i$ for $1\le i \le l-1$. We will prove the
claim for $k=l$. Since $S(b_{1},b_l) =
\frac{c_{1}}{c_l}b_l-x^{d_l-d_{1}}b_{1}$ can be reduced to zero by
$B$. We have $\frac{c_{1}}{c_l}b_l-x^{d_l-d_{1}}b_{1} =
\sum_{i=1}^{l-1} f_ib_i$ with $f_i\in\Zx$ and $\deg(f_ib_i) \le
d_l-1$. Then, $\frac{c_{1}}{c_l}b_l = x^{d_l-d_{1}}b_{1} +
\sum_{i=1}^{l-1} f_ib_i$. By induction, $\widetilde{b}_1$ is a
factor of the right hand side of the above equation. Thus
$\widetilde{b}_1 | b_l$. Let $b_i = s_ib_{1}'$ for $1\le i \le l$,
we have $\frac{c_{1}}{c_l}s_{l} = x^{d_l-d_1}s_1+\sum_{i=1}^{l-1}
f_is_i$ where $\deg(s_i) = d_i-d_1$ and $s_1\in\Z$. Since
$\deg(f_is_i)\le d_l-d_1-1$, we have $\frac{c_1}{c_l} | s_1$ and
$\frac{c_1}{c_l} | b_1$.
For any $1\le i < j < l$, assume $\frac{c_i}{c_l} | b_i$. We will
show that $\frac{c_j}{c_l} | b_j$. Since $S(b_{j-1},b_{j}) =
\frac{c_{j-1}}{c_j}b_{j} -x^{d_j-d_{j-1}}b_{j-1} = \sum_{i=1}^{j-1}
f_ib_i$, we have $\frac{c_{j-1}}{c_l}$ is a factor of the right hand
side of the above equation, for $c_{j-1} | c_{j-2} | \cdots |
c_1$. Then, $\frac{c_{j-1}}{c_l} | \frac{c_{j-1}}{c_j}b_{j}$ and $
\frac{c_{j}}{c_l} | b_{j}$. The claim is proved. \qedd

\begin{example}\label{ex-1}
Here are three Gr\"{o}bner bases in $\Z[x]$\hbox{\rm:} $\{ 2, x \}$,
$ \{ 12, 6x+6, 3x^2+3x, x^3+x^2 \}$, $ \{ 9x+3, 3x^2+4x+1 \}$.
\end{example}

To give the structure of a reduced Gr\"{o}bner basis similar to that
in Example \ref{ex-1}, we introduce the concept of  generalized
Hermite normal form. Let
{\small
\begin{equation}\label{ghf}\C=\left[\begin{array}{lllllllllll}
c_{1,1}   & \ldots & c_{1,l_1}     &c_{1,l_1+1}   &\ldots       &\ldots       &\ldots &\ldots &\ldots            &\ldots  &\ldots   \\
\ldots    & \ldots & \ldots        & \ldots       &\ldots       &\ldots       &\ldots &\ldots &\ldots            &\ldots  &\ldots \\
c_{r_1, 1}& \ldots & c_{r_1, l_1}  &c_{r_1,l_1+1} & \ldots      &\ldots       &\ldots &\ldots & \ldots           &\ldots  &\ldots \\
0         & \ldots & 0             &c_{r_1+1,1}   &    \ldots   &c_{r_1+1,l_2}&\ldots &\ldots &  \ldots          &\ldots  &\ldots  \\
\ldots    & \ldots & \ldots        & \ldots       &\ldots       &\ldots       &\ldots &\ldots &\ldots            &\ldots  &\ldots  \\
0         & \ldots & 0             & c_{r_2,1}    &    \ldots   & c_{r_2,l_2} &\ldots &\ldots & \ldots           &\ldots  &\ldots   \\
\ldots    & \ldots & \ldots        &\ldots        &    \ldots   &\ldots       &\ldots &\ldots & \ldots           &\ldots  &\ldots   \\
0         & \ldots & 0             &0             &    \ldots   &0            &\ldots &0      &c_{r_{t-1}+1,1}   &\ldots  & c_{r_{t-1}+1,l_t} \\
\ldots    & \ldots & \ldots        &\ldots        &    \ldots   &\ldots       &\ldots &\ldots &\ldots            &\ldots  & \ldots \\
0         & \ldots & 0             &0             &    \ldots   &0            &\ldots &0      &c_{r_t,1}         &\ldots  & c_{r_t,l_t} \\
\end{array}\right]_{m\times s}\end{equation}}
whose elements are in $\Zx$.
It is clear that  $m = r_t$ and $s=\sum_{i=1}^t l_i$. We denote by
$\c_k=\c_{r_i,j}$ to be the $k$-th column of the matrix $\C$, where
$k = l_1+l_2+\ldots+l_{i-1}+j, 1\le j\le l_{i}$.
Assume
 \begin{equation}\label{eq-cij}
 c_{i,j}=c_{i,j,0}x^{d_{ij}}+\cdots+c_{i,j,d_{ij}}.\end{equation}
Then the leading monomial of  $\c_{r_i,j}$ is
$c_{r_i,j,0}x^{d_{r_i,j}}{{\beps}_{r_i}}$.

\begin{definition}\label{def-ghf}
The matrix $\C$ in~\bref{ghf} is called a {\em generalized Hermite
normal form} if it satisfies the following conditions:
\begin{description}
\item[1)] $0\le d_{r_i,1} < d_{r_i,2} < \cdots < d_{r_i,l_i}$ for any $i$.

\item[2)] $ c_{r_i,l_i,0} |\cdots | c_{r_i,2,0}  | c_{r_i,1,0}$.

\item[3)]
$S(\c_{r_i,j_1},\c_{r_i,j_2})=
x^{d_{r_i,j_2}-d_{r_i,j_1}}\c_{r_i,j_1}-\frac{c_{r_i,j_1,0}}{c_{r_i,j_2,0}}\c_{r_i,j_2}$
can be reduced to zero by the column vectors of the matrix for any
$1\le i\le t, 1\le j_1<j_2\le l_i$.

\item[4)]
$\c_{r_i,j}$ is G-reduced w.r.t. the column vectors of the matrix
other than $\c_{r_i,j}$, for any $1\le i\le t, 1\le j\le l_i$.
\end{description}
\end{definition}

It is clear that $\{c_{r_i,1},\ldots,c_{r_i,l_i}\}$ is a reduced
Gr\"obner basis in $\Zx$. Then, as a consequence of Theorem
\ref{th-buch} and Lemma \ref{lm-21}, we have
\begin{theorem}\label{th-gbh}
$\fb = \{\f_1,\ldots,\f_s \}\subset\Z[x]^n$ is a reduced Gr\"{o}bner
basis such that $\f_1< \f_2< \cdots < \f_s$ if and only if the
matrix representation for $\fb$ is a generalized Hermite normal
form.
\end{theorem}

The following property of $\C$ will be used later.
\begin{lemma}\label{lm-L1}
Let $\C$ be given in \bref{ghf} be a generalized Hermite normal form
and $L=(\C)$. Then $\rk(L) = t$ and $L$ is a free $\Zx$-module if
$l_i=1$ for $i=1,\ldots,t$.
\end{lemma}

\begin{example}\label{ex-L1}
The following matrices are generalized Hermite normal forms
 \[M_1=\left [\begin{array}{llll}
x   & 2 & 0   \\
0   & 2 & x    \\
\end{array}\right],\,
 M_2=\left[\begin{array}{llll}
2   & x-1 & 0 & 0  \\
0   & 0 & 2      & x-1      \\
\end{array}\right],\,
 M_3=\left[\begin{array}{llll}
-x  & -1 & 0  \\
2   & 0 & -1      \\
0         & 2 & x\\
\end{array}\right].\]
whose columns constitute the reduced  Gr\"{o}bner bases of the $\Zx$-lattices.
For instance, $\fb_1 = \{[x,0]^\tau,[2,2]^\tau,[0,x]^\tau\}$ is the
reduced Gr\"{o}bner basis of $(\fb_1)$.
\end{example}

\begin{example}\label{ex-L2}
The number of generators for a $\Zx$-lattice depends on how to
arrange the row elements. For instance, if we move the third row of
$M_3$ in Example \ref{ex-L1} to be the first row, then we have
$\widetilde{M}_3=\left[\begin{array}{llll}
0   & 2 & x\\
-x  & -1 & 0  \\
2   & 0 & -1      \\
\end{array}\right]$.
Then $(M_3)$ has another set of generators $(2,-1,0)^T$ and
$(x,0,-1)^T$. In other words, if $N_3=\left[\begin{array}{lll}
 2 & x\\
 -1 & 0  \\
 0 & -1      \\
\end{array}\right]$,
then $(\widetilde{M}_3)=(N_3)$ which is a free $\Zx$-lattice.
\end{example}

Let $\C$ be defined in \bref{ghf} and $k\in\N$. Introduce the
following notations:
 \begin{eqnarray}
 \C_{-} &=& \cup_{i=1}^t \cup_{k=1}^{l_i-1}    \{\c_{r_i,k}, x\c_{r_i,k}, \ldots,
    x^{\deg(c_{r_i,k+1})-\deg(c_{r_t,k})-1}\c_{r_i,k}\},\nonumber\\
 \C^{+} &=& \cup_{i=1}^t \cup_{k=0}^{\infty}  \{x^k \c_{r_i,l_i}\}.\label{eq-inf} \\
 \C_{\infty} &=& \C_{-}\cup   \C^{+} \nonumber
 \end{eqnarray}
$\C_{\infty}$ is called the {\em extension} of $\C$.

\begin{example}\label{ex-ext}
Let $\C =\left [\begin{array}{lllll}
6  & 3x              & 0             &3       & 2x         \\
0   & 0            &    6           & 3x       & x^3+x      \\
\end{array}\right].$
Then $\C_{-} =\left [\begin{array}{lllll}
6        & 0          &3  &3x              \\
0        & 6           & 3x  &3x^2    \\
\end{array}\right]$ and
\[ \C_{\infty} = \left[ \begin{array}{lllllllllll}
6   & 3x    & 3x^2    &  3x^3 &\cdots      & 0      &3      &3x           & 2x   & 2x^2 & \cdots  \\
0    & 0    & 0       &   0   &\cdots      &  6    & 3x    &3x^2        & x^3+x  & x^4+x^2 &  \cdots  \\
\end{array} \right]. \]
\end{example}

We need the following properties about the extension of $\C$. By
saying the infinite set $\C_{\infty}$ is linear dependent over $\Z$,
we mean any finite subset of $\C_{\infty}$ is linear dependent over
$\Z$. Otherwise, $\C_{\infty}$ is said to be linear independent.
\begin{lemma}\label{lm-ext1} Let $\C$ be a generalized Hermite normal form.
The columns of $\C_{\infty}$ are linear independent over $\Z$.
\end{lemma}
\proof Suppose $\C$ is given in \bref{ghf}. The leading term of
$\c\in\C_{\infty}$ is $\LT(\c) = ax^l{\beps}_{r_i}$ for
$i=1,\ldots,t$ and $l\in\N$. Furthermore, for two different $\c_1$
and $\c_2$ in $\C_S$ such that $\LT(\c_1) = ax^{l_1}{\beps}_{r_i}$
and $\LT(\c_2) = bx^{l_2}{\beps}_{r_i}$, we have $l_1\ne l_2$. Then
$\LT(\C_{\infty})=\{a_{il_i}x^{l_i}{\beps}_{r_i}\,|\, i=1,\ldots,t;
l_i= d_{i1},d_{i1}+1,\ldots;a_{il_i}\in\Z\}$ are linear independent
over $\Z$, where $d_{i1}$ is from \bref{eq-cij}. Then $\C_{\infty}$
are also linear independent over $\Z$.\qedd
%
%

\begin{lemma}\label{lm-ext2}
Let $\C$ be a generalized Hermite normal form. Then any $\g\in (\C)$
can be written uniquely as a linear combination of finitely many
elements of $\C_\infty$ over $\Z$.
\end{lemma}
\proof $\g\in (\C)$ can be written as a linear combination of
elements of $\C_\infty$ over $\Z$ by the procedure to compute
$\grem(\g,\C)=0$ \cite{cox-1998}. The uniqueness is a consequence of
Lemma~\ref{lm-ext1}.\qedd

Note that syzygies among elements of $\C$ are not linear
combinations over $\Z$. For instance, let
 $\f_1 = [2,0]^\tau$ and  $\f_2 = [x-1,0]^\tau$.
Then we have $(x-1) \f_1 - 2 \f_2={\bf{0}}$ which is not a linear
representation of ${\bf{0}}$ over $\Z$.

\subsection{Kernel of a matrix in $\Z[x]^{n\times s}$}
Let $F = [\f_1,\ldots,\f_s]_{n\times s}\in\Zx^{n\times s}$, where
$\f_i\in\Zxn$. The {\em kernel} of $F$ is
 $$\ker(F) = \{ X\in\Z[x]^s | FX = {\bf 0} \}.$$
It is clear that $\ker(F)$ is a $\Zx$-lattice in $\Zx^s$, which has
the following important property.
\begin{lemma}\label{lm-free}
Let $F = [\f_1,\ldots,\f_s]_{n\times s}\in\Zx^{n\times s}$. Then
$\ker(F)$ is a free  $\Zx$-module.
\end{lemma}
\proof Without loss of generality, we assume $\fb=\{
\f_1,\ldots,\f_s\}$ is a the generalized Hermite normal form which
means $\fb$ is a Gr\"{o}bner basis.
Let $S(\f_i,\f_j)=m_{ij}{\f}_i - m_{ji}{\f}_j $ and
$\grem(S(\f_i,\f_j),\fb) = \sum_k c_k \f_k$ be the normal
representation of in terms of the Gr\"obner basis $\fb$. Then the
{\em syzygy polynomial} $\widetilde{S}(\f_i, \f_j)$
 $$\widetilde{S}(\f_i, \f_j)=m_{ij}{\beps}_i -m_{ji}{\beps}_j - \sum_k c_k {\beps}_k,$$
is an element in $\Z[x]^s$, where ${\beps}_k$ is the $k$-th standard
basis vector of $\Z[x]^{s}$.
Define an order in $\Z[x]^s$ as follows: $ax^\alpha {\beps}_i \prec
bx^\beta {\beps}_j$ if $\LT (ax^\alpha\f_i)
> \LT(bx^\beta \f_j)$ in $\Z[x]^n$.
We will use Schreyer's Theorem on page 224 of~\cite{cox-1998}, which
says that if $\{\f_1,\ldots,\f_s\}$ is a reduce Gr\"obner basis,
then the syzygy polynomials $\widetilde{S}= \{ \widetilde{S}(\f_i,
\f_j)\} $ constitute a Gr\"{o}bner basis of $\ker(M)\subset\Zx^s$
under the newly defined order $\prec$.
By the proof of Schreyer's Theorem, the order $\prec$ is a monomial
order.
Rewrite $\f_i$ as the form~\bref{ghf}, that is, $\{\f_1,\ldots,\f_s
\} = \{\c_{r_i,j} | 1\le i\le t, 1\le j\le l_i \}$.
Let $\s_{ijk} = \widetilde{S}(\c_{r_i,j},\c_{r_i,k})$ with $j<k$. By
(2) of Definition \ref{def-ghf}, $\LT_{\succ}(s_{ijk}) =
\frac{c_{r_i,j,0}}{c_{r_i,k,0}}{\beps}_{\omega_{ik}}$ where
$\omega_{ik} = l_1+l_2+\cdots+l_{i-1}+k$ and
$\frac{c_{r_i,j,0}}{c_{r_i,k,0}}\in\Z$.

We claim that
 \begin{equation}\label{eq-tgb}
 \H=\{\s_{i,k-1,k},1\le i\le t, 2\le k\le l_i\}\end{equation}
is a reduced Gr\"obner basis for $\ker(M)$.
By Schreyer's Theorem, the set of $\s_{ijk}$ forms a Gr\"obner basis
for $\ker(M)$. The claim will be proved if we can show that for
$j=1,\ldots,k-2$, $\LT_{\succ}(\s_{ijk})=h_{ijk}
\LT_{\succ}(\s_{ik-1k})$, where $h_{ijk}\in\Z$. We have
$\LT_{\succ}(s_{ijk}) =
\frac{c_{r_i,j,0}}{c_{r_i,k,0}}{\beps}_{\omega_{ik}}$,
$\LT_{\succ}(s_{ik-1k})$ $ =
\frac{c_{r_i,k-1,0}}{c_{r_i,k,0}}{\beps}_{\omega_{ik}}$. Since
$\frac{c_{r_i,j,0}}{c_{r_i,k,0}} = \frac{c_{r_i,j,0}}{c_{r_i,k-1,0}}
\frac{c_{r_i,k-1,0}}{c_{r_i,k,0}}$, by (2) of Definition
\ref{def-ghf}, $\frac{c_{r_i,j,0}}{c_{r_i,k-1,0}}$ and
$\frac{c_{r_i,k-1,0}}{c_{r_i,k,0}}$ are in $\Z$. The claim is
proved.
%

By Lemma \ref{lm-L1} and the claim, the $\Zx$-lattice $(\H)$ is
free. By Proposition~3.10 of \cite{cox-1998}(page 229), since $(\H)$
is free and $\H$ is set of generators of $\ker(M)$, $\ker(M)$ is
also free. \qedd

\begin{cor}\label{cor-L3}
If $F = [\f_1,\ldots,\f_s]_{n\times s}\in\Zx^{n\times s}$ is a
generalized Hermite normal form, then \bref{eq-tgb} is a reduced
Gr\"{o}bner basis of $\ker(F)$ under the order $\prec$ and
$\rk(\ker(F))  = s - \rk(F)$.
\end{cor}
\proof It suffices to show $\rk(\ker(F))  = s - \rk(F)$. Suppose $F$
has the form \bref{ghf}. Then $\rk(F)=t$. $\H$ in \bref{eq-tgb} is
the reduced Gr\"ober basis for $\ker(F)$. Note that the number of
elements in $\H$ is $s-t$ and these elements are linearly
independent since each of them has the different standard basis
vector. Then $\rk(\ker(F))=s-t=s-\rk(F)$. \qedd

\begin{example}\label{ex-L3}
In Example \ref{ex-L1},
$\f_1=[x,0]^\tau,\f_2=[2,2]^\tau,\f_3=[0,x]^\tau$. Then
 $S(f_2,f_3) = x \f_2 - 2 \f_3 = [2x,0]^\tau = 2\f_1$ and the syzygy
polynomial is  $\s_{212}= x {\beps}_2 - 2 {\beps}_3 - 2 {\beps}_1$.
That is, $\ker(M)$ is generated by $[-2,x,-2]^\tau$.
\end{example}


\begin{remark}
$\ker(F)$ is the syzygy module of $M$, which can be computed by the
Gr\"{o}bner basis method in Chapter 6 of \cite{cox-1998} or
Proposition 3.8 in \cite[p. 227]{cox-1998}.
\end{remark}

\section{Laurent binomial $\sigma$-ideal}
\label{sec-lbi} In this section, we will prove some basic facts
about Laurent binomial $\sigma$-ideals.
As a tool, the characteristic set method for Laurent binomial
$\sigma$-ideals will be presented.

\subsection{Laurent binomial $\sigma$-ideal}
\label{sec-lbi1} In this section, several basic properties of
binomial $\sigma$-ideals will be proved.

By a {\em Laurent $\sigma$-binomial} in  $\Y$,  we mean a
$\sigma$-polynomial with at most two terms,  that is,
$a\Y^{\g}+b\Y^{\h}$ where $a, b\in \F$ and $\g, \h \in\Z[x]^n$.
A Laurent $\sigma$-binomial of the following form is said to be in
{\em normal form}
 $$p=\Y^{\f}- c_\f$$
where $c_\f\in \F^{\ast}=\F\setminus\{0\}$ and $\f\in\Z[x]^n$ is
normal. The vector $\f$ is called the {\em support} of $p$. For
$p=\Y^{\f}- c_\f$,  we denote $\widehat{p}= -c_\f^{-1}\Y^{-\f} p=
\Y^{-\f}- c_\f^{-1}$ which is called the {\em inverse} of $p$. Note
that if $p$ is in a Laurent binomial $\sigma$-ideal $\I$, then
$\widehat{p}$ is also in $\I$.
It is clear that any Laurent $\sigma$-binomial $f$ which is not a
unit can be written uniquely as
$$ f = aM(\Y^{\f}- c_\f)$$
where $a\in \F^{\ast}$, $M$ is a Laurent $\sigma$-monomial, and
$\Y^{\f}- c_\f$ is in normal form. Since $aM$ is a unit in
$\F\{\Y^{\pm}\}$, we can use the normal $\sigma$-binomial $\Y^{\f}-
c_\f$ to represent $f$, and when we say a Laurent $\sigma$-binomial
we always use its normal representation.

A Laurent $\sigma$-ideal is called binomial if it is generated by
Laurent $\sigma$-binomials.

\begin{lemma}\label{lm-bi1}
Let $\Y^{\f_{i}}-c_i,i=1, \ldots,s$ be contained in a Laurent
binomial $\sigma$-ideal $\I$ and $\f=a_1\f_1+ \cdots +a_s\f_s$,
where $a_i\in\Zx$. Then $\Y^{\f} - \prod_{i=1}^s c_i^{a_i}$ is in
$\I$.
\end{lemma}
\proof It is suffices to show that if $p_1=\Y^{{\f}_1}-c_1\in\I$ and
$p_2=\Y^{{\f}_2}-c_2\in\I$, then $\Y^{n{\f}_1}-c_1^n \in\I$ for
$n\in\N$, $\Y^{-{\f}_1}-c_1^{-1}\in\I$,
$\Y^{x{\f}_1}-\sigma(c_1)\in\I$, and
$\Y^{{\f}_1+{\f}_2}-c_1c_2\in\I$, which are indeed true since
$\Y^{n{\f}_1}-c_1^n=(\Y^{{\f}_1})^n-c_1^n$ contains $p_1$ as a
factor, $\Y^{-{\f}_1}-c_1^{-1}=-c_1^{-1}\Y^{-{\f}_1}(\Y^{\f_1}-c_1)
\in\I$, $\Y^{x{\f}_1}-\sigma(c_1)=\sigma(\Y^{{\f}_1}-c_1)\in\I$, and
$\Y^{{\f}_1+{\f}_2}-c_1c_2= \Y^{{\f}_1}(\Y^{{\f}_2}-c_2) +
c_2(\Y^{{\f}_1}-c_1)\in\I$. \qedd

As a direct consequence, we have
\begin{cor}\label{cor-bi1}
Let $\I$ be a Laurent binomial $\sigma$-ideal and
\begin{equation}\label{eq-bi1}
\L(\I):=\{\f\in \Z[x]^{n}\, |\, \exists c_\f\in\F^*\,s.t.\,
\Y^{\f}-c_\f\in \I\}.
\end{equation}
Then $\L(\I)$ is a $\Zx$-lattice, which is called the {\em support
lattice} of $\I$, and a {\em matrix representation} for $\L(\I)$ is
called a matrix representation for $\I$.
\end{cor}

The following lemma shows that a new set of generators  of the
support lattice of a Laurent binomial $\sigma$-ideal leads to a new
set of generators for the $\sigma$-ideal.
\begin{lemma}\label{lm-bi3}
Let $\I=[\Y^{\f_{1}}-c_{1}, \ldots, \Y^{\f_{s}}-c_{s}]$ be a proper
Laurent binomial $\sigma$-ideal and let $\h_1, \ldots, \h_r$ be
another set of generators for the $\Z[x]$-lattice
$(\f_1,\ldots,\f_s)$, and there exist $a_{i,k}\in\Zx$ such that
 $$\h_i = \sum_{k=1}^s a_{i,k} \f_k, i=1,\ldots,r.$$
Then $\I = [\Y^{\h_1}-\prod_{i=1}^s c_i^{a_{1, i}} ,\ldots,
 \Y^{\h_r}-\prod_{i=1}^s c_i^{a_{r, i}}]$.
\end{lemma}
\proof Let $f_t=\Y^{\f_{t}}-c_{t},t=1,\ldots,s$,
$g_l=\Y^{\h_l}-\prod_{i=1}^s c_i^{a_{l, i}},l=1,\ldots,r$, and $\I_1
= [g_1,\ldots,g_r]$.
From Lemma \ref{lm-bi1}, $\I_1\subset \I$.
We now prove $\I\subset \I_1$. Since  $(\h_1, \ldots,
\h_r)=(\f_1,\ldots,\f_s)$, there exist $b_{i,k}\in\Zx$ such that
$\f_i = \sum_{k=1}^r b_{i,k} \h_k.$
Then, $\Y^{\f_i} = \Y^{\sum_{k=1}^r b_{i,k}\h_k}
  = \prod_{k=1}^r (\Y^{\h_k})^{b_{i,k}}$.
Replacing $\Y^{\h_k}$ by $g_k + \prod_{i=1}^s c_i^{a_{l, i}}$, we
have
 $\Y^{\f_i}
  = \prod_{k=1}^s (g_k + \prod_{j=1}^s c_j^{a_{k, j}})^{b_{i,k}}
  = \widetilde{g} + \widetilde{c}$, where
   $\widetilde{c}=\prod_{j=1}^s c_j^{\sum_{k=1}^s b_{i,k}a_{k,j}}\in\F$ and
   $\widetilde{g}\in \I_1\subset\I$.
As a consequence, $\Y^{\f_i} - \widetilde{c} = \widetilde{g} \in\I$.
Then $f_i - (\Y^{\f_i} - \widetilde{c}) = \widetilde{c}- c_i\in\I$.
Since $\I$ is proper, we have $\widetilde{c}=c_i$ and hence $f_i
=\Y^{\f_i} - \widetilde{c} \in\I_1$. The lemma is proved.\qedd

The following lemma can be used to check whether a Laurent binomial
$\sigma$-ideal is proper.
\begin{lemma}\label{lm-bi2}
Let $\I=[\Y^{\f_{1}}-c_{1}, \ldots, \Y^{\f_{s}}-c_{s}]$ be a Laurent
binomial $\sigma$-ideal and let $M$ be the $n\times s$ matrix with
columns $\f_1, \ldots, \f_s$. Furthermore, let $\ker(M)$ be
generated by $\bu_1, \ldots, \bu_t$,  where $\bu_i = (u_{i, 1},
\ldots, u_{is})^\tau\in\Zx^s$.
Then $\I\ne[1]$ if and only if $\prod_{i=1}^s c_i^{u_{l, i}}=1$ for
$l=1, \ldots, t$.
\end{lemma}
\proof ``$\Rightarrow$" Let $f_i = \Y^{\f_i}-c_i$. Suppose
$c=\prod_{i=1}^s c_i^{u_{l, i}}\ne1$ for some $l$. Replacing $c_i$
by $\Y^{\f_{i}} - f_i$ in the above equation and noting that
$\bu_l\in\ker(M)$, we have
$c=\prod_{i=1}^s c_i^{u_{l, i}} = \prod_{i=1}^s (y^{\f_{i}} - f_i)
^{u_{l, i}} = \prod_{i=1}^s \Y^{M\cdot \bu_l} + g = 1 + g$ where
$g\in\I$. Then $0\ne c -1 \in\I$ and $\I=[1]$, a contradiction.

``$\Leftarrow$" Suppose the contrary. Then there exist
$g_i\in\F\{\Y^{\pm}\}$ such that
 \begin{equation}\label{eq-lg1} g_1 f_1+\cdots+ g_s f_s=
 1.\end{equation}
Let $l$ be the maximal $c$ such that $y_c^{(k)}$ occurs in some
$f_i$, $o$ the largest $j$ such that $y_l^{(j)}$ occurs in some
$f_k$, and $d=\max_{k=1}^s\deg(f_k,y_l^{(o)})$. Let $f_k =
\Y^{\f_k}-c_k = I_k y_l^{dx^o} - c_k$.
Since \bref{eq-lg1} is an identity about the algebraic variables
$y_i^{x^j}$, we can set $y_l^{dx^o} = c_k/I_k$ in \bref{eq-lg1} to
obtain a new identity. In the new identity, $f_k$ becomes zero and
the left hand side of \bref{eq-lg1} has at most $s-1$ summands.
We will show that this procedure can be continued for the new
identity. Then the left hand side of \bref{eq-lg1} will eventually
becomes zero, and a contradiction is obtained and the lemma is
proved.

If $\ord(f_i,y_l) < o$ or $\ord(f_i,y_l) = o$ and
$\deg(f_i,y_l^{x^o}) < d$ for some $i$, then $f_i$ is not changed in
the above procedure. Let us assume that for some $v$,
$\deg(f_v,y_l^{x^o}) = d$ and $f_v = \Y^{\f_v}-c_v = I_v y_l^{dx^o}
- c_v$. Then after the substitution,
$f_v = c_kI_v/I_k - c_v = c_k \widetilde{f}_v$ where
$\widetilde{f}_v=I_v/I_k - c_v/c_k$.
We claim that either $\widetilde{f}_v=0$ or $I_v/I_k$ is a proper
monomial, and as a consequence, the above substitution can continue.
To prove the claim, it suffices to show that if $I_v=I_k$ then
$c_v=c_k$. If $I_v=I_k$, then $\f_v=\f_k$, that is $\f_v-\f_k=0$ is
a syzygy among $\f_i$ and let ${\beps}_{vk}$ be the corresponding
syzygy vector. Then ${\beps}_{vk}\in\ker(M)$ can be written as a
linear combination of $\bu_1,\ldots,\bu_s$. Let
${\bf{c}}=(c_1,\ldots,c_s)^\tau$. Then
$c_vc_k^{-1}={\bf{c}}^{{\beps}_{vk}}$ can be written as a product of
${\bf{c}}^{\bu_l}=\prod_{i=1}^s c_i^{u_{l, i}}=1$, and thus
$c_vc_k^{-1}=1$. \qedd

\subsection{Characteristic set of Laurent binomial $\sigma$-ideal}
\label{sec-lbi2}
We show that the characteristic set method presented in section
\ref{sec-cs} can be modified to the case of Laurent binomial
$\sigma$-ideals.
First, assume that all Laurent $\sigma$-binomials are in normal
form, which makes the concepts of order, leading  variables, etc.
unique.

Second, when defining the concepts of rank and $q$ to be reduced
w.r.t. $p$, we need to replace $\deg(p, y_j^{(o)})$ by $|\deg(p,
y_j^{(o)})|$. Precisely,  $q$ is said to be {\em reduced} w.r.t. $p$
if  $|\deg(q, y_j^{(k+l)})|<|\deg(p, y_j^{(k)})|$  for all $l\in\N$,
where $\lead(p)=y_j^{(k)}$.
For instance, $y_1^{-2x}y_2-1$ is not reduced w.r.t. $y_1^2-1$. With
these changes, the concepts of $\sigma$-chain and characteristic set
can be defined in the Laurent $\sigma$-binomial case.
For instance, the $\sigma$-chains in Example \ref{ex-L11} become the
following Laurent normal form:
\[
 \begin{array}{llllll}
 \A_1& =& y_1^{x}-1, &y_1^2y_2^2-1, &y_2^{x}-1& \\
 \A_2& =& y_1^2-1, &y_1^{-1}y_1^{x}-1, &y_2^2-1, &y_2^{-1}y_2^{x}-1\\
 \A_3& =& y_1^{-x}y_2^2-1,& y_1^{-1}y_3^{2}-1, &y_2^{-1}y_3^{x}-1&\\
\end{array}\]

Third, the $\sigma$-remainder for two Laurent $\sigma$-binomials
need to be modified as follows.
We first consider how to compute $\prem(f,g)$ in the simple case:
$o=\ord(f,y_l)=\ord(g,y_l)$, where $y_l=\lvar(g)$.
Let $g = I_g (y_{l}^{(o)})^{d} - c_g$, where $d=\deg(g,y_{l}^{(o)})$
and $I_g$ is the initial of $g$. As mentioned above,  $g$ is in
normal form,  that is $d>0$. Let $d_f = \deg(f, y_{l}^{(o)})$ and $f
= I_f (y_{l}^{(o)})^{d_f} - c_f$.
We consider two cases.

In the first case,  let us assume  $d_f\ge0$. If $d_f < d_g$,  then
set $r=\prem_1(f, g)$ to be $f$.
Otherwise,  perform the following basic step
 \begin{eqnarray}
r :=\prem_1(f,g)=
 (f - g \frac{I_f}{I_g}(y_{l}^{(o)})^{d_f-d_g})/c_g
  = \frac{I_f}{I_g} (y_{l}^{(o)})^{d_f-d_g} -\frac{c_f}{c_g}. \label{eq-prem1}
 \end{eqnarray}
Let  $\h_r,\h_f,\f_g$ be the supports of $r,f,g,$ respectively. Then
 \begin{eqnarray}\label{eq-prem2}
 \h_r=\h_f-\h_g.
 \end{eqnarray}
Set $f=r$ and repeat the procedure  $\prem_1$ for $f$ and $g$. Since
$d_f$ decreases strictly after each iteration, the procedure will
end and return $r$ which satisfies
 \begin{eqnarray}
 r &=&  \frac{f}{c_g^k} -hg
     =  \frac{I_f}{I_g^{k}} (y_{l}^{(o)})^{d_f-kd_g} -\frac{c_f}{c_g^k}
 \label{eq-prem3}\\
 \h_r&=&\h_f-k \h_g\label{eq-prem4}
 \end{eqnarray}
where $k=\lfloor\frac{d_f}{d_g}\rfloor$ and $h\in\F\{\Y^{\pm}\}$.
Let $\prem(f,g) = r$ or the inverse of $r$ in the case that $r$
is not in normal form.

In the second case,  we assume  $d_f<0$. The $\sigma$-remainder can
be computed similar to the first case. Instead of $g$, we consider
$\widehat{g} = (I_g)^{-1} (y_{l}^{(o)})^{-d_g} - c_g^{-1}$.
If $|d_f| < d_g$,  then set $r=\prem_1(f, g)$ to be $f$.
Otherwise,  perform the following basic step
 \begin{eqnarray*}
 r := \prem_1(f, g)=c_g(f - \widehat{g} I_gI_f (y_{l}^{(o)})^{d_f+d_g})
 = I_f I_g (y_{l_g}^{(o)})^{d_f+d_g} -c_fc_g.
 \end{eqnarray*}
In this case, equation \bref{eq-prem2} becomes $\h_r=\h_f+\h_g$. To
compute $\prem(f, g)$, repeat the above basic step for $f=r$ until
$|d_f| < d_g$.

For two general $\sigma$-binomials $f$ and $g$, $\prem(f,g)$ is
defined as follows: if $f$ is reduced w.r.t $g$, set $\prem(f,g)=f$.
Otherwise, let $y_l=\lvar(g)$, $o_f=\ord(f,y_l)$, and
$o_g=\ord(g,y_l)$. Define
 $$\prem(f,g)=\prem(\ldots,\prem(\prem(f,g^{(o_f-o_g)}),g^{(o_f-o_g-1)}), \ldots,g ).$$
Let $\A: A_1,\ldots,A_s$ be a Laurent binomial $\sigma$-chain and
$f$ a $\sigma$-binomial. Then define
 $$ \prem(f,\A)=\prem(\ldots,\prem(\prem(f,A_s),A_{s-1}),
 \ldots,A_1).$$
In summary,  we have
\begin{lemma}\label{lm-rem1}
Let $\A= A_1,\ldots,A_s$ be a Laurent binomial $\sigma$-chain, $f$ a
$\sigma$-binomial, and $r=\prem(f,\A)$. Then $r$ is  reduced w.r.t.
$\mathcal {A}$ and   satisfies
   \begin{equation}\label{eq-rem2}
   cf \equiv r,  \mod \,  [\mathcal {A}],
   \end{equation}
where $c\in\F^*$. Furthermore, let the supports of $r$ and $f$ be
$\h_r$ and $\h_f$, respectively. Then $\h_f-\h_r$ is in the $\Z[x]$-lattice generated by the supports of $A_i$.
\end{lemma}

Now, the concepts of Laurent binomial regular and coherent
$\sigma$-chains and the characteristic set for Laurent binomial
$\sigma$-ideals can be defined and Theorem \ref{th-rp} can be
extended to the Laurent binomial case.
\begin{theorem}\label{th-rp1}
A Laurent binomial $\sigma$-chain $\A$ is a characteristic set of\,
$\sat(A)$ if and only if $\A$ is regular and coherent.
\end{theorem}

\vskip5pt In the rest of this section, we will establish a relation
between Gr\"obner bases of $\Zx$-lattices and characteristic sets of
Laurent binomial $\sigma$-ideals.
Let $\fb=\{ \f_{1}, \ldots, \f_{s}\}$ be a reduced Gr\"obner basis
of a $\Zx$-lattice such that $\f_1 < \f_2 < \cdots <\f_s$.

\begin{lemma}\label{lm-t2}
If $\fb=\{ \f_{1}, \ldots, \f_{s}\}$ is a reduced Gr\"obner basis of
a $\Zx$-lattice, then $ \Y^{\f_{1}}-c_{1}, \ldots,
\Y^{\f_{s}}-c_{s}$ is a $\sigma$-chain for any $c_i\in\F^*$.
\end{lemma}
\proof Let $A_i = \Y^{\f_{i}}-c_{i}$.
Then for $i< j$,  $A_j$ is of higher rank than $A_i$.  $A_j$ is
reduced w.r.t. $A_i$ if and only if $\f_j$ is G-reduced w.r.t.
$\f_i$, which is valid due to Definition \ref{def-ghf}.\qedd

Note that for a Laurent binomial $\sigma$-chain $\A$, we have
$\sat(\A)=[\A]$.
\begin{lemma}\label{lm-t3}
Let $\fb=\{ \f_{1}, \ldots, \f_{s}\}$ be a reduced Gr\"obner basis
of a $\Zx$-lattice and $\A$ the $\sigma$-chain $\Y^{\f_{1}}-c_{1},
\ldots, \Y^{\f_{s}}-c_{s}$.
If $\I=[\A]$ is a proper Laurent binomial $\sigma$-ideal, then for
$f=\Y^{\f}-c_\f\in\I$, $\grem(\f,\fb)=0$ if and only if
$\prem(f,\A)=0$.
\end{lemma}
\proof Let us first consider  $\prem_1$ in \bref{eq-prem1} for $f$
and $A_i=\Y^{\f_i}-c_i=I_i (y_{l_i}^{(o_i)})^{d_i} - c_i$, where
$\lead(A_i)=y_{l_i}$ and $I_i$ is the initial of $A_i$.
From \bref{eq-prem2}, the support of $r=\prem_1(f,A_i)$ is  $\f -
\f_i$.

It is clear that $\LT(\f_i) = d_i x^{o_i} {\beps}_{l_i}$. Let $\f_i
=  d_i x^{o_i} {\beps}_{l_i}+ \overline{\f}_i$. Similarly, write $\f
= d_f x^{o_i}{\beps}_{l_i} +\overline{m}$ where $d_f
x^{o_i}{\beps}_{l_i}$ is the leading term of $\f$ w.r.t.
${\beps}_{l_i}$ and $d_f\ge d_i\ge0$.
Then a basic step to compute $\grem(\f,\f_i)$ is to compute
$\grem_1(\f,\f_i)= \f - \f_i = (d_f-d_i) x^{o_i}{\beps}_{l_i}
+\overline{\f}- \overline{\f}_i$, which is the support of
$\prem_1(f,A_i)$.

As a consequence, using the basic step $\grem_1$ to compute
$\grem(\f,\fb)$, we have a sequence of elements in $\Z[x]^n$:
$\g_0=\f, \g_1,\ldots,\g_t=\grem(\f,\fb)$.
Correspondingly,  using the basic step $\prem_1$ to compute
$\prem(\f,\A)$, we have a sequence of $\sigma$-binomials
$f_0=f,f_1,\ldots,f_t=\prem(f,\A)$ such that the support of $f_i$ is
$\g_i$.
Let $f_t = \Y^{\g_t} - \widetilde{c}_t$. Since $\I$ is proper,
$f_t=\prem(f,\A)=\Y^{\g_t} - \widetilde{c}_t=0$ if and only if $\g_t
=\grem(\f,\fb)= 0$.\qedd

\begin{lemma}\label{lm-t4}
Let $\fb=\{ \f_{1}, \ldots, \f_{s}\}$ be a reduced Gr\"obner basis
of a $\Zx$-lattice and $\A$ the $\sigma$-chain $\Y^{\f_{1}}-c_{1},
\ldots, \Y^{\f_{s}}-c_{s}$.
If $\I=[\A]$ is a proper Laurent binomial $\sigma$-ideal, then $\A$
is a regular and coherent $\sigma$-chain.
\end{lemma}
\proof By Lemma \ref{lm-t2}, $\A$ is a $\sigma$-chain. Since the
initials of $\A$ are $\sigma$-monomials which are units in
$\F\{\Y^{\pm}\}$, $\A$ is regular. We still need to prove that $\A$
is coherent.

%
%
Let $A_i=\Y^{\f_{i}}-c_{i}$ and $A_j=\Y^{\f_{j}}-c_{j}$ ($i<j$) have
the same leading variable $y_l$, and
    $A_i=I_i y_l^{d_i x^{o_i}}-c_i$,
    $A_j=I_j y_l^{d_j x^{o_j}}-c_j$.
From Definition \ref{def-ghf},  we have $o_i < o_j$ and $d_j| d_i$.
Let $d_i = t d_j$ where $t\in\N$.
According to \bref{eq-prem3}, we have
 \begin{equation}\label{eq-delta1}
 \Delta(A_i,A_j)
=\prem((A_i)^{x^{o_j-o_i}},A_j) =
       \frac{(I_i)^{x^{o_j-o_i}}}{I_j^t} -
         \frac{(c_i)^{x^{o_j-o_i}}}{c_j^t}.
 \end{equation}
Then the support of $\Delta(A_i,A_j)$ is $x^{o_j-o_i} \f_i -
\frac{d_i}{d_j} \f_j$.

Since $\LT(A_i)=d_i x^{o_i} {\beps}_l$ and  $\LT(A_j)=d_j x^{o_j}
{\beps}_l$,  we have $N=\lcm(d_i x^{o_i}, d_j x^{o_j}) = d_i
x^{o_j}$.
According to Definition \ref{def-sv}, the S-vector of $\f_i$ and
$\f_j$ is
 $$S(\f_i,\f_j) = x^{o_j-o_i} \f_i - \frac{d_i}{d_j}  \f_j.$$
Since $\fb$ is a Gr\"obner basis, we have
$\grem(S(\f_i,\f_j),\fb)=0$. Also note that the support of
$\Delta(A_i,A_j)$ is $S(\f_i,\f_j)$. Then by Lemma \ref{lm-t3},
$\prem(\Delta(A_i,A_j),\A)=0$, that is, $\A$ is coherent.\qedd

We summarize the results in this section as the following theorems.
\begin{theorem}\label{th-t1}
Let $\fb=\{ \f_{1}, \ldots, \f_{s}\}$ be a reduced Gr\"obner basis
of a $\Zx$-lattice and $\A$ the $\sigma$-chain $\Y^{\f_{1}}-c_{1},
\ldots, \Y^{\f_{s}}-c_{s}$.
If $\I=[\A]$ is a proper Laurent binomial $\sigma$-ideal, then $\A$
is a characteristic set of $\I$.
\end{theorem}
\proof  This is a direct consequence of Lemma \ref{lm-t4} and
Theorem \ref{th-rp1}.\qedd

In the following theorem, we show how to compute the characteristic
set for a Laurent binomial $\sigma$-ideal using Gr\"obner bases of
$\Zx$-lattices.
\begin{theorem}\label{th-t2}
Let $\I=[\Y^{\f_{1}}-c_{1}, \ldots, \Y^{\f_{s}}-c_{s}]$ be a proper
Laurent binomial $\sigma$-ideal, $\hb=[\h_1, \ldots, \h_r]$ the
reduced Gr\"obner basis of the $\Zx$-lattice $(\f_1,\ldots,\f_s)$,
and
 $\h_i = \sum_{k=1}^s a_{i,k} \f_k$
for $a_{i,k}\in\Zx$. 
 Let $\A$ be $\Y^{\h_1}-\prod_{i=1}^s c_i^{a_{1, i}} ,$ $\ldots,
  \Y^{\h_r}-\prod_{i=1}^s c_i^{a_{r, i}}$. Then $\I=[\A]$ and $\A$ is a characteristic set of
   $\I$.
\end{theorem}
\proof By Lemma \ref{lm-bi3}, $\I = [\A]$. By Theorem \ref{th-t1},
$\A$ is a characteristic set of $\I$.\qedd

\begin{cor}\label{cor-dim1}
Let $\I$ be a Laurent reflexive prime  binomial $\sigma$-ideal in
$\F\{\Y^{\pm}\}$. Then $\dim(\I)=n-\rk(\L(\I))$.
\end{cor}
\proof Suppose $\C=[\c_1,\ldots,\c_s]$ in \bref{ghf} is the matrix
representation for $\I$. Since $\I$ is reflexive and prime, $\I$ has
a characteristic set of form $\A:\Y^{\c_1}-c_1,\ldots\Y^{\c_s}-c_s$.
By Theorem 4.3 of \cite{gao-dcs}, $\dim(\I)=n-t=n-\rk(\L(\I))$.\qedd

\begin{cor}\label{cor-rad}
A Laurent binomial $\sigma$-ideal is radical.
\end{cor}
\proof By Theorem \ref{th-t2}, $\I=[\Y^{\h_1}-c_1, \ldots,
\Y^{\h_r}-c_r]$, where $\A: \Y^{\h_1}-c_1, \ldots, \Y^{\h_r}-c_r$ is
the characteristic set of $\I$.
Let $A_i=\Y^{\h_i}-c_i$ and $y_{l_i}^{(o_i)} = \lead(A_i)$. $\A$ is
also {\em saturated} in the sense that its separant $\frac{\partial
A_i}{\partial y_{l_i}^{(o_i)}}$ are $\sigma$-monomials and  hence
units in $\F\{\Y^{\pm}\}$. Then similar to the differential case
\cite{bouziane}, it can be shown that $\sat(\A) = [\A]$ is a radical
$\sigma$-ideal.\qed

Finally, we show that the converse of Lemma \ref{lm-t4} is also
true.
\begin{theorem}\label{th-t4}
Let $\f_i\in\Zxn$ and $c_i\in\F^*$ for $i=1,\ldots,s$.
$\A:\Y^{\f_{1}}-c_1, \ldots, \Y^{\f_{s}}-c_s$  is a regular and
coherent $\sigma$-chain if and only if $\fb=\{\f_{1}, \ldots,
\f_{s}\}$ is a reduced Gr\"obner basis and $[\A]$ is a proper
Laurent binomial $\sigma$-ideal.
In this case, the support lattice of $[\A]$ is $(\fb)$.
\end{theorem}
\proof Lemma \ref{lm-t4} proves one side of the theorem.
For the other direction, let $\A$ be a regular and coherent
$\sigma$-chain. From the proof of Lemma \ref{lm-t2}, $\f_i$ is
G-reduced to $\f_i$ for $i\ne j$.
By Theorem \ref{th-rp1}, $\A$ is a characteristic set of
$\sat(\A)=[\A]$, which means $[\A]$ is proper.
Use the notations introduced in the proof of Lemma \ref{lm-t4}.
Since $S(\f_i,\f_j)$ is the support of $\Delta(A_i,A_j)$, by the
proof of Lemma \ref{lm-t3}, $\f_{ij}=\grem(S(\f_i,\f_j),\fb)$ is the
support of $\prem(\Delta(A_i,A_j),\A)$.
Since $\A$ is coherent, $\prem(\Delta(A_i,A_j),\A)=\Y^{\f_{ij}}-c=0$
for any $i$ and $j$, which implies
$\f_{ij}=\grem(S(\f_i,\f_j),\fb)=0$ and hence $\fb$ is a reduced
Gr\"obner basis.

It remains to show that the support lattice of $[\A]$ is $(\fb)$. By
Theorem \ref{th-t1} and Lemma \ref{lm-t3}, $f=\Y^{\g}-c\in[\A]$ if
and only if $\prem(f,\A)=0$, which is equivalent to
$\grem(\g,\fb)=0$, that is $\L([\A]) = (\fb)$.\qedd
%

\subsection{Partial character and Laurent binomial $\sigma$-ideal}
\label{sec-lbi3} In this section,  we will show that proper Laurent
binomial $\sigma$-ideals can be described uniquely with their
partial characters.

\begin{definition}
A {\em partial character} $\rho$ on $\Z[x]^{n}$ is a homomorphism
from a $\Z[x]$-lattice $L_{\rho}$ to the multiplicative group
$\F^{\ast}$ satisfying $\rho(x\f)=\sigma(\rho(\f))$ for $\f\in
L_\rho$.
\end{definition}
Partial characters can be defined on any $\Z[x]$-lattice $L$. A {\em
trivial} partial character on $L$ is defined by setting $\rho(\f)=1$ for any $\f\in L$.

Let $\rho$ be a partial character over $\Zxn$ and $L_\rho=(\f_{1},
\ldots, \f_{s})$, where $\fb=\{\f_{1}, \ldots, \f_{s}\}$ is a
reduced Gr\"obner basis. Define
\begin{eqnarray}
 \I(\rho)&:=&[\Y^{\f}-\rho(\f)\,|\, \f\in L_{\rho}].\label{eq-I}\\
 \A(\rho)&:=& \Y^{\f_{1}}-\rho(\f_{1}), \ldots, y^{\f_{s}}-\rho(\f_{s}).\label{eq-Arho}
\end{eqnarray}
The Laurent binomial $\sigma$-ideal $\I(\rho)$ has the following
properties.


\begin{lemma}\label{lm-l1}
For $\rho$\,　and $\A$ defined above,
$\I(\rho)=[\A(\rho)]\ne [1]$ and $\A(\rho)$ is a characteristic set
of $\I(\rho)$.
\end{lemma}
\proof By Lemma \ref{lm-bi1}, $\I(\rho)=[\A(\rho)]$.
%
%
By Lemma \ref{lm-bi2}, in order to prove $\I(\rho)\ne [1]$, it
suffices to show that for any syzygy $\sum_i a_i \f_i = 0$ among
$\f_i$, we have $\prod_i \rho(\f_i)^{a_i}=1$. Indeed,  $\rho(\sum_i
a_i \f_i)=\prod_i\rho(\f_i)^{a_i}=1$, since $\rho$ is a homomorphism
from the $\Zx$-module $L_\rho$ to $\F^{*}$.
Since $\fb$ is a reduced Gr\"ober basis, by Theorem \ref{th-t1},
$\A$ is a characteristic set of $\I(\rho)$.\qedd

\begin{lemma}\label{lm-l3}
A Laurent $\sigma$-binomial $\Y^{\f}-c_\f$ is in $\I(\rho)$ if and
only if $\f\in L_\rho$ and $c_\f=\rho(\f)$.
\end{lemma}
\proof By Lemma \ref{lm-l1}, $\A(\rho)$ is a characteristic set of
$\I(\rho)$.
Since $f=y^{\f}-c_{\f}$ is a $\sigma$-binomial in $\I(\rho)$, we
have $r=\prem(f,\A)=0$. By Lemma \ref{lm-rem1}, $\f$ is in the
$\Z[x]$-module $L_\rho$. The other side is obviously true and the
lemma is proved.\qedd

We now show that all Laurent binomial $\sigma$-ideals are defined by
partial characters.
\begin{theorem}\label{th-l1}
The map $\rho\Rightarrow\I(\rho)$  gives a one to one correspondence
between the set of proper Laurent binomial $\sigma$-ideals and
partial characters on $\Z[X]^n$.
\end{theorem}
\proof By Lemma \ref{lm-l1}, a partial character defined a proper
Laurent binomial $\sigma$-ideal.
For the other side, let $\I\subseteq \F\{\Y^{\pm}\}$ be a proper
Laurent binomial $\sigma$-ideal.
$\I$ is generated by its members of the form $y^{\f}-c_{\f}$ for
$\f\in \Z[x]^{n}$ and $c_\f\in \F^{\ast}$. Let $L_\rho = \L(\I)$
which is defined in  \bref{eq-bi1} and  $\rho(\f)=c_\f$. Since $\I$
is proper, $c_{\f}$ is uniquely determined by $\f$. By Lemma
\ref{lm-bi1} and Corollary \ref{cor-bi1}, $\rho$ is a partial
character which is uniquely determined by $\I$. It is clear
$\I(\rho)=\I$.
To show the correspondence is one to one, it suffices to show
$\rho(\I(\rho)) = \rho$ which is a consequence of Lemma \ref{lm-l3}.
The theorem is proved.\qedd

From Lemma \ref{lm-l1} and Theorem \ref{th-l1},  we have
\begin{cor}\label{cor-l1}
Any Laurent binomial $\sigma$-ideal is finitely generated.
\end{cor}

\begin{cor}\label{cor-l11}
Let $\f_i\in\Zxn$ and $c_i\in\F^*$ for $i=1,\ldots,s$.
If $\A:\Y^{\f_{1}}-c_1, \ldots, \Y^{\f_{s}}-c_s$ is a regular and
coherent $\sigma$-chain, then there exists a partial character
$\rho$ over $\Zxn$ such that $L_\rho=(\f_1,\ldots,\f_s)$,
$\rho(\f_i)=c_i$, and $\I(\rho)=[\A]$.
\end{cor}
\proof By Theorem \ref{th-t4}, $[\A]$ is proper with support lattice
$(\f_1,\ldots,\f_s)$. By Theorem \ref{th-l1}, the corresponding
partial character of $[\A]$ satisfies properties in the
corollary.\qedd

Now, we can prove Theorem \ref{th-m1} easily.

{\em Proof of Theorem \ref{th-m1}}. By Theorem \ref{th-l1}, $\I$ is
a proper Laurent binomial $\sigma$-ideal in $\F\{\Y^{\pm}\}$ if and
only if $\I=\I(\rho)$ for a partial character $\rho$, hence (3).
By Theorem \ref{th-t4}, (1) and (2) of Theorem \ref{th-m1} are
equivalent.
From Lemma \ref{lm-l1} and Corollary \ref{lm-t4}, (3) implies (2).
From Corollary \ref{cor-l11}, (2) implies (3).
\qedd

\subsection{Reflexive and prime  Laurent binomial $\sigma$-ideals}
\label{sec-pr1}

In this section, we first give criteria for reflexive and prime
Laurent binomial $\sigma$-ideals and then give a decomposition
theorem for perfect Laurent binomial $\sigma$-ideals.
Let $a$ be an element in an over field of $\F$, $S$ a set of
elements in an over field of $\F$, and $k\in\N$. Denote $a^{[k]} =
\{a,a^x,\ldots,a^{x^k} \}$ and $S^{[k]} = \cup_{b\in S} b^{[k]}$.

For the $\sigma$-indeterminates $\Y=\{y_1,\ldots,y_n\}$ and
$t\in\N$, we will treat the elements of $\Y^{[t]}$ as algebraic
indeterminates, and $\F[\Y^{[\pm t]}]$ is the Laurent polynomial
ring in $\Y^{[t]}$.
Let $\I$ be a Laurent binomial $\sigma$-ideal in $\F\{\Y^{\pm}\}$.
Then it is easy to check that
$$\I_t =\I\cap \F[\Y^{[\pm t]}]$$
is a Laurent binomial ideal in $\F[\Y^{[\pm t]}]$.

Denote $\Zx_t$ to be the set of elements in $\Zx$ with degree $\le
t$. Then $\Zxn_t$ is the $\Z$-module generated by $x^i{\beps}_l$ for
$i=0,\ldots,t,l=1,\ldots,n$. It is clear that $\Zxn_t$ is isomorphic
to $\Z^{n(t+1)}$ as $\Z$-modules by mapping $x^i{\beps}_l$ to the
$((l-1)(t+1) + i+1)$-th standard basis vector in $\Z^{n(t+1)}$. Hence,
we treat them as the same in this section.
Let $L$ be a $\Zx$-lattice and $t\in\N$. Then
 $$L_t =L\cap \Zxn_t =L\cap \Z^{n(t+1)}$$
is a $\Z$-module in $\Z^{n(t+1)}$.
Similarly, it can be shown that when restricting to $\Zxn_t$, a
partial character $\rho$ on $\Zxn$  becomes a partial character
$\rho_t$ on $\Z^{n(t+1)}$.

\begin{lemma}\label{lm-pr1}
With the notations introduced above,  we have $\I_t =\I\cap
\F[\Y^{[\pm t]}] =\I(\rho_t)$.
\end{lemma}
\proof It suffices to show that the support lattice of $\I_t$ is
$L_{\rho_t} = L_t$. By Lemma \ref{lm-l3}, $\Y^\f-c_m\in\I_t$ if and
only if $\f\in L\cap\Zxn_t$, or equivalently, $\max_{m\in\f}
\deg(m,x)\le t$, which is equivalent to $\f\in L_t$.\qedd
%
%

\begin{definition} Let $L$ be a $\Z[x]$-module in $\Z[x]^n$.
\begin{itemize}
\item $L$ is called {\em $\Z$-saturated} if,  for any
$a\in\Z$ and $\f\in\Z[x]^n$,  $a\f\in L$ implies $\f\in L$.

\item
$L$ is called {\em $x$-saturated} if,  for any $\f\in\Z[x]^n$,
$x\f\in L$ implies $\f\in L$.

\item
$L$ is called {\em saturated} if it is both $\Z$- and $x$-
saturated.
\end{itemize}
\end{definition}

We now prove (1)-(2) of Theorem \ref{th-m2}, that is
\begin{theorem}\label{th-pr1}
Let $\rho$ be a partial character over $\Zxn$. If $\F$ is
algebraically closed and inversive, then
\begin{description}
\item[(a)] $L_\rho$ is $\Z$-saturated if and only if $\I(\rho)$ is prime;
\item[(b)] $L_\rho$ is $x$-saturated if and only if $\I(\rho)$ is reflexive;
\item[(c)] $L_\rho$ is saturated if and only if $\I(\rho)$ is reflexive prime.
\end{description}
\end{theorem}
\proof  It is clear that (c) comes from (a) and (b). Let
$\I=\I(\rho)$ and $L=L_\rho$.

$(a)$: $\I$ is a Laurent prime $\sigma$-ideal if and only if $\I_t$
is a Laurent prime ideal for all  $t$.
From Lemma \ref{lm-pr1}, the support of $\I_t$ is $L_t$. Then by
\cite[Thm 2.1]{es-bi}, $\I_t$ is a Laurent prime ideal if and only
if $L_t$ is a $\Z$-saturated $\Z$-module.
Furthermore, a $\Zx$-lattice $L$ is $\Z$-saturated if and only if
$L_t$ is a $\Z$-saturated $\Z$-module for all $t$. Thus, (a) is
valid.

(b): Suppose $\I$ is reflexive. For $x\f\in L$, by Lemma
\ref{lm-l3}, there is a $\Y^{x\f}-c \in \I$. Since $\F$ is
reflexive, $c = d^x$ for $d\in\F$. Then $\sigma(\Y^{\f}-d) \in \I$
and hence $\Y^{\f} - d\in\I$ since $\I$ is reflexive. By Lemma
\ref{lm-l3} again, $\f\in L$ and $L$ is $x$-saturated.
To prove the other direction, assume $L$ is $x$-saturated.
For $f^{x} \in \I$, we have an expression
\begin{equation}\label{eq-pr2}
f^{x}=\sum_{i=1}^{s}f_{i}(\Y^{\f_{i}}-c_i)
\end{equation}
where $\Y^{\f_{i}}-c_i\in\I$ and $f_{i}\in \F\{\Y^{\pm}\}$.
Let $d=\max_{i=1}^s\deg(\Y^{\f_{i}}-c_i,y_1)$ and assume
$\Y^{\f_{1}} =M_1 y_1^d$.
Replace $y_1^d$ by $c_1/M_1$ in \bref{eq-pr2}. Since \bref{eq-pr2}
is an identity for the variables $y_i^{(j)}$, this replacement is
meaningful and we obtain a new identity. $\Y^{\f_{1}}-c_1$ becomes
zero after the replacement. Due to the way to chose $d$, if another
summand, say $\Y^{\f_{2}}-c_2$, is affected by the replacement, then
$\Y^{\f_{2}} = M_2 y_1^d$. After the replacement, $\Y^{\f_{2}}-c_2$
becomes $c_1(M_2/M_1-c_2/c_1)$ which is also in $\I$ by Lemma
\ref{lm-l3}.
In summary, after the replacement, the right hand side of
\bref{eq-pr2} has less than $s$ summands and the left hand side of
\bref{eq-pr2} does not changed. Repeat the above procedure, we will
eventually obtain a new indenity
\begin{equation}\label{eq-pr21}
f^{x}=\sum_{i=1}^{\bar{s}} \bar{f}_{i}(\Y^{x\g_{i}}-\bar{c}_i)
\end{equation}
where $\Y^{x\g_{i}}-\bar{c}_i\in\I$ and $\bar{f}_{i}\in
\F\{\Y^{\pm}\}$.
We may assume that any $y_i$ does not appear in $\bar{f}_{i}$.
Otherwise, by setting $y_i$ to be $1$, the left hand side of
\bref{eq-pr21} is not changes and a new identity is obtained.
Since $\F$ is inversive, $\bar{c}_i = e_i^x$ and $\bar{f}_{i} =
g_i^x$ for $e_i\in\F$ and $g_i\in\F\{\Y^{\pm}\}$.
By Lemma \ref{lm-l3}, $\Y^{x\g_{i}}-e_i^x \in I$ implies $x\g_i\in
L$. Since $L$ is $x$-saturated, $x\g_i\in L$ implies $\g_i\in L$ and
hence $\Y^{\g_{i}}-e_i \in I$ by Lemma \ref{lm-l3} again.
From \bref{eq-pr21}, $\sigma(f-\sum_{i=1}^{\bar{s}}
g_{i}(\Y^{\g_{i}}-e_i))=0$ and hence $f=\sum_{i=1}^{\bar{s}}
g_{i}(\Y^{\g_{i}}-e_i)\in \I$. (b) is proved.\qedd

Decision procedures for whether a $\Zx$-lattice is $\Z$-saturated or
$x$-saturated will be given in Section \ref{sec-alg}.

%
From Corollary \ref{cor-rad}, every proper Laurent binomial
$\sigma$-ideal is radical. The following example shows that a
binomial $\sigma$-ideal is not necessarily perfect.
\begin{example}\label{ex-31}
Let $\I=[\A]\subset\Q\{y_1,y_2\}$,  where $\A = \{y_1^2+1,
y_1^x-y_1, y_2^2+1, y_2^x+y_2\}$ is a $\sigma$-chain. We have
$y_2^2+1-(y_1^2+1) = (y_2-y_1)(y_2+y_1)$. Then $\{\I\} =\{\I,
y_2-y_1\}\cap\{\I, y_2+y_1\} = [1]$ since $[\I, y_2-y_1]=[\I,
y_2+y_1]=[1]$. On the other hand,  $1\not\in\I$, since $\A$ is a
characteristic set of $\I$ by Theorem \ref{th-t1}.
%
\end{example}

\begin{definition}
Let $L\subset\Z[x]^{n}$ be a $\Z[x]$-lattice.
The {\em $\Z$-saturation} of $L$ is $\sat_{\Z}(L)=\{\f\in
\Z[x]^{n}\, |\, \exists a \in \Z \st a\f\in L\}$.
The {\em $x$-saturation} of $L$ is $\sat_{x}(L)=\{\f\in \Z[x]^{n}\,
|\, x \f\in L\}$.
The {\em saturation} of $L$ is $\sat(L)=\{\f\in \Z[x]^{n}\, |\,
\exists a\in\Z, \exists k \in \N \st ax^k \f\in L\}$.
\end{definition}
It is clear that the $\Z$-saturation ($x$-saturation) of $L$ is
$\Z$-saturated ($x$-saturated) and
 $$\sat(L) =
\sat_{\Z}(\sat_{x}(L))=\sat_{x}(\sat_{\Z}(L)).$$
Algorithms to compute the $\Z$-saturation and $x$-saturation of a
$\Zx$-lattice will be given in Section \ref{sec-alg}.

\begin{theorem}\label{lm-ld0}
Let $\I$ be a Laurent binomial $\sigma$-ideal and $L$ the support
lattice of $\I$. If $\F$ is inversive, then the reflexive closure of
$\I$ is also a Laurent binomial $\sigma$-ideal whose support lattice
is the $x$-saturation of $L$.
\end{theorem}
\proof Let $\I_x$ be the reflexive closure of $\I$ and
$L_x=\sat_x(L)$. Suppose $\I=[f_1,\ldots,f_r]$, where $f_i =
\Y^{\f_i}-c_i$. Then $L=(\f_1,\ldots,\f_r)$.
If $L$ is $x$-saturated, by Theorem \ref{th-pr1}, $\I$ is reflexive.
Otherwise,  there exist $k_1\in\N$, $b_{i}\in\Zx$, and $\h_1\in\Zxn$
such that $\h_1\not\in L$ and
 \begin{equation}\label{eq-ldec1}
  x^{k_1} \h_1 = \sum_{i=1}^r b_{i} \f_i\in L.\end{equation}
By Lemma \ref{lm-bi1}, $\Y^{x^{k_1} \h_1} - \widetilde{a}$ is in
$\I$, where $\widetilde{a} = \prod_{i=1}^r c_i^{b_{i}}$. Since $\F$
is inversive, $\overline{a}=\sigma^{-k_1}(\widetilde{a})\in\F$.
Then, $\sigma^{k_1} (\Y^{\h_1} - \overline{a})\in\I$, and hence
$\Y^{\h_1} - \overline{a}\in \I_x$.
Let $\I_1=[f_1,\ldots,f_r,\Y^{\h_1} - \overline{a}]$. It is clear
that $L_1=(\f_1,\ldots,\f_r,\h_1)$ is the support lattice of $\I_1$.
Then $\I\varsubsetneq\I_1\subset \I_x$ and $L\varsubsetneq
L_1\subset L_x$.
Repeating the above procedure for $\I_1$ and $L_1$, we obtain $\I_2$
and $L_2=(\f_1,\ldots,\f_r,\h_1,\h_2)$ such that $\h_2\not\in L_1$
and $x^{k_2}\h_2\in L_1$.
We claim that $L_2\subset L_x$. Indeed, let $x^{k_2}\h_2 =
\sum_{i=1}^r e_i \f_i + e_0 \h_1$. Then by \bref{eq-ldec1},
$x^{k_1+k_2}\h_2=x^{k_1} (x^{k_2}\h_2)
 =x^{k_1}\sum_{i=1}^r e_i \f_i + e_0 (x^{k_1}\h_1)
 =x^{k_1}\sum_{i=1}^r e_i \f_i + e_0 \sum_{i=1}^r b_i \f_i\in L$ and
 the claim is proved.
As a consequence, $\I_2\subset \I_x$.

Continuing the process, we have
 $\I\varsubsetneq\I_1\varsubsetneq\cdots\varsubsetneq\I_t\subset \I_x$ and $L\varsubsetneq
L_1\varsubsetneq\cdots\varsubsetneq L_t\subset L_x$ such that $L_i$
is the support lattice of $\I_i$.
The process will terminate, since $\Zxn$ is Northerian. The final
$\Zx$-lattice $L_t$ is $x$-saturated and hence $\I_t$ is reflexive
by Theorem \ref{th-pr1}. Since $L_x$ is the smallest $x$-saturated
$\Zx$-lattice containing $L$ and $L\subset L_t\subset L_x$, we have
$L_t=L_x$ and $\I_t=\I_x$.\qedd

\begin{cor}\label{cor-ld0}
Let $L\subset\Z[x]^{n}$ be a $\Z[x]$-lattice. Then $L_x=\sat_x(L)$
is generated by $L$ and a finite number of elements which are linear
combinations of elements of $L$ divided by certain $x^d$.
Furthermore, $\rk(L)=\rk(L_x)$.
\end{cor}
\proof From the proof of Theorem \ref{lm-ld0}, $\sat_x(L)=
(L,\h_1,\ldots,\h_t)$ and  for each $\h_i$, there is a positive
integer $n_i$ such that $x^{n_i}\h_i\in L$.
Let $A$ be a representation matrix of $L$. Then a representation
matrix $B$ of $L_x$ can be obtained by adding to $A$ a finite number
of new columns which are linear combinations of columns of $A$
divided by some $x^d$. Therefore, $\rk(A)=\rk(B)$.\qedd

Similarly, we can show
\begin{lemma}\label{lm-ld1}
Let $L\subset\Z[x]^{n}$ be a $\Z[x]$-lattice. Then $L_\Z=\sat_\Z(L)$
is generated by $L$ and a finite number of elements which are linear
combinations of elements of $L$ divided by an integer. Furthermore,
$\rk(L)=\rk(L_\Z)$.
\end{lemma}


By Corollary \ref{cor-rad}, a Laurent binomial $\sigma$-ideal $\I$
is radical. By Example \ref{ex-31}, $\I$ is not necessarily perfect.
We now give a decomposition theorem for perfect $\sigma$-ideals.
\begin{theorem}\label{th-l2}
Let $\I$ be a Laurent binomial $\sigma$-ideal,  $L$ the support
lattice of $\I$, and $L_S$ the saturation of $L$. If $\F$ is
algebraically closed and inversive,  then $\{\I\}$ is either $[1]$
or can be written as the intersection of Laurent reflexive prime
binomial $\sigma$-ideals whose support lattice is $L_S$.
\end{theorem}
\proof  Let $\I_x$ be the reflexive closure of $\I$ and
$L_x=\sat_x(L)$. By Theorem \ref{lm-ld0}, $L_x$ is the support lattice
of $\I_x$.
Suppose $\I_x=[f_1,\ldots,f_r]$, $f_i = \Y^{\f_i} -
c_i,i=1,\ldots,r$, and $L_x=(\f_1,\ldots,\f_r)$. If $L_x$ is
$\Z$-saturated, then by Theorem \ref{th-pr1}, $\I_x$ is reflexive
prime. Otherwise, there exist $k_1\in \N$, $a_i\in\Z[x]$, and
$\h_1\in \Z[x]^n$ such that\ $\h_1\not\in L_x$ and
 \begin{equation}\label{eq-ldec2}
  k_1\h_1 = a_1\f_1 + \cdots + a_r \f_r \in L_x.
 \end{equation}
By Lemma \ref{lm-bi1}, $\Y^{k_1\h_1}- \widetilde{a}\in\I$, where
$\widetilde{a}=\prod_{i=1}^r c_i^{a_i}$. Since $\F$ is algebraically
closed, we have
$$\Y^{k_1\h_1}- \widetilde{a}= \prod_{l=1}^{k_1} (\Y^{\h_1}-\widetilde{a}_{l})\in\I_x$$
where $\widetilde{a}_{l}, l=1, \ldots, k_1$ are the $k_1$ roots of
$\widetilde{a}$. By the difference Nullstellensatz
\cite[p.87]{cohn}, we have the following decomposition
$$\{\I\} = \cap_{l_1=1}^{k_1} \{\I_{1l_1}\}$$
where $\I_{1l} = [f_1,\ldots,f_r, \Y^{\h_1}-\widetilde{a}_{l}]$.
Check whether $\I_{1l_1}=[1]$ with Lemma \ref{lm-bi2} and discard
those trivial ones. Then the support lattice for any of $\I_{1l}$ is
$L_1 =(\f_1,\ldots,\f_r,\h_1)$. Similar to the proof of Theorem
\ref{lm-ld0}, we can show that $\I_x\varsubsetneq\I_{1l}$ and
$L_x\varsubsetneq L_1 \subset  L_S$.

Repeating the process, we have
 $\I_x\varsubsetneq\I_{1l_1}\varsubsetneq\cdots\varsubsetneq\I_{tl_t}$ for $l_i=1,\ldots,k_i$
and $L_x\varsubsetneq L_1\varsubsetneq
L_2\varsubsetneq\cdots\varsubsetneq L_t\subset L_S$ such that $L_i$
is the support lattice of $\I_{il_i}$ for $l_i=1,\ldots,k_i$ and
 $$\{\I\}=\cap_{l_i=1}^{k_i}\{\I_{il_i}\}, i=1,\ldots,t.$$
The process will terminate, since $\Zxn$ is Northerian. Since $L_S$
is the smallest $\Z$-saturated $\Zx$-lattice containing $L_x$ and
$L_x\subset L_t\subset L_S$, we have $L_t=\sat_{\Z}(L_x)
=\sat_{\Z}(\sat_{x}(L)) = L_S$. Then $\I_{tl_t}$ is reflexive prime
and the theorem is proved. \qedd

Since the reflexive prime components of $\I$ have the same support
lattice, they also have the same dimension.
\begin{cor}\label{cor-l2}
Any Laurent binomial $\sigma$-ideal $\I$ is dimensionally unmixed.
\end{cor}

The condition for $\F$ to be algebraically closed and inversive is
necessary for Theorem \ref{th-l2} to be valid. For instance, if
$\F=\Q(\lambda)$ and $\sigma(f(\lambda)) = f(\lambda^2)$, then
$[y_1^2y_2^2-\lambda]$ is prime and its support lattice
$([2,2]^\tau)$ is not $\Z$-saturated.
Furthermore, $[y_1^xy_2^x-\lambda]$ is reflexive and its support
lattice $([x,x]^\tau)$ is not $x$-saturated.

\subsection{Well-mixed and perfect Laurent binomial $\sigma$-ideals}
\label{sec-pr3}

In this section, we give a criterion for a Laurent binomial
$\sigma$-ideal to be well-mixed and perfect in terms of its support
lattice and show that the well-mixed and perfect closures of a
Laurent binomial $\sigma$-ideal are still binomial.

Recall that a $\sigma$-ideal $\I$ is called {\em well-mixed} if
$fg\in\I$ implies $fg^x\in\I$ for $f,g\in\F\{\Y^{\pm}\}$. The
smallest well-mixed $\sigma$-ideal containing
$S\subset\F\{\Y^{\pm}\}$ is denoted by $\langle S\rangle$. Let $S' =
\{ fg^x | fg\in S \}$. We define inductively: $S_0 = S, S_n =
[S_{n-1}]', n=1, 2, \ldots$. The union of the $S_n$ is clearly a
well-mixed $\sigma$-ideal and is contained in every well-mixed
$\sigma$-ideal containing $S$. Hence this union is $\langle
S\rangle$. If $\I\subset\F\{\Y^{\pm}\}$ is a Laurent $\sigma$-ideal,
then $\langle\I\rangle$ ia called the well-mixed closure of $\I$.
We first prove some basic properties of well-mixed $\sigma$-ideals
which will be used later. Note that these properties are also valid
in $\F\{\Y\}$.

\begin{lemma}\label{lm-wm}
Let $\I_1, \ldots,\I_s$ be prime $\sigma$-ideals.
Then $\I = \cap_{i=1}^s \I_i$ is a well-mixed $\sigma$-ideal.
\end{lemma}
\proof It is obvious.\qedd

\begin{lemma}\label{lm-wm1}
Let $S_1, S_2$ be two subsets of $\F\{\Y^{\pm}\}$ which satisfy $a\in S_i$ implies $\sigma(a)\in S_i, i=1,2$. Then $[S_1]_n [S_2]_n \subset [S_1S_2]_{n}$.
\end{lemma}
\proof Let $s\in [S_1]_1$ and $t\in [S_2]_1$. Then $s= f_1g_1^x$ and
$t=f_2g_2^x$ where $f_1g_1, f_2g_2$ are linear combinations of
members of $S_1$ and $S_2$, respectively. Then, $f_1g_1f_2g_2\in
[S_1S_2]$, and $st=f_1f_2(g_1g_2)^x\in [S_1S_2]_1$. Hence,
$[S_1]_1[S_2]_1 \subset [S_1S_2]_1$. By induction, $[S_1]_n [S_2]_n
\subset [S_1S_2]_{n}$. \qedd

\begin{lemma}\label{lm-wm2}
Let $S_1, S_2$ be two subsets of $\F\{\Y^{\pm}\}$ which satisfy
$a\in S_i$ implies $\sigma(a)\in S_i, i=1,2$. Then
$\sqrt{[S_1S_2]_n} = \sqrt{[S_1]_n\cap [S_2]_n}$ for $n\ge 1$, and
$\sqrt{\langle S_1 \rangle} \cap \sqrt{\langle S_2 \rangle} =
\sqrt{\langle S_1S_2 \rangle}$.
\end{lemma}
\proof  The last statement is an immediate consequence of the first one. Since $[S_1S_2] \subset [S_i]$ , we have $[S_1S_2]_n \subset [S_i]_n$ for $i=1,2$, and $[S_1S_2]_n \subset [S_1]_n\cap [S_2]_n$ follows. Hence, $\sqrt{[S_1S_2]_n} \subset \sqrt{[S_1]_n\cap [S_2]_n}$. Let $a\in [S_1]_n\cap [S_2]_n$ we have
$a^2 \in [S_1]_n[S_2]_n$. By Lemma~\ref{lm-wm1}, $a^2 \in [S_1S_2]_{n} $. Hence $a\in \sqrt{[S_1S_2]_n}$, and $ \sqrt{[S_1]_n\cap [S_2]_n} \subset \sqrt{[S_1S_2]_n}$ follows.
\qedd

\begin{lemma}\label{lm-wm3}
Let $\I_1,\ldots,\I_m$ be Laurent $\sigma$-ideals. Then
$\sqrt{\langle \cap_{i=1}^m \I_i \rangle} = \cap_{i=1}^m
\sqrt{\langle \I_i \rangle} $.
\end{lemma}
\proof  Let $\I = \cap_{i=1}^m \I_i$. Then $\sqrt{\I} = \sqrt{[\prod_{i=1}^m\I_i]}$.
By Lemma~\ref{lm-wm2}, we have $\sqrt{\langle \prod_{i=1}^m\I_i \rangle } = \sqrt{\prod_{i=1}^{n-1} \langle \I_i \rangle} \cap \sqrt{\langle \I_n \rangle} = \ldots = \cap_{i=1}^m \sqrt{\langle \I_i \rangle}$. Now we show that $\sqrt{ \langle \I \rangle } = \sqrt{\langle \prod_{i=1}^m\I_i \rangle  }$. Since $\prod_{i=1}^m\I_i \subset \I$, we have $\sqrt{\langle \prod_{i=1}^m\I_i \rangle} \subset \sqrt{ \langle \I \rangle }$. By Lemma~\ref{lm-wm2},  $\sqrt{\langle I \rangle} = \sqrt{\langle I \rangle} \cap \cdots\cap \sqrt{\langle I \rangle} =  \sqrt{\langle I^m \rangle} \subset \sqrt{ \langle \prod_{i=1}^m\I_i \rangle }$, and hence $\sqrt{ \langle \I \rangle } = \sqrt{\langle \prod_{i=1}^m\I_i \rangle  }$. Then, $\sqrt{ \langle \I \rangle } =  \cap_{i=1}^m \sqrt{\langle \I_i \rangle}$.
\qedd

Now, we prove a basic property for a $\sigma$-field $\F$.
\begin{lemma}\label{lm-per1}
Let $\zeta_m = e^{\frac{2\pi \ix}{m}}$ be the primitive $m$-th root
of unity, where $\ix=\sqrt{-1}$ and $m\in\Z_{\ge2}$. If $\F$ is an
algebraically closed $\sigma$-field, then there exists an
$o_m\in[0,m-1]$ such that $\gcd(o_m,m)=1$ and
$\sigma(\zeta_m)=\zeta_m^{o_m}$. Furthermore, the perfect
$\sigma$-ideal $\{y^m-1\}$ in $\F\{y\}$ is
 \begin{equation}\label{eq-per1}
 \{y^m-1\}=[y^m-1,y^{x}-y^{o_m}]\end{equation}
where $y$ is a $\sigma$-indeterminate.
\end{lemma}
\proof Since  $\F$ is algebraically closed, $\zeta_m$ is in $\F$.
From $y^m-1=\prod_{j=0}^{m-1} (y-\zeta_m^j)=0$, we have
$\sigma(y)^m -1=\prod_{j=0}^{m-1} (\sigma(y)-\zeta_m^j)=0$.
Then, there exists an $o_m$ such that $0\le o_m\le m-1$ and
$\sigma(\zeta_m) = \zeta_m^{o_m}$.
Suppose $\gcd(o_m,m)=d>1$ and let $o_m=dk, m=ds$, where
$s\in[1,m-1]$. Then
$\sigma(\zeta_m^s)=\zeta_m^{o_ms}=\zeta_m^{dks}=\zeta_m^{km}=1$,
which implies $\zeta_m^s=1$, a contradiction.

By the difference Nullstellensatz \cite[p.87]{cohn}, we have
$\{y^m-1\}=\cap_{j=0}^{m-1} [y-\zeta_m^j]$. In order to show
\bref{eq-per1}, it suffices to show
 $\cap_{j=0}^{m-1} [y-\zeta_m^j]=[y^m-1,y^x-y^{o_m}]$.
 Since $y^{x}-y^{o_m} = (y- \zeta_m^j)^x +
\zeta_m^{xj} - y^{o_m} = (y- \zeta_m^j)^x + \zeta_m^{jo_m} - y^{o_m}
\in[y -\zeta_m^j]$ for any $0\le j\le m-1$, we have $y^{x}-y^{o_m}
\in \cap_{j=0}^{m-1} [y- \zeta_m^j]$ and hence
$[y^m-1,y^x-y^{o_m}]\subset \cap_{j=0}^{m-1} [y-\zeta_m^j]$.
Let $f\in\cap_{j=0}^{m-1} [y-\zeta_m^j]$. Since
$y^{x}-y^{o_m}\in[y-\zeta_m^j]$,  for $j=0,\ldots,m-1$, from
$f\in[y-\zeta_m^j]$, we have $f= g_j(y-\zeta_m^j)+\sum_k
h_{jk}(y^{x}-y^{o_m})^{x^k}$, where $g_j,h_{jk}\in\Q\{y\}$.
Then $f^m = \prod_{j=0}^{m-1} (g_j(y-\zeta_m^j)+\sum_k
h_{jk}(y^{x}-y^{o_m})^{kx}) = \prod_{j=0}^{m-1}g_j(y^m-1)+p$, where
$p\in[y^{x}-y^{o_m}]$. Hence, $f\in[y^m-1,y^x-y^{o_m}]$ and
$\cap_{j=0}^{m-1} [y-\zeta_m^j]\subset[y^m-1,y^x-y^{o_m}]$. The
lemma is proved.\qedd

The number $o_m$ introduced in Lemma \ref{lm-per1} depends on $\F$
only and is called the {\em $n$-th transforming degree} of unity. In
the following corollaries, $\F$ is assumed to be algebraically
closed and hence $o_m$ is fixed for any $m\in\N$. From the proof of
Lemma \ref{lm-per1}, we have

\begin{cor}\label{cor-per1}
$y^{x}-y^{o_m}\in \cap_{j=0}^{m-1} [y-\zeta_m^j]$.
\end{cor}

\begin{cor}\label{cor-per2}
For $n,m,k$ in $\N$, if $n=km$ then $o_n = o_m\,\mod\, m$.
\end{cor}
\proof By definition, $\zeta_n^{k} =\zeta_m$. Then,
$\sigma(\zeta_n^{k})=\zeta_n^{ko_n}= \zeta_m^{o_n}$. From,
$\sigma(\zeta_n^{k})=\sigma(\zeta_m)=\zeta_m^{o_m}$, we have
$\zeta_m^{o_n}=\zeta_m^{o_m}$. Then $o_n=o_m\,\mod\, m$.\qedd

\begin{lemma}\label{lm-per0}
$\langle y^m-1\rangle=\{y^m-1\}=[y^m-1,y^{x}-y^{o_m}]$.
\end{lemma}
\proof By Lemma~\ref{lm-per1}, it suffices to show
$y^{x}-y^{o_m} \in \langle y^m-1\rangle$. Since
$y^m-1 = \prod_{j=0}^{m-1} (y-\zeta_m^j)$ and $(y-\zeta_m^{i})^x=(y^x-\zeta_m^{o_m i})$, we have
$f_i = (y^x-\zeta_m^{o_m i})\prod_{0\le j \le m-1, j\ne i} (y-\zeta_m^j) \in \langle y^m-1\rangle$
for $i=0, \ldots, m-1$.
We will show that $y^{x}-y^{o_m} \in (f_0, \ldots, f_{m-1})$.
To show this, we need the formula  $\frac{1}{y^m-1} = \sum_{i=0}^{m-1} \frac{1}{m(\zeta_m^i)^{m-1}(y-\zeta_m^i)}  = \frac{1}{m}\sum_{i=0}^{m-1}\frac{\zeta_m^{i}}{y-\zeta_m^i}$ from \cite[p. 494]{GCG1992}.
We have \begin{eqnarray*}
\frac{1}{m}\sum_{i=0}^{m-1} \zeta_m^{i}f_i  & =& \frac{1}{m}\sum_{i=0}^{m-1} \zeta_m^{i}\frac{y^m-1}{y-\zeta_m^i} (y^x-\zeta_m^{o_m i}) \\
& =& \frac{1}{m}\sum_{i=0}^{m-1} \zeta_m^{i}\frac{y^m-1}{y-\zeta_m^i} y^x - \frac{1}{m}\sum_{i=0}^{m-1} \zeta_m^{i}\frac{y^m-1}{y-\zeta_m^i}\zeta_m^{o_m i} \\
& =& y^x - \frac{1}{m}\sum_{i=0}^{m-1} \frac{y^m-1}{y-\zeta_m^i}\zeta_m^{(o_m+1)i}.
\end{eqnarray*}
Let $g(y) = \frac{1}{m}\sum_{i=0}^{m-1} \frac{y^m-1}{y-\zeta_m^i}\zeta_m^{(o_m+1)i}$.
Then, $g(\zeta_m^j) = \frac{1}{m}\zeta_m^{(o_m+1)j} \frac{y^m-1}{y-\zeta_m^j}|_{y=\zeta_m^j} $ $=  \frac{1}{m}\zeta_m^{(o_m+1)j} $\,\,\,\, $\prod_{0\le i\le m-1, i\ne j}$ $(\zeta_m^i-\zeta_m^j) = \frac{1}{m}\zeta_m^{(o_m+1)j}\zeta_m^{j(m-1)}\prod_{1\le i\le m-1} (\zeta_m^i -1) = \frac{1}{m}\zeta_m^{o_m j}m =(\zeta_m^j)^{o_m}$. Since $\deg(g(y)) \le m-1$ and $g(\zeta_m^j) = (\zeta_m^j)^{o_m}$ for $j=0,\ldots,m-1$, we have $g(y) = y^{o_m}$.
Hence $y^x-y^{o_m} \in (f_0, \ldots, f_{m-1}) \subset \langle y^m-1\rangle$.
\qedd

\begin{cor}\label{cor-per3}
For $m\in\N$, $a\in\F^*$, and $\f\in\Zxn$, we have
$\Y^{(x-o_m)\f}-a^{x-o_m} \in \langle\Y^{m\f} - a^m \rangle$.
\end{cor}
\proof Let $z = \frac{\Y^{\f}}{a}$ and $\I=[\Y^{m\f} - a^m]$. Then
$z^{m} - 1 \in \I$. By Lemma \ref{lm-per0},  $z^{x-o_m}-1
\in\langle z^m-1\rangle\subset \langle \I \rangle$.
Then $(\frac{\Y^\f}{a})^{x-o_m}-1 \in \langle\I\rangle$ or
$\Y^{(x-o_m)\f}-a^{x-o_m} \in \langle \I \rangle$.\qedd

The following example shows that the generators of a well-mixed or
perfect ideal depend on the difference field $\F$.
\begin{example}\label{ex-per1}
%
Let $\F=\Q(\sqrt{-3})$ and $p=y_1^{3}-1$. Following Lemma
\ref{lm-per1}, if $\sigma(\sqrt{-3})=\sqrt{-3}$, we have $o_3=1$ and
$\langle p\rangle =\{p\}=[p,y_1^x - y_1]$.
If $\sigma(\sqrt{-3})=-\sqrt{-3}$, we have $o_3=2$ and
$\langle p\rangle =\{p\}=[p,y_1^x - y_1^2]$.
If $\F=\Q$ then
 $\{p\}=[y_1-1]\cap[y_1^2+y_1+1,y_1^x-y_1]\cap[y_1^2+y_1+1,y_1^x-y_1^2]$
does not has a set of simple generators.
%
%
%
%
%
%
\end{example}

Motivated by Corollary \ref{cor-per3}, we have the following
definition.
\begin{definition}
If a $\Zx$-lattice $L$ satisfies
  \begin{equation}\label{eq-wmc}
  m\f \in L \Rightarrow (x-o_m)\f\in  L\end{equation}
where $m\in\N$, $\f\in\Zxn$, and  $o_m$ is defined in Lemma
\ref{lm-per1}, then it is called {\em {\rm{M}}-saturated}. For any
$\Zx$-lattice $L$, the smallest {\rm{M}}-saturated $\Zx$-lattice
containing $L$ is called the {\rm{M}}-saturation of $L$ and is
denoted by $\sat_M(L)$.
\end{definition}

The following result gives an effective version for condition
\bref{eq-wmc}.
\begin{lemma}\label{lm-wmc}
A $\Zxn$-lattice $L$ is M-saturated if and only
if the following condition is true:
Let $L_1 = \sat_\Z (L) = (\g_1,\ldots,\g_s)$  such that
$m_i\g_i\in L$ for $ m_i\in\N$. Then $(x-o_{m_i})\g_i\in
L$.
\end{lemma}
\proof We need only to show $(x-o_{m_i})\g_i\in L$ implies
\bref{eq-wmc}.
For any $ m\f \in L$, we have $\f\in L_1$ and hence $\f =
\sum_{i=1}^r q_i\g_i$, where $q_i\in\Zx$. Let
$t=\lcm(m,m_1,\ldots,m_s)$. By Corollary \ref{cor-per2}, we have
$o_t=o_{m_i} +c_i m_i$, where $c_i\in\Z$.
Then $(x-o_t)\f = \sum_{i=1}^r q_i(x-o_t)\g_i = \sum_{i=1}^r
q_i(x-o_{m_i})\g_i - \sum_{i=1}^r q_ic_i m_i\g_i \in L$.
By Corollary \ref{cor-per2}, $o_t = o_m + c m$, where $c\in\Z$. Then
$(x-o_m)\f = (x-o_t)\f +c m\f\in L$.\qedd

We now give a criterion for a Laurent binomial $\sigma$-ideal to be
well-mixed in terms of its support lattice.
\begin{theorem}\label{th-wm}
Let $\rho$ be a partial character and  $\F$ an algebraically closed
and inversive $\sigma$-field. If $\I(\rho)$ is well-mixed, then
$L_\rho$ is {\rm{M}}-saturated. Conversely, if $L_\rho$ is
{\rm{M}}-saturated, then either $\langle \I(\rho)\rangle=[1]$ or $\I(\rho)$ is
well-mixed.
\end{theorem}
\proof Suppose that $\I(\rho)$ is well-mixed.
If there exists an $m\in\N$ such that $m\f \in
L_\rho$, then by Lemma \ref{lm-l3}, there exists a $c\in\F^*$ such
that $\Y^{m\f} - c\in\I(\rho)$. Since $\F$ is algebraically closed,
there exists an $a\in\F^*$ such that $c = a^m$. Then, $\Y^{m\f} -
a^m\in\I(\rho)$.
Since $\I(\rho)$ is well-mixed, by Corollary \ref{cor-per3},
$\Y^{(x-o_m)\f}-a^{x-o_m} \in \I(\rho)$, and by Lemma \ref{lm-l3}
again, $(x-o_m)\f\in L_\rho$ follows and $L_\rho$ is {\rm{M}}-saturated.

Conversely, let $L_\rho$ be {\rm{M}}-saturated. If $L_\rho$ is
$\Z$-saturated, then by Theorem \ref{th-pr1}, $\I(\rho)$ is
prime and hence well-mixed by Lemma~\ref{lm-wm}.
Otherwise, there exists an
$m_1\in \N$, and $\f\in \Z[x]^n$ such that $\f\not\in L_\rho$ and
$m_1\f\in L_\rho$.
By Lemma \ref{lm-l3}, there exists an $a\in\F^*$ such that
$\Y^{m_1\f} - a^{m_1}\in\I(\rho)$.  We claim that either
$\langle \I(\rho) \rangle=[1]$ or
 \begin{equation}\label{eq-wmt1}\I(\rho) = \cap_{l_1=0}^{m_1-1} \I_{l_1}\end{equation}
where $\I_{l_1} = [\I(\rho), \Y^{\f}-a\zeta_{m_1}^{l_1}]$ and
$\zeta_{m_1} = e^{\frac{2\pi \ix}{m_1}}$.
By \bref{eq-wmc}, $(x-o_{m_1})\f\in L_\rho$. By Lemma \ref{lm-l3},
there exists a $b\in\F^*$ such that $\Y^{(x-o_{m_1})\f} -
b\in\I(\rho)$.
Since $\Y^{m_1\f} - a^{m_1}\in\I(\rho)$, by Corollary
\ref{cor-per1}, we have $\Y^{(x-o_{m_1})\f} -
a^{x-o_{m_1}}\in[\Y^{\f}-a\zeta_{m_1}^{l_1}]$ for any $l_1$. Then
$b-a^{x-o_{m_1}}= \Y^{(x-o_{m_1})\f} - a^{x-o_{m_1}} -
(\Y^{(x-o_{m_1})\f} - b) \in\I_{l_1}$ for any $l_1$.
If $b\ne a^{x-o_{m_1}}$, $\I_{l_1}=[1]$ for all $l_1$, and hence $1\in
\cap_{l_1=0}^{m_1-1} \I_{l_1} \subset \langle \I(\rho) \rangle $ by Lemma~\ref{lm-per0}
and $\langle \I(\rho) \rangle=[1]$ follows.
Now suppose $b=a^{x-o_{m_1}}$ or $a^{x}=ba^{o_{m_1}}$. To prove \bref{eq-wmt1},
it suffices to show  $\cap_{l_1=0}^{m_1-1} \I_{l_1} \subset \I(\rho)$.
Let $f\in \cap_{l_1=0}^{{m_1}-1} \I_{l_1}$. From $f\in\I_{l_1}$, we
have $f = f_{l_1} + \sum_{j=0}^s
p_{j}\sigma^j(\Y^{\f}-a\zeta_{m_1}^{l_1})$, where $f_{l_1}\in\I(\rho)$.
By Lemma \ref{lm-per1}, $\sigma(\zeta_{m_1})=\zeta_{m_1}^{o_{m_1}}$.
We thus have
\begin{eqnarray*}
\sigma(\Y^{\f}-a\zeta_{m_1}^{l_1})
 &=& \Y^{x\f} - b\Y^{o_{m_1}\f} + b\Y^{o_{m_1}\f}- \sigma(a\zeta_{m_1}^{l_1})\\
 &=& \Y^{o_{m_1}\f}(\Y^{(x-o_{m_1})\f}-b) + b(\Y^{o_{m_1}\f}- a^{o_{m_1}} \zeta_{m_1}^{l_1o_{m_1}})
 + (ba^{o_{m_1}}  -\sigma(a))\zeta_{m_1}^{l_1o_{m_1}}.
\end{eqnarray*}
Since  $\Y^{(x-o_{m_1})\f}-b\in\I(\rho)$ and
 $ba^{o_{m_1}}-\sigma(a)=ba^{o_{m_1}}-a^x=0$, we have $\sigma(\Y^{\f}-a\zeta_{m_1}^{l_1}) =
g_{l_1} + q_{l_1}(\Y^{\f}- a\zeta_{m_1}^{l_1})$, where
$g_{l_1}\in\I(\rho)$. Using the above equation repeatedly, we have
$f = h_{l_1} + p_{l_1}(\Y^{\f}- a\zeta_{m_1}^{l_1})$, where
$h_{l_1}\in\I(\rho)$.
Then, $f^{m_1}= \prod_{l_1=0}^{{m_1}-1}(h_{l_1}+p_{l_1}(\Y^{\f}-
a\zeta_{m_1}^{l_1})) = s +
\prod_{l_1=0}^{{m_1}-1}p_{l_1}(\Y^{\f}-a\zeta_{m_1}^{l_1}) = s +
(\Y^{{m_1}\f}-a^{m_1})\prod_{l_1=0}^{{m_1}-1}p_{l_1}\in \I(\rho)$,
where $s$ is in $\I(\rho)$. By Corollary~\ref{cor-rad}, we have
$f\in \I(\rho)$. The claim is proved.

The support lattice for any of $[\I_{l_1}]$ is $L_1 =(L_\rho,\f)$.
Similar to the proof of Theorem \ref{lm-ld0}, we can show that
$\I(\rho)\varsubsetneq\I_{l_1}$ and $L_\rho\varsubsetneq L_1$.
If $L_1$ is not $\Z$-saturated, there exists a $k>1$ and $\g\in\Zxn$
such that $\g\not\in\L_1$ and $k\g\in L_1$. Let $m_2=k{m_1}$. We
have $m_2\g=k{m_1}\g\in L_\rho$ and there exists a $c\in\F^*$
 such that $\Y^{m_2\g}-c^{m_2}\in\I(\rho)$.
Hence, $(x-o_{m_2})\g \in L_\rho\subset L_1$ and there exists a
$d\in\F^*$, such that $\Y^{(x-o_{m_2})\g}-d\in I(\rho)$.
Let $L_2 = (L_1, \g)$ and
$\I_{l_1,l_2}=[\I_{l_1},\Y^{\g}-c\zeta_{m_2}^{l_2}]$,
$l_2=0,\ldots,m_2-1$. Then $L_1\varsubsetneq L_2$ and $L_2$ is the
support lattice for all $\I_{l_1,l_2}$ provided
$\I_{l_1,l_2}\ne[1]$. Similar to the above, it can be shown that
$d-c^{x-o_{m_2}} \in \I_{l_1,l_2}$ for any $l_1,l_2$.
If $d-c^{x-o_{m_2}}\ne0$, then $\I_{l_1,l_2}=[1]$ for any $l_1,l_2$
and $\langle \I_{l_1}\rangle =[1]$ by Lemma~\ref{lm-per0}.
Since Laurent binomial $\sigma$-ideals are radical, $\langle
\I(\rho)\rangle = \cap_{l_1=0}^{m_1-1} \langle \I_{l_1}\rangle =[1]$
by Lemma~\ref{lm-wm3} and \bref{eq-wmt1}.
If $d-c^{x-o_{m_2}}=0$, it can be similarly proved that
 $\I_{l_1} = \cap_{l_2=0}^{m_2-1}\I_{l_1,l_2}$ for any $l_1$.
%
%
As a consequence, we have either $\langle I(\rho)\rangle =[1]$ or
$\I(\rho)=\cap_{l_1=0}^{m_1-1} \I_{l_1} =
\cap_{l_1=0}^{m_1-1}\cap_{l_2=0}^{m_2-1} \I_{l_1,l_2}.$

Repeating the process, we have either $\langle I(\rho)\rangle =[1]$ or
 $$\I(\rho)=\cap_{l_1=0}^{m_1-1} \I_{l_1}= \cdots = \cap_{l_1=0}^{m_1-1}\cdots\cap_{l_t=0}^{m_t-1} \I_{l_1,\ldots,l_t}$$
where $L_\rho\varsubsetneq L_1\varsubsetneq\cdots\varsubsetneq L_t \subset sat_\Z(L_\rho)$.
 Since $\Zxn$ is Notherian, the procedure will terminate and
$L_t$ is  $\Z$-saturated. Since each $\I_{l_1,\ldots,l_t}$ is
either $[1]$ or a prime $\sigma$-ideal, and hence
either $\langle \I(\rho)\rangle=[1]$ or $ \I(\rho)$ is well-mixed by Lemma~\ref{lm-wm}. \qedd

The following example shows that $\langle\I(\rho)\rangle=[1]$ can
indeed happen in Theorem \ref{th-wm}.
\begin{example}\label{ex-311}
Let $\I=[\A]$,  where $\A = \{y_1^2+1, y_1^x-y_1, y_2^2+1,
y_2^x+y_2\}$ is a $\sigma$-chain. Notice that the support lattice of
$\I$ is M-saturated.  We have $y_2^2-y_1^2=y_2^2+1-(y_1^2+1)\in\I$.
Then by Corollary \ref{cor-per3},
$y_1y_2^x-y_1^xy_2\in\langle\I\rangle$. From $y_1^x-y_1,
y_2^x+y_2\in\I$, we have $y_1y_2\in\langle\I\rangle$ and hence $1\in
\langle\I\rangle$.
\end{example}

The following example shows that Theorem \ref{th-per1} is not valid
if condition \bref{eq-wmc} is replaced by
 $m\f \in L_\rho \Rightarrow (x-k)\f\in  L_\rho$ for some $k\in\N$.
\begin{example}\label{eq-per2}
Let $\F=\Q(\sqrt{-3})$ and $\sigma(\sqrt{-3})=-\sqrt{-3}$. Then
$o_3=2$. Let $p_1=y_1^3-1, p_2=y_1^{x-1}-1$. Then $p_1,p_2$ consist
of a Laurent regular and coherent $\sigma$-chain, and by Theorem
\ref{th-t4}, $\I=\sat(p_1,p_2)=[p_1,p_2]$ is proper and with support
lattice $([3],[x-1])$.
By Example \ref{ex-per1}, $y_1^{x-2}-1\in\langle\I\rangle$ and then
$y_1-1 = -y_1^2((y_1^{x-2}-1)-(y_1^{x-1}-1))/y^x\in\langle\I\rangle$
and thus $\langle\I\rangle=[y_1-1]$.
\end{example}

We now prove the main result of this section.
\begin{theorem}\label{th-wm2}
Let $\F$ be an algebraically closed and inversive $\sigma$-field and
$\I=\I(\rho)$ a Laurent binomial $\sigma$-ideal. Then the well-mixed
closure of $\I$ is either $[1]$ or a Laurent binomial $\sigma$-ideal
whose support lattice is $\sat_M(L_\rho)$.
\end{theorem}
\proof
Let us suppose that $\langle\I(\rho)\rangle\ne[1]$.
If $L$ is not M-saturated, then there exists an $m\in\N$ and
$\f\in\Zxn$ such that $\f\not\in L$, $m\f\in L$, and
$(x-o_m)\f\not\in L$.
By Lemma \ref{lm-l3}, there exists a $c\in\F^*$ such that $\Y^{m\f}
- c^m\in\I(\rho)$.
Let $\I_1 =[\I, \Y^{(x-o_m)\f}-c^{x-o_m}]$ and $L_1 = (L,(x-o_m)\f)$. By
Corollary \ref{cor-per3}, $\Y^{(x-o_m)\f}-c^{x-o_m}\in\langle \I(\rho)\rangle$.
Let $L_M=\sat_M(L)$. Then $\I\varsubsetneq \I_1\subset\langle\I\rangle$
and $L\varsubsetneq L_1\subset L_M$.
Repeat the procedure to construct $I_i$ and $L_i$ for $i=2,\ldots,t$
such that $\I\varsubsetneq
\I_1\varsubsetneq\cdots\varsubsetneq\I_t \subset\{\I\}$ and
$L\varsubsetneq L_1\varsubsetneq\cdots\varsubsetneq L_t\subset
L_M$. Since $\Zxn$ is Notherian, the procedure will terminate at,
say $t$. Then $L_t=L_M$ is M-saturated. By Lemma \ref{lm-pxsat},
$L_t$ is also $x$-saturated. By Theorem \ref{th-wm},
$\I_t\subset\langle\I\rangle$ is well-mixed and hence $\I_t=\langle\I\rangle$. \qedd

By Corollary \ref{cor-ld0} and the proof of Theorem \ref{th-wm2},
we have
\begin{cor}
A $\Zx$-lattice and its M-saturation have the same rank. Hence, for
a Laurent binomial $\sigma$-ideal $\I$, either $\langle \I
\rangle=[1]$ or $\dim(\langle\I\rangle) = \dim(\I)$.
\end{cor}

\begin{example}\label{ex-perc1}
 Let $p=y_2^{2}-y_1^{2}$. Following the proof of Theorem \ref{th-wm2},
 it can be shown that
 $\langle p\rangle=\{p\}=[y_1^{-2}y_2^{2}-1,y_1^{1-x}y_2^{x-1}-1]=[y_2^{2}-y_1^{2},y_1y_2^{x}-y_1^{x}y_2]$.
\end{example}

\vskip10pt
In the rest of this section, we prove similar results for
the perfect closure of Laurent binomial $\sigma$-ideals.
We first give a definition.

\begin{definition}
If a $\Zx$-lattice is both $x$-saturated and M-saturated, then it is
called {\em {\rm{P}}-saturated}. For any  $\Zx$-lattice $L$, the
smallest {\rm{P}}-saturated $\Zx$-lattice containing $L$ is called
the {\rm{P}}-saturation of $L$ and is denoted by $\sat_P(L)$.
\end{definition}

\begin{lemma}\label{lm-pxsat}
For any  $\Zx$-lattice $L$, $\sat_P(L)=\sat_x(\sat_M(L))=\sat_M(\sat_x(L))$.
\end{lemma}
\proof Let $L_1=\sat_x(\sat_M(L))$ and $L_2 = \sat_M(\sat_x(L))$.
It suffices to show $L_1=L_2$.
We claim that $L_1$ is P-saturated. Let
$m\f\in L_1$ for $m\in\N$. Then $mx^a\f\in\sat_M(L)$ for some
$a\in\N$, which implies $(x-o_m)x^a\f\in L\subset
\sat_x(\sat_M(L))=L_1$. Since $L_1$ is $x$-saturated, $(x-o_m)\f\in
L_1$ and the claim is proved.
Since $L\subset \sat_M(L)$, $\sat_x(L)\subset\sat_x(\sat_M(L))=L_1$.
From the claim, $L_1$ is P-saturated and hence
$L_2\subset\sat_M(L_1)=L_1$.

For the other direction, we claim that  $L_2$ is
$x$-saturated. Let $x\f\in
\sat_M(\sat_x(L))\subset\sat_\Z(\sat_x(L))$. Then there exists an
$m\in\N$, such that $m\f\in \sat_x(L)$ which implies $(x-o_m)\f\in
\sat_M(\sat_x(L))$ and hence $o_m\f=x\f-(x-o_m)\f\in
\sat_M(\sat_x(L))$ follows. By Lemma \ref{lm-per1}, $\gcd(o_m,m)=1$.
Then $\f \in \sat_M(\sat_x(L))$, and the claim is true. Since
$\sat_M(L)\subset \sat_M(\sat_x(L)) = L_2 = \sat_x(\sat_M(\sat_x(L))
)$, we have $L_1\subset L_2$. \qedd

A $\sigma$-ideal $\I$ is perfect if and only if $\I$ is reflexive,
radical, and well-mixed. Since a Laurent binomial $\sigma$-ideal
$\I$ is always radical, $\I$ is perfect if and only if $\I$ is
reflexive and well-mixed. Following from this observation, we can
deduce the following results about perfect Laurent binomial
$\sigma$-ideal ideals.

\begin{theorem}\label{th-per1}
Let $\rho$ be a partial character and $\F$ an algebraically closed
and inversive $\sigma$-field. If $\I(\rho)$ is perfect, then
$L_\rho$ is {\rm{P}}-saturated. Conversely, if $L_\rho$ is
{\rm{P}}-saturated, then either $\{\I(\rho)\}=[1]$ or $\I(\rho)$ is
perfect.
\end{theorem}
\proof If $\I(\rho)$ is perfect, then it is well-mixed and
reflexive. By Theorems~\ref{th-wm} and Theorem~\ref{th-pr1},
$L_\rho$ is {\rm{M}}-saturated and $x$-saturated, and hence
{\rm{P}}-saturated.
Conversely, if $L_\rho$ is {\rm{P}}-saturated, it is
{\rm{M}}-saturated and $x$-saturated. By Theorem~\ref{th-wm},
 either $\langle \I(\rho)\rangle=[1]$ or $\I(\rho)$ is well-mixed.
If $\langle \I(\rho)\rangle=[1]$, $\{\I(\rho)\}=[1]$. Otherwise, by
Theorem~\ref{th-pr1}, $\I(\rho)$ is reflexive. By
Corollary~\ref{cor-rad}, $\I(\rho)$ is radical. Then $\I(\rho)$ is
perfect.\qedd

\begin{theorem}\label{th-per2}
Let $\F$ be an algebraically closed and inversive $\sigma$-field and
$\I=\I(\rho)$ a Laurent binomial $\sigma$-ideal. Then the perfect
closure of $\I$ is either $[1]$ or a Laurent binomial $\sigma$-ideal
whose support lattice is $\sat_P(L_\rho)$.
\end{theorem}
\proof Let $\I_x$ be the reflexive closure of $\I$ and
$L_x=\sat_x(L_\rho)$. By Theorem~\ref{lm-ld0}, $L_x$ is the support
lattice of $\I_x$ and $\I_x$ is reflexive.
Let $\I_p$ be the well-mixed closure of $\I_x$ and
$L_p=\sat_M(\sat_x(L_\rho))$. By Theorem~\ref{th-wm2}, $\I_p$ is
either $[1]$ or well-mixed, and in the latter case $L_p$ is the
support lattice of $\I_p$. By Lemma \ref{lm-pxsat},
$L_p=\sat_P(L_\rho)$. Then $\I_p$ is either $[1]$ or perfect by
Theorem \ref{th-per1}. \qedd

With Theorem \ref{th-wm2} and Theorem \ref{th-per1}, (3) and (4) of
Theorem \ref{th-m2} are proved.

\section{Binomial $\sigma$-ideal}\label{sec-bi}
In this section, properties of binomial $\sigma$-ideals will be
proved. First, certain results from \cite{es-bi} will be extended to
the difference case using the theory of infinite Gr\"obner basis
\cite{Kei}. Then, properties proved in Section \ref{sec-lbi} will be
extended to so-called normal binomial $\sigma$-ideals.

\subsection{Basic properties of binomial $\sigma$-ideal}

A {\em $\sigma$-binomial}   in  $\Y$ is  a $\sigma$-polynomial with
at most two terms,  that is, $a\Y^{\a}+b\Y^{\b}$ where $a, b\in \F$
and $\a, \b \in\N[x]^n$.
For $\f\in\Z[x]^n$,  let $\f^{+}, \f^{-}\in \N^{n}[x]$ denote the
positive part and the negative part of $\f$ such that
$\f=\f^{+}-\f^{-}$.
Consider a $\sigma$-binomial $f=a\Y^{\a}+b\Y^{\b}$, where
$a,b\in\F^*$. Without loss of generality, assume $\a > \b$ according
to the order defined in Section \ref{subsec-Gro}. Then $f$ has the
following canonical representation
 \begin{equation}\label{eq-nfb}
 f=a\Y^{\a}+b\Y^{\b} = a\Y^\g(\Y^{\f+}-c\Y^{\f^-})
 \end{equation}
where $c=\frac{-b}{a}$, $\f=\a-\b\in\Zxn$  is a normal vector, and
$\g=\a-\f^+\in\N[x]$. The normal vector $\f$ is called the {\em
support} of $f$. Note that $\gcd(\Y^{\f+},\Y^{\f^-})=1$.

A $\sigma$-ideal in $\F\{\Y\}$ is called {\em binomial} if it is
generated by, possibly infinitely many, $\sigma$-binomials.

In the following, $\F\{\Y\}$ is considered as a polynomial ring in
infinitely many variables $\Theta(\Y)=\{y_i^{x^j}, i=1,\ldots,n;
j\ge 0\}$ and denoted by $S=\F[\Theta(\Y)]$. A theory of Gr\"obner
basis in the case of infinitely many variables is developed in
\cite{Kei} and will be used in this section.
For any $m\in\N$, denote $\Theta^{\lr{m}}(\Y)=\{y_i^{x^j},
i=1,\ldots,n; j=0,1,\ldots,m\}$ and
$S^{\lr{m}}=\F[\Theta^{\lr{m}}(\Y)]$.

A monomial order in  $S$ is called {\em compatible} with the
difference structure, if $y_i^{x^{k_1}}<y_i^{x^{k_2}}$ for $k_1 <
k_2$.
%
%
Only compatible monomial orders are considered in this section.

Let $\I$ be a  $\sigma$-ideal in $\F\{\Y\}$. Then $\I$ is an
algebraic ideal in $S$. By \cite{Kei}, we have
\begin{lemma}\label{lm-gbbi1}
Let $\I$ be a binomial $\sigma$-ideal in $\F\{\Y\}$. Then for a
compatible monomial order, the reduced Gr\"{o}bner basis $\G$ of
$\I$ exists and satisfies
 \begin{equation}\label{eq-gbbi1} \G = \cup_{m=0}^{\infty} \G^{\lr{m}}\end{equation}
where $\G^{\lr{m}}=\G\cap S^{\lr{m}}$ is the reduced Gr\"{o}bner
basis of $\I^{\lr{m}} = \I \cap S^{\lr{m}}$ in the polynomial ring
$S^{\lr{m}}$.
\end{lemma}

Contrary to the Laurent case, a binomial $\sigma$-ideal may be
infinitely generated, as shown by the following example.
\begin{example}\label{ex-id}
Let  $\I=[y_1^{x^{i}}y_2^{x^{j}}-y_1^{x^{j}}y_2^{x^{i}} : 0 \leq i <
j \in \N]$. It is clear that $\I$ does not have a finite set of
generators and hence a finite Gr\"obner basis.
%
%
The Gr\"obner basis of
 $$\I^{\lr{m}}=\I\cap \Q[y_1,
y_2;y^{x}_1,y^{x}_2;\ldots;y^{x^{m}}_1,y^{x^{m}}_2]$$ is
$\{y_1^{x^{i}}y_2^{x^{j}}-y_1^{x^{i}}y_2^{x^{i}} : 0 \leq i < j \leq
m\}$ with a monomial order satisfying $y_1<
y_2<y^{x}_1<y^{x}_2<\cdots<y^{x^m}_1<y^{x^m}_2$. Then
 $\{y_1^{x^{i}}y_2^{x^{j}}-y_1^{x^{j}}y_2^{x^{i}}
: 0 \leq i < j \in \N\}$ is an infinite reduced Gr\"obner basis for
$\I$ in the sense of \cite{Kei} when $y^{x^m}_1$ and $y^{x^m}_2$ are
treated as independent variables.
\end{example}
\begin{remark}\label{rem-bi11}
The above concept of Gr\"obner basis does not consider the
difference structure.
The concept may be refined by introducing the reduced {\em
$\sigma$-Gr\"obner basis} \cite{cohn}. A $\sigma$-monomial $M_1$ is
called {\em reduced} w.r.t. another $\sigma$-monomial $M_2$ if there
do not exist a $\sigma$-monomial $M_0$ and a $k\in\N$ such that
$M_1=M_0 M_1^{x^k}$.
Then the reduced $\sigma$-Gr\"obner basis of $\I$ in Example
\ref{ex-id} is $\{y_1y_2^{x^{i}}-y_1^{x^{i}}y_2: i \in \Z_{\ge1}\}$
which is still infinite.
Since the purpose of Gr\"obner basis in this paper is theoretic and
not computational, we will use the version of infinite Gr\"obner
basis in the sense of \cite{Kei}.
\end{remark}

With Lemma \ref{lm-gbbi1}, a large portion of the properties for
algebraic binomial ideals proved by Eisenbud and Sturmfels in
\cite{es-bi} can be extended to the difference case.
The proofs follow the same pattern: to prove a property for $\I$, we
first show that the property is valid for $\I$ if and only if it is
valid for all $\I^{\lr{m}}$, and then the corresponding statement
from \cite{es-bi} will be used to show that the property is indeed
valid for $\I^{\lr{m}}$. We will illustrate the procedure in the
following corollary. For other results, we omit the proofs.

\begin{cor}\label{term}
Let $\I \subset \F\{\Y\}$ be a binomial $\sigma$-ideal. Then the
Gr\"{o}bner basis $\G$ of $\I$ consists of $\sigma$-binomials and
the normal form of any $\sigma$-term modulo $\G$ is again a
$\sigma$-term.
\end{cor}
\proof By a $\sigma$-term, we mean the multiplication of an element
from $\F^*$ and a $\sigma$-monomial. By \bref{eq-gbbi1}, it suffices
to show that corollary is valid for all $\G^{\lr{m}}$, that is, the
Gr\"{o}bner basis $\G^{\lr{m}}$ of $\I^{\lr{m}}$ consists of
binomials and the normal form of any term modulo $\G^{\lr{m}}$ is
again a term.
Since $\G^{\lr{m}}$ is the Gr\"{o}bner basis of $\I^{\lr{m}} = \I
\cap S^{\lr{m}}$ and $\I^{\lr{m}}$ is a binomial ideal in a
polynomial ring with finitely many variables, the corollary follows
from Proposition 1.1 in \cite{es-bi}.\qedd
%
%

\begin{cor}
A $\sigma$-ideal $\I$ is binomial if and only if the reduced
Gr\"{o}bner basis for $\I$ consists of $\sigma$-binomials.
\end{cor}

%
\begin{cor}\label{eli}
If $\I$ is a  binomial $\sigma$-ideal, then the elimination ideal
$\I \cap \F\{y_1,y_2,\ldots,y_r\}$ is binomial for every $r\leq n$.
\end{cor}
%
%

Now we consider the intersection of two ideals. The following lemma
can be proved similar to its algebraic counterpart.
\begin{lemma}\label{lm-di11}
If $\I$ and $\J$ are binomial $\sigma$-ideals in $\F\{\Y\}$ then we
have $\I \cap\J=[t\I+(1-t)\J]\cap \F\{\Y\}$ where $t$ is a new
$\sigma$-indeterminate.
\end{lemma}

The intersection of binomial $\sigma$-ideals is not binomial in
general, but from Lemma \ref{lm-di11} and \cite{es-bi} we have
\begin{cor}\label{cor-gbi8}
If  $\I$ and $\I^{'}$ are binomial $\sigma$-ideals and
$\J_{1},\ldots,\J_{s}$ are $\sigma$-ideals generated by
$\sigma$-monomials, then $[\I+\I^{'}]\cap[ \I+\J_{1}]\cap\ldots
\cap[\I+\J_{s}]$ is binomial.
\end{cor}

\begin{cor}\label{cor-tttt1}
Let $\I$ be a binomial $\sigma$-ideal  and let
$\J_{1},\ldots,\J_{s}$ be monomial $\sigma$-ideals.

(a) The intersection $[ \I+\J_{1}]\cap\cdots \cap[\I+\J_{s}]$ is
generated by $\sigma$-monomials modulo $\I$.

(b) Any $\sigma$-monomial in the sum $\I+\J_{1}+\cdots+\J_{s}$ lies
in one of the $\sigma$-ideals $\I+\J_{i}$.

\end{cor}

\begin{cor}\label{diff-quotient}
If $\I$ is a binomial $\sigma$-ideal, then for any $\sigma$-monomial
$M$, the $\sigma$-ideal quotients $[\I:M]$ and $[\I:M^{\infty}]$ are
binomial.
\end{cor}

\begin{cor}
Let $\I$ be a binomial $\sigma$-ideal and $\J$ a monomial
$\sigma$-ideal. If $f\in \I+\J$ and $g$ is the sum of those terms of
$f$ that are not individually contained in $\I+\J$, then $g\in \J$.
\end{cor}

Finally, from \cite[Theorem 3.1]{es-bi}, we have
\begin{theorem}\label{th-bi11}
If $\I$ is a  binomial $\sigma$-ideal, then the radical of $\I$ is
binomial.
\end{theorem}

We now consider whether a $\sigma$-ideal is reflexive. We first give
a criterion for reflexiveness.

\begin{lemma}\label{lm-biref1}
A $\sigma$-ideal $\I$ is reflexive if and only if  $b^x\in\I\imply
b\in\I$ for any $\sigma$-binomial $b\in\F\{\Y\}$.
\end{lemma}
\proof We need only to prove one side of the statement, that is, if
$b^x\in\I\imply b\in\I$ for any $\sigma$-binomial $b$ then $\I$ is
reflexive.
Let $p$ be a $\sigma$-polynomial such that $p^x\in\I$. Then, there
exists an $m\in\N$ such that $p^x\in\I^{\lr{m}} = \I \cap
S^{\lr{m}}$. Let $\G$ be the (finite) reduced Gr\"obner basis of
$\I^{\lr{m}}$ in $S^{\lr{m}}$ under the variable order $y_i^{x^j}<
y_k$ for any $i,j,k$. By Proposition 1.1 in \cite{es-bi}, $\G$
consists of binomials. $p^x$ can be reduced to zero by $\G$. Due to
the chosen variable order, we have $p^x =\sum_i e_i^x g_i^x$, where
$e_i^x\in S^{\lr{m}}$ and $g_i^x$ is a $\sigma$-binomial in
$S^{\lr{m}}$.
Since $g_i^x$ are $\sigma$-binomials in $\I$, we have $g_i\in\I$.
Then, $p =\sum_i e_i g_i\in\I$ and $\I$ is reflexive.\qedd

\begin{theorem}\label{th-bi12}
If $\I$ is a  binomial $\sigma$-ideal, then the reflexive closure of
$\I$ is binomial.
\end{theorem}
\proof Let $\I_1$ be the $\sigma$-ideal generated by the
$\sigma$-binomials $p$ such that $p^{x^{k}} \in\I$ for some
$k\in\N$. We claim that $\I_1$ is the reflexive closure of $\I$ and
it suffices to show that $\I_1$ is reflexive. Let $p$ be a
$\sigma$-binomial such that $p^x\in\I_1$. Then for some $k\in\N$,
$(p^x)^{x^k} = p^{x+1}\in\I$. Thus $p\in\I_1$ and $\I_1$ is
reflexive by Lemma \ref{lm-biref1}.\qedd

\subsection{Normal binomial $\sigma$-ideal}
In this section, most of the results about Laurent binomial
$\sigma$-ideals proved in Section \ref{sec-lbi} will be extended to
normal binomial $\sigma$-ideals.

Let $\mb$ be the multiplicative set generated by $y_i^{x^j}$ for
$i=1, \ldots, n,  j\in\N$.
A $\sigma$-ideal $\I$ is called {\em normal} if for any $M\in\mb$
and $p\in\F\{\Y\}$,  $Mp\in\I$ implies $p\in\I$. For any
$\sigma$-ideal $\I$,
$$\I:\mb=\{f\in \F\{\Y\}\,|\,\exists M\in \mb \st M f\in \I\}$$
is a normal $\sigma$-ideal.
For any $\sigma$-ideal $\I$ in $\F\{\Y\}$, it is easy to check that
 \begin{equation}\label{eq-nbi1}
 \F\{\Y^{\pm}\}\I\cap\F\{\Y\}=\I:\mb.\end{equation}
We first prove a property for general normal $\sigma$-ideals.
\begin{lemma}\label{lm-gni1}
Let $\I$ be a normal $\sigma$-ideal in $\F\{\Y\}$. Then $\I$ is
reflexive (radical, perfect, prime) if and only if
$\F\{\Y^{\pm}\}\I$ is reflexive (radical, perfect, prime) in
$\F\{\Y^{\pm}\}$.
\end{lemma}
\proof Let $\overline{\I}=\F\{\Y^{\pm}\}\I$ be a Laurent
$\sigma$-ideal. Since $\I$ is normal, from \bref{eq-nbi1} we have
$\overline{\I}\cap\F\{\Y\}=\I$.
If $\overline{\I}$ is reflexive, it is clear that $\I$ is reflexive.
For the other direction, if $f^{x}\in \overline{\I}$, then by
clearing denominators of $f^{x}$, there exists a $\sigma$-monomial
$M^{x}$ in $\Y$ such that $M^{x}f^{x}\in\overline{\I} \cap \F\{\Y\}=
\I$. Since $\I$ is reflexive, $Mf\in \I$ and hence $f\in
\overline{\I}$, that is, $\overline{\I}$ is reflexive. The results
about radical and perfect $\sigma$-ideals can be proved similarly.

We now show that $\I$ is prime if and only if $\overline{\I}$ is
prime.
If $\overline{\I}$ is prime, it is clear that  $\I$ is also prime.
For the other side, let $fg\in \overline{\I}$. Then  there exist
$\sigma$-monomials $N_1,N_2$ such that $N_1f\in\F\{\Y\}$,
$N_2g\in\F\{\Y\}$, and hence $N_1f N_2g\in \I$. Since $\I$ is prime,
$N_1f$ or $N_2g$ is in $\I$ that is $f$ or $g$ is in
$\overline{\I}$.\qedd

We now consider normal binomial $\sigma$-ideals.  Given a partial
character $\rho$ on $\Z[x]^{n}$, we define the following binomial
$\sigma$-ideal in $\F\{\Y\}$
 \begin{equation}\label{eq-I+}
 \I^{+}(\rho)=[{\Y^{\f^{+}}-\rho(\f)\Y^{\f^{-}}:\f\in
 L_{\rho}}].\end{equation}
We will show that any normal binomial $\sigma$-ideal can be written
as the form \bref{eq-I+}.
%


\begin{lemma}\label{lm-nl1}
Let $\rho$ be a partial character on $\Zxn$ and $\I(\rho)$ defined
in \bref{eq-I}. Then $\I^{+}(\rho)=\I(\rho)\cap \F\{\Y\}$. As a
consequence, $\I^{+}(\rho)$ is proper and normal.
\end{lemma}
\proof It is clear that  $\I^{+}(\rho)\subset \I(\rho)\cap
\F\{\Y\}$.
If $f\in \I(\rho) \cap \F\{\Y\}$, then $f=\sum_{i=1}^s f_i
M_i(\Y^{\f_i}-\rho(\f_i))$ where $f_i\in\F$, $\f_i\in L_{\rho}$, and
$M_i$ are Laurent $\sigma$-monomials in $\Y$.
There exists a $\sigma$-monomial $M$ in $\Y$ such that
 \begin{equation}\label{eq-nl1}
 Mf=\sum_{i=1}^s f_i N_i(\Y^{\f_i^{+}}-\rho(\f_i)\Y^{\f_i^{-}})\in \I^{+}(\rho),\end{equation}
where $N_i$ is a $\sigma$-monomial in $\Y$. We will prove $f\in
\I^{+}(\rho)$ from the above equation.
Without loss of generality, we may assume that $M=y_c^{x^{o}}$ for
some $c$ and $o\in\N$. Note that \bref{eq-nl1} is an algebraic
identity in $y_i^{x^{k}}, i=1,\ldots,n, k\in\N$.
If $N_i$ contains $y_c^{x^{o}}$ as a factor, we move
$F_i=f_iN_i(\Y^{\f_i^{+}}-\rho(\f_i)\Y^{\f_i^{-}})$ to the left hand
side of \bref{eq-nl1} and let $f_1=f-F_i/y_c^{x^{o}}$. Then
$f\in\I^{+}(\rho)$ if and only if $f_1\in\I^{+}(\rho)$.
Repeat the above procedure until no $N_i$ contains $y_c^{x^{o}}$ as
a factor.

If $s=0$ in \bref{eq-nl1}, then $f=0$ and the lemma is proved.
Since $\gcd(\Y^{\f_i^{+}},\Y^{\f_i^{-}})=1$, $y_c^{x^{o}}$ cannot be
a factor of both $\Y^{\f_i^{+}}$ and $\Y^{\f_i^{-}}$.
Let $\Y^{\f_i^{+}}$ be the largest $\sigma$-monomial in
\bref{eq-nl1} not containing $y_c^{x^{o}}$ under a given
$\sigma$-monomial total order . If $\Y^{\f_i^{-}}$ is the largest
$\sigma$-monomial in \bref{eq-nl1} not containing $y_c^{x^{o}}$, the
proving process is similar.
There must exists another $\sigma$-binomial
$f_jN_j(\Y^{\f_j^{+}}-\rho(\f_j)\Y^{\f_j^{-}})$ such that
$N_i\Y^{\f_i^{+}}=N_j\Y^{\f_j^{-}}$.
Let $N_i=\Y^{\p_i}, N_j=\Y^{\p_j}$. Then
$\Y^{\f_i^{+}+\p_i}=\Y^{\f_j^{-}+\p_j}$ and
$\f_i^{+}+\p_i=\f_j^{-}+\p_j$.
We have $p=f_iN_i (\Y^{\f_i^{+}}-\rho(\f_i)\Y^{\f_i^{-}})+
 f_jN_j(\Y^{\f_j^{+}}-\rho(\f_j)\Y^{\f_j^{-}})
  =\frac{f_i}{\rho(\f_j)}(\Y^{\f_j^{+}+\p_j}-\rho(\f_i)\rho(\f_j)\Y^{\f_i^{-}+\p_i})
 +(f_j-\frac{f_i}{\rho(\f_j)})N_j(\Y^{\f_j^{+}}-\rho(\f_j)\Y^{\f_j^{-}})$.
%
%
Since
$\f=\f_j^{+}+\p_j-(\f_i^{-}+\p_i)=\f_i^{+}-\f_i^{-}+\f_j^{+}-\f_j^{-}=\f_i+\f_j\in
L_{\rho}$, we have
$\Y^{\f_j^{+}+\p_j}-\rho(\f_i)\rho(\f_j)\Y^{\f_i^{-}+\p_i}
 =N(\Y^{\f^+}-\rho(\f)\Y^{\f^{-}})\in \I^{+}(\rho)$, where $N$ is a
 $\sigma$-monomial.
As a consequence, $p\in \I^{+}(\rho)$.
%
%
If $N$ contains $y_c^{x^{o}}$, move the term
$\frac{f_i}{\rho(\f_j)}N(\Y^{\f^{+}}-\rho(\f)\Y^{\f^{-}})$ to the left hand
side of \bref{eq-nl1} as we did in the first phase of the proof.
After the above procedure, equation \bref{eq-nl1} is still valid.
Furthermore, the number of $\sigma$-binomials in \bref{eq-nl1} does
not increase, no $N_i$ contains $y_c^{x^{o}}$, and the largest
$\sigma$-monomial $\Y^{\f_i^{+}}$ or $\Y^{\f_i^{-}}$ not containing
$y_c^{x^{o}}$ becomes smaller.
The above procedure will stop after a finite number of steps, which
means $s=0$ in \bref{eq-nl1} and hence $y_c^{x^{o}}f =0$ which means
the original $f$  is in $\I^{+}(\rho)$. Then
$\I^{+}(\rho)=\I(\rho)\cap \F\{\Y\}$.

$\I^{+}(\rho)=\I(\rho)\cap \F\{\Y\}$ is proper. For otherwise
$\I(\rho)=[1]$,  contradicting to Lemma \ref{lm-l1}.
Note that $\I^{+}(\rho)\F\{\Y^{\pm}\}=\I(\rho)$. Then
$\I^{+}(\rho)=\I(\rho)\cap \F\{\Y\}=\I^{+}(\rho)\F\{\Y^{\pm}\}\cap
\F\{\Y\}=\I^{+}(\rho):\mb$, and $\I^{+}(\rho)$ is normal.\qedd

\begin{lemma}  \label{lm-nl2}
Let $\rho$ be a partial character over $\Zxn$. Then
$\Y^{\f^+}-c\Y^{\f^-} \in \I^{+}(\rho)$ if and only if $\f\in
L_{\rho}$ and $c=\rho(\f)$.
\end{lemma}
\proof By Lemma \ref{lm-nl1}, $\Y^{\f^+}-c\Y^{\f^-} \in
\I^{+}(\rho)$ if and only if $\Y^{\f}-c \in\I(\rho)$ which is
equivalent to $\f\in L_{\rho}$ and $c=\rho(\f)$ by Lemma
\ref{lm-l3}. \qedd

\begin{lemma}  \label{lm-nl3}
If $\I$ is a normal binomial $\sigma$-ideal,  then there exists a
unique partial character $\rho$ on $\Z[x]^{n}$ such that
$\I=\I^{+}(\rho)$ and $L_{\rho}=\{\f\in\Zxn\,|\,
\Y^{\f^+}-\rho(\f)\Y^{\f^-} \in \I\}$ which is called the {\em
support lattice} of $\I$.
\end{lemma}
\proof We have $\I\cdot \F\{\Y^{\pm}\}\cap \F\{\Y\} =\I:\mb$.
According to Theorem \ref{th-l1},
 there exists a partial character $\rho$ such that
$\I\cdot \F\{\Y^{\pm}\}=\I(\rho)$. Then by Lemma \ref{lm-nl1},
$\I=(\I:\mb)=\I\cdot \F\{\Y^{\pm}\}\cap \F\{\Y\}=\I(\rho)\cap
\F\{\Y\}=\I^{+}(\rho)$.
By Lemma \ref{lm-nl2}, we have $L_{\rho}=\{\f\in\Zxn\,|\,
\Y^{\f^+}-\rho(\f)\Y^{\f^-} \in \I=\I^{+}(\rho)\}$. The uniqueness
of $\rho$ comes from the fact that $L_{\rho}$ is uniquely determined
by $\I$.\qedd


By Lemmas \ref{lm-nl1} and  \ref{lm-nl3}, we have
\begin{theorem}\label{th-nbi2}
The map $\I(\rho) \Rightarrow \I^{+}(\rho)$ gives a one to one
correspondence between Laurent binomial $\sigma$-ideals and normal
binomial $\sigma$-ideals.
\end{theorem}

Due to Lemma \ref{lm-nl1} and Theorem \ref{th-nbi2}, most properties
of Laurent binomial $\sigma$-ideals can be extended to normal
binomial $\sigma$-ideals.
As a consequence of Corollary \ref{cor-rad}, Lemma \ref{lm-gni1},
and Lemma \ref{lm-nl1}, we have
\begin{cor}\label{cor-rad1}
A normal binomial $\sigma$-ideal is radical.
\end{cor}

As a consequence of Theorem \ref{lm-ld0}, Lemma \ref{lm-gni1}, and
Theorem \ref{th-nbi2}.
\begin{cor}\label{cor-ld01}
Let $\I=\I^+(\rho)$ be a normal binomial $\sigma$-ideal. If $\F$ is
inversive, then the reflexive closure of $\I$ is also a normal
binomial $\sigma$-ideal whose support lattice is the $x$-saturation
of $L_\rho$.
\end{cor}

%
\begin{cor}\label{th-nl1}
Let $\I=\I+(\rho)$ be a normal binomial $\sigma$-ideal. If $\F$ is
algebraically closed and inversive,  then
\begin{description}
\item[(a)] $L_\rho$ is $\Z$-saturated if and only if $\I$ is prime;
\item[(b)] $L_\rho$ is $x$-saturated if and only if $\I$ is reflexive;
\item[(c)] $L_\rho$ is saturated if and only if $\I$ is reflexive prime.
\end{description}
\end{cor}
\proof It is easy to show that
$\I(\rho)=\I^{+}(\rho)\F\{\Y^{\pm}\}$. Then the corollary is a
consequence of Theorem $\ref{th-pr1}$, Lemma \ref{lm-gni1}, and
Lemma \ref{lm-nl1}.\qedd


For properties related with perfect $\sigma$-ideals, it becomes more
complicated. Direct extension of Theorems \ref{th-l2},
\ref{th-per1}, and \ref{th-per2} to the normal binomial case is not
correct as shown by the following example.
\begin{example}\label{ex-bi-per}
Let $\I = [y_1^x - y_1, y_2^2 - y_1^2, y_2^x + y_2]$ which is a
normal binomial $\sigma$-ideal whose representation matrix is
$L=\left [\begin{array}{ccc}
x-1  & -2    & 0   \\
0    & 2     & x-1 \\
\end{array}\right]$.
Since $o_2=1$, $L$ is $P$-saturated. Also, $L_s = \sat(L) =
 \left [\begin{array}{ccc}
  x-1  & -1 \\
 0     & 1 \\ \end{array}\right]$.
We have $\{\I\} = \{\I,y_2-y_1\}\cap\{\I,y_2+y_1\} = [y_1,y_2]$.
Then $\{\I\}\ne[1]$ and $\I$ is not perfect and hence Theorems
\ref{th-per1} and \ref{th-per2} are not correct. Theorem \ref{th-l2}
is also not correct, since the supporting lattice of the prime
component of $\I$ is not $L_s$.
\end{example}

It can be seen that the problem  is due to the occurrence of
$\sigma$-monomials. For any partial character $\rho$, it can be
shown that
 \begin{equation}\label{eq-per1-1}
  \{\I^{+}(\rho)\}:\mb=\{\I(\rho)\}\cap \F\{\Y\}.\end{equation}
We thus have the following modifications for Theorems \ref{th-per1},
\ref{th-per2}, and \ref{th-l2}.

\begin{cor}\label{th-nl1p}
Let  $\I=\I^+(\rho)$ be a normal binomial $\sigma$-ideal and $\F$ an
inversive and algebraically closed $\sigma$-field.
If $\I$ is perfect, then $L_\rho$ is $P$-saturated. Conversely, if
$L_\rho$ is $P$-saturated and $x$-saturated, then either
$\{\I\}:\mb=[1]$ or $\I$ is perfect.
\end{cor}
\proof If $\I$ is perfect, by Lemma \ref{lm-gni1}, $\I(\rho) =
\I\F\{\Y^{\pm}\}$ is also perfect. By Theorem \ref{th-per1},
$L_\rho$ is $P$-saturated.
If $L_\rho$ is $P$-saturated and $x$-saturated, by Theorem
\ref{th-per1}, either $\I(\rho)=[1]$ or $\I(\rho)$ is perfect. If
$\I(\rho)=[1]$, by \bref{eq-per1-1}, $\{\I\}:\mb=[1]$. If $\I(\rho)$
is perfect, by Lemma \ref{lm-gni1}, $\I=\I^+(\rho)$ is also
perfect.\qedd

Similarly, as a consequence  of Theorem \ref{th-per2}, Lemma
\ref{lm-gni1}, and Theorem \ref{th-nbi2}, we have
\begin{cor}\label{th-nl2p}
Let $\I=\I^+(\rho)$ be a normal binomial $\sigma$-ideal and  $\F$ an
inversive and algebraically closed $\sigma$-field.
Then either $\{\I\}:\mb=[1]$ or $\{\I\}:\mb$ is a binomial
$\sigma$-ideal whose support lattice is the $P$-saturation of
$L_\rho$.
\end{cor}

%
%

In the rest of this section, we give decomposition theorems for
perfect binomial $\sigma$-ideals. We first consider normal binomial
$\sigma$-ideals. By Corollary \ref{cor-rad1} and Example
\ref{ex-31}, a normal binomial $\sigma$-ideal is radical but may not
be perfect.
\begin{theorem}\label{th-nl2}
Let $\I=\I^{+}(\rho)$ be a normal binomial $\sigma$-ideal and $\F$
an inversive and algebraically closed $\sigma$-field. Then
$\{\I\}:\mb$ is either $[1]$ or can be written as the intersection
of reflexive prime binomial $\sigma$-ideals whose support lattice is
the saturation lattice of $L_\rho$.
\end{theorem}
\proof By Theorem \ref{th-l2}, either $\{\I(\rho)\}=[1]$ or
$\{\I(\rho)\}=\bigcap^{s}_{i=1}\I(\rho_i)$, where $\I(\rho_i)$ are
reflexive prime $\sigma$-ideals whose support lattices are
$\sat(L_\rho)$.
By \bref{eq-per1-1} and Lemma \ref{lm-nl1}, either
$\{\I^{+}(\rho)\}:\mb=[1]$ or $\{\I^{+}(\rho)\}:\mb=\{\I(\rho)\}\cap
\F\{\Y\}=(\bigcap^{s}_{i=1}\I(\rho_i))\cap \F\{\Y\}
=\bigcap^{s}_{i=1}(\I(\rho_i)\cap
\F\{\Y\})=\bigcap^{s}_{i=1}\I^{+}(\rho_i)$. By Corollary
\ref{th-nl1}, $ \I^{+}(\rho_i)$ is reflexive and prime whose support
lattices are the saturation of $L_\rho$. \qedd


Now, consider general binomial $\sigma$-ideals.
\begin{lemma}\label{lm-nl5}
$\I\subset\F\{\Y\}$ is a reflexive prime binomial $\sigma$-ideal if
and only if $\I=[y_{i_{1}},\ldots,y_{i_{s}}]$ $+\I_1$, where
$\{y_{i_{1}},\ldots,y_{i_{s}}\}=\Y\cap \I$, $\{z_1,\ldots,z_t\}=\Y
\backslash \I$, and $\I_1$ is a normal binomial reflexive prime
$\sigma$-ideal in $\F\{z_1,\ldots,z_t\}$.
\end{lemma}
\proof If $\I$ is reflexive and prime, then $(y_i^{x^j})^d\in\I$ if
and only if $y_i\in\I$. Let $\I_1 = \I\cap\F\{z_1,\ldots,z_t\}$.
Then $\I=[y_{i_{1}},\ldots,y_{i_{s}}]+\I_1$. $\I_1$ is clearly
reflexive and prime. We still need to show that $\I_1$ is normal.
Let $N f\in\I_1$ for a $\sigma$-monomial $N$ in $\{z_1,\ldots,z_t\}$
and $f\in\F\{z_1,\ldots,z_t\}$. $N$ cannot be in $\I_1$. Otherwise,
some $z_i$ is in $\I_1$ since $\I_1$ is reflexive and prime, which
contradicts to $\{z_1,\ldots,z_t\}=\Y \backslash \I$. Therefore,
$f\in\I_1$ and $\I_1$ is normal.
The other direction is trivial.\qedd

The $\sigma$-ideal $\I$ in Lemma \ref{lm-nl5} is said to be {\em
quasi-normal}. The following result can be proved similarly to
Theorem \ref{th-l2}.
\begin{theorem}\label{th-nl5}
Let $\I$ be a binomial $\sigma$-ideal. If $\F$ is algebraically
closed and inversive, then the perfect $\sigma$-ideal $\{\I\}$ is
either $[1]$ or the intersection of quasi-normal reflexive prime
binomial $\sigma$-ideals.
\end{theorem}
\proof We prove the theorem by induction on $n$.
Let $\I_1=\{\I\}:\mb$. Then $\{\I\}=   \I_1 \cap \cap_{i=1}^n
\{\I,y_i\}$. It is easy to check  $\I_1=\{\I:\mb\}:\mb$. Since
$\I:\mb$ is normal, by Theorem \ref{th-nl2}, $\I_1$ is either $[1]$
or intersection of normal reflexive prime $\sigma$-ideals.
If $n=1$, then $\{\I,y_i\}$ must be either $[y_1]$ or $[1]$. Then
the theorem is proved for $n=1$. Suppose the theorem is valid for
$n=1,\ldots,k-1$.
Still use $\{\I\}=   \I_1 \cap \cap_{i=1}^n \{\I,y_i\}$.
Let $\I_i$ be the $\sigma$-ideal obtained by setting $y_i$ to $0$ in
$\I$. By the induction hypothesis, $\I_i$ can be written as
intersection of quasi-normal reflexive prime $\sigma$-ideals in
$\F\{\Y\setminus\{y_i\}\}$. So the theorem is also valid for
$\{\I,y_i\}=\{\I_i,y_i\}$. The theorem is proved.\qedd

\subsection{Characteristic set for normal binomial $\sigma$-ideal}
The theory of characteristic set given in Section \ref{sec-lbi2} can
be extended to the normal $\sigma$-binomial case.

Let $\rho$ be a partial character over $\Zxn$, $L_{\rho}
=(\f_1,\ldots,\f_s)$ where $\fb=\{\f_1,\ldots,\f_s\}$ is a reduced
Gr\"obner basis, and
 \begin{equation}\label{eq-A+}
 \A^{+}(\rho): \Y^{\f_1^+}-\rho(\f_1)\Y^{\f_1^-},
      \ldots,
      \Y^{\f_s^+}-\rho(\f_s)\Y^{\f_s^-}.\end{equation}

We have the following canonical representation for normal binomial
$\sigma$-ideals.
\begin{theorem}\label{th-nl3}
Use the notations in \bref{eq-A+}.
Then $\I^{+}(\rho)=\sat(\A^{+}(\rho))$. Furthermore,
$\A^{+}(\rho)$ is a regular and coherent $\sigma$-chain and hence is
a characteristic set of $\I^{+}(\rho)$.
\end{theorem}
\proof
Let $\I_1=[\A^{+}(\rho)] :\mb$. We claim $\I_1=\sat(\A^{+}(\rho))$.
It is clear that $\sat(\A^{+}(\rho))\subset
[\A^{+}(\rho)]:\mb=\I_1$.
For the other direction, let $p\in\I_1$ and
$p_1=\prem(p,\A^{+}(\rho))$ which is reduced w.r.t. $\A^{+}(\rho)$.
By \bref{eq-prem}, $p_1\in\I_1$. As a consequence,
$p_1\in[\A(\rho)]$ as Laurent $\sigma$-polynomials in
$\F\{\Y^{\pm}\}$. By Lemma \ref{lm-l1}, $\A(\rho)$ is a
characteristic set of $[\A(\rho)]$. Since $p_1$ is reduced w.r.t.
$\A^+(\rho)$, it is also reduced w.r.t. $\A(\rho)$. Then $p_1=0$ and
hence the claim is proved.

We now prove $\I^+(\rho) = \sat(\A^+(\rho))$.
By the above claim, Lemma \ref{lm-l1}, and Lemma \ref{lm-nl1},
$\sat(\A^{+}(\rho))=[\A^{+}(\rho)]
:\mb=[\A^{+}(\rho)]\F\{\Y^{\pm}\}\cap \F\{\Y\}=[\A(\rho)]\cap
\F\{\Y\}=\I(\rho)\cap \F\{\Y\}=\I^{+}(\rho) $.

It remains to prove that $\A^{+}(\rho)$ is a characteristic set of
$\I_1=[\A^{+}(\rho)] :\mb$. By definition, it suffices to show that
if $p\in\I_1$ is reduced w.r.t. $\A^{+}(\rho)$ then $p=0$. Let
$A_i=\Y^{\f_i}-\rho(\f_i)$ and
$A_i^+=\Y^{\f_i^+}-\rho(\f_i)\Y^{\f_i^-}$. Since $p\in\I_1$, there
exist a $\sigma$-monomial $M$ and $f_{i,j}\in\F\{\Y\}$ such that
$Mp=\sum_{i,j}f_{i,j} (A_{i}^+)^{x^j}$. Then in $\F\{\Y^{\pm}\}$, we
have $p=\sum_{i,j}g_{i,j} A_{i}^{x^j}\in[\A(\rho)]$, where
$g_{i,j}\in\F\{\Y^{\pm}\}$. Since $p$ is reduced w.r.t.
$\A^+(\rho)$, it is also reduced w.r.t. $\A(\rho)$. By Lemma
\ref{lm-l1}, $\A(\rho)$ is a characteristic set of $[\A(\rho)]$ and
hence $p=0$. The claim is proved.

Since $\I_1=\sat(\A^+(\rho))$, $\A^{+}(\rho)$ is also a
characteristic set of $\sat(\A^+(\rho))$. By Theorem \ref{th-rp},
 $\A^{+}(\rho)$ is regular and coherent.\qedd

\begin{example}\label{ex-ig}
Let $L=([1-x,x-1]^\tau)$ be a $\Zx$-module and $\rho$ the trivial
partial character on $L$, that is, $\rho(\f)=1$ for $\f\in$L.
By Theorem \ref{th-nl3}, $\I^{+}(\rho)=\sat[y_1y_2^{x}-y_1^{x}y_2]
\subseteq \Q\{y_1, y_2\}$. By Theorem \ref{lm-nl1}, $\I^{+}(\rho)$
is a reflexive prime $\sigma$-ideal.
We can show that
$\I^{+}(\rho)=[y_1^{x^{i}}y_2^{x^{j}}-y_1^{x^{j}}y_2^{x^{i}}\,|\, 0
\leq i < j \leq m]$, which is an infinitely generated
$\sigma$-ideal.
\end{example}

As a consequence of Theorem \ref{th-nbi2}, Theorem \ref{th-nl3}, and
Lemma \ref{lm-l1}, we have
\begin{cor}\label{cor-ascr}
Let $\A(\rho)$ and  $\A^{+}(\rho)$ be defined in \bref{eq-Arho} and
\bref{eq-A+}, respectively. Then\newline
$([\A(\rho)]\F\{\Y^{\pm}\})\cap\F\{\Y\}=\sat(\A^{+}(\rho))$.
\end{cor}


In order to prove the converse of Theorem \ref{th-nl3}, we need the
following lemmas.
Let $\f_i\in\Zxn$ and $c_i\in\F^*$, $i=1,\ldots,s$. Consider the
following $\sigma$-chains
\begin{eqnarray}
 \A&:&\Y^{\f_{1}}-c_{1}, \ldots, \Y^{\f_{s}}-c_{s}\label{eq-aa}\\
 \A^{+}&:&\Y^{\f_{1}^+}-c_{1}\Y^{\f_{1}^-}, \ldots,
 \Y^{\Y^{\f_{s}^+}}-c_{s}\Y^{\f_{s}^-}\nonumber
 \end{eqnarray}
in $\F\{\Y^{\pm}\}$ and $\F\{\Y\}$, respectively. Since $\f_i$ are
assumed to be normal, $\A^+$ is a $\sigma$-chain if and only if $\A$
is a Laurent $\sigma$-chain.

\begin{lemma}\label{lm-asc1}
Use the notations in \bref{eq-aa}. Let
$p=a\Y^{\a}+b\Y^{\b}=aN(\Y^{\f}-c)\in\F\{\Y\}$, where
$\a,\b\in\N[x]^n$, $\f\in\Zxn$, $N$ is a $\sigma$-monomial,
$c\in\F^*$.
If $\A^+$ is coherent and regular, then $\prem(p,\A^+)=0$ implies
$\prem(\Y^{\f}-c,\A)=0$.
\end{lemma}
\proof
%
%
%
Since $\prem(p,\A^+)=0$, there exists a $\sigma$-monomial $M_1$ such
that $M_1p\in[\A^+]$.
Let $p_1 = \Y^{\f}-c$. Since $r_1=\prem(p_1,\A)=\Y^\g-c_\g$, by
Lemma \ref{lm-rem1}, there exists a $c_1\in\F^*$ such that
$r_1-c_1p_1\in[\A]$. Then, there exists a $\sigma$-monomial $M_2$
such that $M_2Nr_1, M_2Np_1\in\F\{\Y\}$ and
$M_2N(r_1-c_1p_1)\in[\A^+]$ and hence
$M_2M_1N(r_1-c_1p_1)=M_2M_1Nr_1-\frac{c_1}{a}M_2M_1p\in[\A^+]$. Let
$M=M_1M_2N$. From $M_1p\in[\A^+]$, we have
$Mr_1\in[\A^+]\subset\sat(A^+)$.

Suppose $A_i= \Y^{\f_{i}^+}-c\Y^{\f_{i}^-}=
 I_i^+ y_{c_i}^{d_ix^o_i}-cI_i^-$, where $y_{c_i}$ is the leading variable of $A_i$.
A variable like $y_{c_i}^{x^{o_i+k}}$ for $k\in\N$ is called a {\em
main variable} of $\A^+$. A variable $y_{i}^{x^j}$ is called a {\em
parameter} of $\A^+$ is it is not a main variable.
If $M$ contains a main variable of $\A^+$ as a factor. Then let
$z=y_{c_i}^{x^{o_i+k}}$ be the largest one appearing in $M$ under
the variable ordering induced by the lexicographical of the index
$(c_i,o_i+k)$. Let $s=\deg(M,z)$ and $M_1=M/(z^s)$.
We may assume that $d_i$ is a factor of $s$. Otherwise, let
$s_1=\lfloor \frac{s}{d_i}\rfloor$, $s_0 = s-s_1d_i$, and
$M=Mz^{d_i-s_0}=M_1z^{d_i(s_1+1)}$. We still have
$Mr_1\in\sat(A^+)$.
We may use $A_i=0$ to eliminate $z$ from $M$ as follows:
 $M_1z^{s-d_i}(cI_i^-)^{x^k} r_1 =
 M_1z^{s-d_i}(I_i^+ y_{c_i}^{d_ix^o_i}-A_i)^{x^k} r_1=
 M(I_i^+)^{x^k}r_1-M_1z^{s-d_i}(A_i)^{x^k} r_1 \in\sat(A^+)$.
Note that $\deg(M_1z^{s-d_i}(cI_i^-)^{x^k},z)=s-d_i$. Repeat the
above procedure, we may find a $\sigma$-monomial $N$ such that $Nr_1
\in\sat(A^+)$, $N$ does not contain $z$ as a factor, and any
variable $y_i^{x^j}$ in $M$ is smaller than $z$ in the given
variable ordering.
Repeat the procedure, we may finally obtain a $\sigma$-monomial $L$
such that $L$ does not contain main variables of $\A^+$ as factors
and $Lr_1\in\sat(A^+)$. Since $L$ contains only parameters of $\A^+$
and $r_1$ is reduced w.r.t. $\A^+$, $Lr_1$ is also reduced w.r.t.
$\A^+$. Since $\A^+$ is regular and coherent, by Lemma \ref{lm-21},
it is the characteristic set of $\sat(A^+)$. Therefore, $Lr_1=0$,
and $r_1=0$. \qedd

The following example shows that if $\prem(p,\A^+)\ne0$ then the
relation between $\prem(p,$ $\A^+)$ and $\prem(\Y^{\f}-c,\A)$ may be
complicated, where $p=a\Y^{\a}+b\Y^{\b}=aN(\Y^{\f}-c)$.
\begin{example}\label{ex-asc2}
Let $p=y_2(y_2-1)$, $A_1=y_1^{-1}y_2^2-1$, and $A_1^+=y_2^2-y_1$.
Then $\prem(p,A_1^+)=y_1-y_2$ in $\F\{y_2\}$. But in
$\F\{y_2^{\pm}\}$, $p$ is represented as $\widetilde{p}=y_2-1$ and
$\prem(\widetilde{p},A_1)=y_2-1$.
\end{example}

\begin{lemma}\label{lm-asc2}
Use the notations in \bref{eq-aa}. $\A$ is a regular and coherent
$\sigma$-chain in $\F\{\Y^{\pm}\}$ if and only if $\A^{+}$ is a
regular and coherent $\sigma$-chain in $\F\{\Y\}$.
\end{lemma}
\proof
If $\A$ is regular and coherent, by Corollary \ref{cor-l11}, there
exists a partial character $\rho$ over $\Zxn$ such that
$L_\rho=(\f_1,\ldots,\f_s)$, $\rho(\f_i)=c_i$, and $\I(\rho)=[\A]$.
By Theorem \ref{th-nl3}, $\A^+=\A^+(\rho)$ is regular and coherent.

Assume that $\A^+$ is regular and coherent. We first show that
$[\A]\ne[1]$ in $\F\{\Y^{\pm}\}$. It suffices to show that
$\sat(\A^+)$ does not contain a $\sigma$-monomial. Suppose the
contrary, there is a $\sigma$-monomial $M\in\sat(\A^+)$. Since
$\A^+$ is a regular and coherent chain, we have $\prem(M,\A^+)=0$.
Now consider the procedure of $\prem$, it can be shown that the
pseudo-remainder of a nonzero $\sigma$-monomial w.r.t. a binomial
$\sigma$-chain is still a nonzero  $\sigma$-monomial, a
contradiction.

%
Note that $\A$ is always regular since $\sigma$-monomials are
invertible in $\F\{\Y^{\pm}\}$.  Then, it suffices to prove that
$\A$ is coherent.

Let $A_i=\Y^{\f_{i}}-c_{i}$ and
$A_i^+=\Y^{\f_{i}^+}-c_{i}\Y^{\f_{i}^-}$.
Assume $A_i^+$ and $A_j^+$ ($i<j$) have the same leading variable
$y_l$, and
    $A_i^+=I_i^+ y_l^{d_i x^{o_i}}-c_iI_i^-$,
    $A_j^+=I_j^+ y_l^{d_j x^{o_j}}-c_jI_j^-$, where $I_i^-=\Y^{\f_{i}^-}$.
From Definition \ref{def-ghf},  we have $o_i < o_j$ and $d_i| d_j$.
Let $d_i = t d_j$ where $t\in\N$.
From \bref{eq-delta},
$$\Delta(A_i^+,A_j^+)
=\prem((A_i^+)^{x^{o_j-o_i}},A_j^+) =
       c_j^t (I_j^-)^t (I_i)^{x^{o_j-o_i}} -
       (I_j^+)^t  (c_iI_i^+)^{x^{o_j-o_i}}.$$
Comparing to \bref{eq-delta1}, if $\Delta(A_i,A_j)=\Y^\h-c_f$, then
$\Delta(A_i^+,A_j^+)=c_j^tM(\Y^{\h^+}-c_f\Y^{\h^-})$, where $M$ is a
$\sigma$-monomial.
Since $\A^+$ is coherent, $\prem(\Delta(A_i^+,A_j^+),\A^{+})$ $=0$.
By Lemma \ref{lm-asc1}, $\prem(\Delta(A_i,A_j),\A)=0$ which implies
that $\A$ is coherent.\qedd

We now prove the converse of Theorem \ref{th-nl3}.
\begin{theorem}\label{th-nl4}
Use the notations in \bref{eq-aa}.
%
%
If $\A^+$ is a regular and coherent $\sigma$-chain, then there is a
partial character $\rho$ over $\Zxn$ such that
$L_\rho=(\f_1,\ldots,\f_s)$, $\rho(\f_i)=c_i$, $\I(\rho)=[\A]$, and
$\I^{+}(\rho)=\sat(\A^{+})$
\end{theorem}
\proof
By Lemma \ref{lm-asc2}, $\A$ is regular and coherent.
By Theorem \ref{th-t4}, $\f$ is a reduced Gr\"obner basis for a
$\Zx$-lattice and $[\A]\subset\F\{\Y^{\pm}\}$ is proper.
By Corollary \ref{cor-l11}, there exists a partial character $\rho$
such that $L_\rho=(\f_1,\ldots,\f_s)$, $\rho(\f_i)=c_i$, and
$\I(\rho)=[\A]$. By Theorem \ref{th-nl3},
$\I^{+}(\rho)=\sat(\A^{+}(\rho))=\sat(\A^{+})$.\qedd

As a consequence of Theorem \ref{th-nl4} and Lemma \ref{lm-nl1}, we
have
\begin{cor}\label{cor-bi41}
Normal binomial $\sigma$-ideals are in a one to one correspondence
with $\sat(\A^{+})$, where $\A^+$ is a regular and coherent chain
given in \bref{eq-aa}.
\end{cor}
As a consequence of Corollary \ref{cor-ascr} and Theorem
\ref{th-nl4}, we have
\begin{cor}\label{cor-ascr1}
$\F\{\Y^{\pm}\}[\A]\cap\F\{\Y\}=\sat(\A^{+})$.
\end{cor}

%

As a consequence of Theorem \ref{th-pr1}, Corollary \ref{th-nl1},
and Theorem \ref{th-nl4}.
\begin{cor}\label{cor-asc2}
%
%
$[\A]$ is a reflexive (prime) $\sigma$-ideal in $\F\{\Y^{\pm}\}$ if
and only if $\sat(\A^+)$ is a reflexive (prime) $\sigma$-ideal in
$\F\{\Y\}$.
\end{cor}

\subsection{Perfect closure of binomial $\sigma$-ideal and binomial $\sigma$-variety}
In this section, we will show that the perfect closure of a binomial
$\sigma$-ideal is also binomial. We will also give a geometric
description of the zero set of a binomial $\sigma$-ideal.
%
%
For the perfect closure of a binomial $\sigma$-ideal, we have
\begin{theorem}\label{diff-perfect1}
Let $\F$ be an algebraically closed and inversive $\sigma$-field.
Then the perfect closure of a binomial $\sigma$-ideal $\I$ is
binomial.
\end{theorem}

We first remark that it is not known wether the well-mixed closure of a binomial $\sigma$-ideal is still binomial.
Before proving Theorem \ref{diff-perfect1}, 
we first prove several lemmas. In the rest of this section, we
assume that $\I\subseteq S=\F\{\Y\}$ and $\m$ the set of
$\sigma$-monomials in $S$.
\begin{lemma}\label{diff-anybi}
If $\I$ is a binomial $\sigma$-ideal, then
$\{\I\}:\m$ is either $[1]$ or a binomial  $\sigma$-ideal.
\end{lemma}
\proof It is easy to check
$\{\I\}\F\{\Y^{\pm}\}=\{\I\F\{\Y^{\pm}\}\}$.
By \bref{eq-nbi1},
$\{\I\}:\m=\{\I\}\F\{\Y^{\pm}\}\cap\F\{\Y\}=\{\I\F\{\Y^{\pm}\}\}\cap\F\{\Y\}$.
Now the lemma follows from  Theorem \ref{th-per2}.\qedd

\begin{lemma}\label{diff-radical}
If $\I$ is a $\sigma$-ideal in $\F\{\Y\}$, then
\begin{equation}\label{eq-radical}
  \{\I\}=\{\I\}:\mb\cap \{\I+y_1\}\cap \cdots \cap \{\I+y_n\}
\end{equation}
\end{lemma}
\proof The right hand side of \bref{eq-radical} clearly contains
$\{\I\}$. It suffices to show that every reflexive prime $P$
containing  $\I$ contains one of the $\sigma$-ideals on the
right-hand side of \bref{eq-radical}. If $\{\I\}:\mb \subseteq P$,
we are done. Otherwise, there exists an element $f\in
(\{\I\}:\mb)\setminus P$ which implies that there exists a
$\sigma$-monomial $M$ such that $Mf \in \{\I\}\subseteq P$. This
implies $y_i\in P$ for some $i$. Thus, $P$ contains $\{\I+y_i\}$ as
required.\qedd

\begin{lemma}\label{diff-bm1}
Let $\I$ be a binomial $\sigma$-ideal in $S=\F\{\Y\}$ and
$S'=\F\{y_1,\ldots,y_{n-1}\}$. If $\I'=\I\cap S'$, then $[\I+y_n]$
is the sum of $[\I'S+y_n]$ and a monomial $\sigma$-ideal in $S'$.
\end{lemma}
\proof  Every $\sigma$-binomial involving $y^{x^{k}}_n$ is either
contained in $[y_n]$ or is congruent modulo $[y_n]$ to a
$\sigma$-monomial in $S'$. Thus, all generators of $\I$ which are
not in $\I'$ may be replaced by $\sigma$-monomials in $S'$ when
forming a generating set for $[\I+y_n]$.\qedd

\begin{lemma}\label{diff-perfect}
Let $\I$ be a perfect binomial $\sigma$-ideal in $S=\F\{\Y\}$. If
$\MI$ is a $\sigma$-monomial $\sigma$-ideal, then
$\{\I+\MI\}=[\I+\MI_1]$ for some monomial $\sigma$-ideal $\MI_1$.
\end{lemma}
\proof
If $1\in\MI$, then the lemma is obviously valid. Otherwise,
$[\I+\MI]:\mb=[1]$. Lemma \ref{diff-radical} yields $\{\I+\MI
\}=\bigcap_{i=1}^{i=n}\{\I+\MI+y_i\}$.
By Corollary \ref{cor-tttt1}, we need only to show that
$\{\I+\MI+y_i\}$ is the sum of $\I$ and a monomial $\sigma$-ideal.
For simplicity, let $i=n$ and write
$S'=\F\{y_1,y_2,\ldots,y_{n-1}\}$. Since $\I$ is perfect, the
$\sigma$-ideal $\I'=\I\cap S'$ is perfect as well. By Lemma
\ref{diff-bm1}, $[\I+\MI+y_n]=[\I'S+\MI'S+y_n]$ where $\MI'$ is a
monomial $\sigma$-ideal in $S'$. By induction on $n$, the perfect
closure of $\I'+\MI'$ in $S'$ has the form $\I'+\MI_1'$, where
$\MI_1'$ is a monomial $\sigma$-ideal of $S'$. Putting this
together, we have
\[
 \begin{array}{llll}
 \{\I+\MI+y_n\} & =& \{\I'S+\MI'S+y_n\} \\
 &=& [\I'S+\MI_1'S+y_n]\\
 &\subseteq & [\I+\MI_1'S+y_n] \\
 &\subseteq&  \{\I+\MI+y_n\}
\end{array}\]
So $\{\I+\MI+y_n\}=[\I+\MI_1'S+y_n]$ is $\I$ plus a monomial
$\sigma$-ideal, as required.\qedd

%
\vskip5pt{\noindent\em Proof of Theorem \ref{diff-perfect1}:} We
will prove the theorem by induction on $n$.
By Lemma \ref{diff-anybi}, $\I_1=\{\I\}:\mb$ is binomial.
For $n=1$, by Lemma \ref{diff-radical}, $ \{\I\}=\I_1\cap
\{\I+y_1\}$. If $\{\I+y_1\}=1$ then $ \{\I\}=\I_1$ is binomial.
Otherwise $\{\I+y\}=[y]$ and hence $\I\subset [y]$. Since $\I\subset
\I_1$, $\{\I\}=\I_1\cap [y]=[\I+\I_1]\cap [\I+y]$ is binomial by
Lemma \ref{cor-gbi8}.
Suppose the lemma is valid for $n-1$ variables and let $\I$ be a
binomial $\sigma$-ideal in $S=\F\{\Y\}$. Let $\I_j:=\I \cap S_j$,
where $S_j=\F\{y_1,\ldots,y_{j-1},y_{j+1},\ldots,y_n\}$. By the
induction hypothesis, we may assume that the perfect closure of each
$\I_j$ is binomial. Adding these binomial $\sigma$-ideals to $\I$,
we may assume that each $\I_j$ is perfect begin with.
By Lemma \ref{diff-anybi}, $\I_1=\{\I\}:\mb$ is binomial. Then there
exists a binomial $\sigma$-ideal $\I'$, say $\I'=\I_1$,  such that
$\I_1=[\I+\I']$. By Lemma \ref{diff-bm1}, $[\I+y_j]=[\I_jS+\J_j
S+y_j]$, where $\J_j$ is a monomial $\sigma$-ideal in $S_j$. Since
$\I_j$ is perfect, the $\sigma$-ideal $\I_jS$ is perfect, so we can
apply Lemma \ref{diff-perfect} with $\MI=[\J_j S+y_i]$ to see that
there exists a monomial $\sigma$-ideal $\MI_j$ in $S$ such that
$\{\I+y_j\}=\{\I_jS+\J_j S+y_j\}=[\I_jS+\MI_j]=[\I+\MI_j]$. By Lemma
\ref{diff-radical} and Corollary \ref{cor-gbi8},
$\{\I\} = [\I+\I']\bigcap\cap_{j=1}^n[\I+\MI_j]$ is binomial.\qedd

\begin{example}\label{ex-perc11}
Let $p=y_2^{2}-y_1^{2}$. Following the proof of Theorem
\ref{diff-perfect1},  $\{p\} = (\{p\}:\m)\cap[y_1,y_2]$.
By Example \ref{ex-perc1} and Corollary \ref{cor-ascr1},
 $\I_1=\{p\}:\m=\sat[y_2^{2}-y_1^{2},y_1y_2^{x}-y_1^{x}y_2]= [y_1y_2^{x^i}-y_1^{x^i}y_2,y_2^{1+x^j}-y_1^{1+x^j}\,|\,i,j\in\N]$.
Thus, $\{p\} = \I_1\cap[y_1,y_2]=\I_1$.
%
%
%
\end{example}

A $\sigma$-ideal $\I$ is called {\em normal} if for any $M\in\mb$
and $p\in\F\{\Y\}$,  $Mp\in\I$ implies $p\in\I$. The following
example shows that the perfect closure of a normal $\sigma$-ideal
could be not normal.
\begin{example}
Let $\I=[y_2^2y_4^2-y_1^2y_3^2,y_1^x-y_1,y_2^x-y_2]$.
By Example \ref{ex-perc11}, we have
 $\{\I\}=[p_i,
  (y_2y_4)^{1+x^j}-(y_1y_3)^{1+x^j},
 y_1^x-y_1,y_2^x-y_2\,|\,i,j\in\N]$, where
$p_i= y_1y_3(y_2y_4)^{x^i}-y_2y_4(y_1y_3)^{x^i}$.
Note that $p_i=y_1y_2(y_3y_4^{x^i}-y_4y_3^{x^i})$  modulo
$[y_1^x-y_1,y_2^x-y_2]$ and $y_3y_4^{x^i}-y_4y_3^{x^i}$ is not in
$\{\I\}$.
%
%
%
\end{example}


In the rest of this section, we give a geometric description of the
zero set of a binomial $\sigma$-ideal, which is a generalization of
Theorem 4.1 in \cite{es-bi} to the difference case. The basic idea
of the proof also follows \cite{es-bi}, except we need to consider
the distinction between the perfect $\sigma$-ideals and radical
ideals.


We decompose the affine $n$-space $(\AH)^{n}$ into the union of
$2^{n}$ $\sigma$-coordinate flats:
$$ (\AH^{*})^{\Omega}:=\{(a_1,a_2,\ldots,a_n)\,|\,a_i\neq 0 , i\in \Omega ; a_i=0, i \notin \Omega  \}$$
where $\Omega$ runs over all subsets of $\{1,2,\ldots,n\}$. The Cohn
closure of $ (\AH^{*})^{\Omega}$ in $(\AH)^{n}$ is defined by the
$\sigma$-ideal
 $$M(\Omega):=[y_i|i\notin \Omega]\subset\F\{\Y\}.$$
 The
$\sigma$-coordinate ring of $ (\AH^{*})^{\Omega}$  is the Laurent
polynomial $\sigma$-ring
$\F\{\Omega^{\pm}\}:=\F\{y_i,y_i^{-1},i\in\Omega\}$. We can define a
coordinate projection $ (\AH^{*})^{\Omega'}\longrightarrow
(\AH^{*})^{\Omega}$ whenever $\Omega\subseteq \Omega'\subseteq
\{1,2,\ldots,n\}$ by setting all those coordinates not in $\Omega$
to zero.

 If $X$ is any $\sigma$-variety of $(\AH)^{n}$ and
 $\I=\IH(X)\subseteq \F\{\Y\}$, then the Cohn closure of the
 intersection of $X$ with $(\AH^{*})^{\Omega}$ corresponds to the
 $\sigma$-ideal
 $$\I_{\Omega}:=[\I+M(\Omega)]:\m_{\Omega}\subset \F\{\Y\}$$
  where $
 \m_{\Omega}=\{\prod_{i\in \Omega}y_i^{m_i(x)}|m_i(x)\in \N[x]\}$.
%
%
Since $\I$ is perfect, by the difference Nullstellsatz
\cite[p.87]{cohn}
 $$\I=\bigcap_{\Omega}\{\I_{\Omega}\}.$$
 If $\I$ is binomial, then by
 Corollary \ref{diff-quotient} the $\sigma$-ideal $\I_{\Omega}$ is
 also binomial.

\begin{lemma}\label{diff-4}
Let $R:=\F\{z_1,z_1^{-1},\ldots,z_t,z_t^{-1}\}\subset
R':=\F\{z_1,z_1^{-1},\ldots,z_t,z_t^{-1},y_1,\ldots,y_s\}$ be a
Laurent polynomial $\sigma$-ring and a polynomial $\sigma$-ring over
it. If $B\subset R'$ is a binomial $\sigma$-ideal and $M\subset R'$
is a monomial $\sigma$-ideal such that $[B+M]$ is a proper
$\sigma$-ideal in $R'$, then $$[B+M]\cap R=B\cap R.$$
\end{lemma}
\proof  This is a $\sigma$-version of \cite[Lemma 4.2]{es-bi}, which
can be proved similarly.\qedd

We can make a classification of all binomial $\sigma$-varieties $X$
by intersecting $X$ with $(\AH^{*})^{\Omega}$, since by Theorem
\ref{diff-perfect1}, the perfect closure of a binomial
$\sigma$-ideal is still binomial.
\begin{theorem}\label{diff-structure}
Let $\F$ be any algebraically closed and inversive $\sigma$-field. A
$\sigma$-variety $X\subset \AH^{n}$ is generated by
$\sigma$-binomials if and only if  the following three conditions
hold.
\begin{description}
\item[(1)] For each $(\AH^{*})^{\Omega}$, the $\sigma$-variety
$X\cap(\AH^{*})^{\Omega}$ is generated by $\sigma$-binomials.

\item[(2)] The family of sets $U=\{\Omega\subseteq
\{1,2,\ldots,n\}|X\cap(\AH^{*})^{\Omega}\neq \emptyset$\} is closed
under taking intersections.

\item[(3)] If $\Omega_1,\Omega_2\in U $ and $\Omega_1 \subset \Omega_2$,
then the coordinate projection $(\AH^{*})^{\Omega_2}\longrightarrow
(\AH^{*})^{\Omega_1} $ maps $X\cap (\AH^{*})^{\Omega_2}$ onto a
subset of $X\cap (\AH^{*})^{\Omega_1}$.
\end{description}

\end{theorem}

The above theorem can be reduced to the following algebraic version.
\begin{theorem}\label{diff-structure-algebra}
Let $\F$ be any algebraically closed and inversive $\sigma$-field. A
perfect $\sigma$-ideal $\I\subset \F\{\Y\}$ is binomial  if and only
if the following three conditions hold.
\begin{description}
\item[(1)]
 For each $\Omega\subseteq\{1,\ldots,n\}$, $\I_{\Omega}$ is
binomial.

\item[(2)] $U=\{\Omega\subseteq
\{1,2,\ldots,n\}\,|\,\{\I_{\Omega}\}\neq [1]$\} is closed under
taking intersections.

\item[(3)] If $\Omega_1,\Omega_2\in U $ and $\Omega_1 \subset \Omega_2$,
then $\I_{\Omega_1}\cap\F\{\Omega_1\}\subset\I_{\Omega_2}$, where
$\F\{\Omega_1\}=\F\{y_i\,|\, y_i\in\Omega_1\}$.
\end{description}
\end{theorem}

\proof
%
Suppose $\I$ is a perfect $\sigma$-ideal in $\F\{\Y\}$.
Since  $\I$ is binomial, by Lemma \ref{diff-quotient} $\I_{\Omega}$
is also binomial and (1) is proved.
To prove (2) by contradiction, assume that for $\Omega_1,\Omega_2\in
U$, $\{\I_{\Omega_1}\}\neq[1], \{\I_{\Omega_2}\}\neq[1],
\{\I_{\Omega_1\cap\Omega_2} \}=[1]$.  We consider two cases.
If $\I_{\Omega_1\cap\Omega_2} =[1]$, then for some $m(x)\in \N[x]$
we have  $(\prod_{i\in\Omega_1\cap\Omega_2
}y_i)^{m(x)}\in[\I+M(\Omega_1)+M(\Omega_2)]$. By Corollary
\ref{cor-tttt1},  $(\prod_{i\in\Omega_1\cap\Omega_2 }y_i)^{m(x)}$ is
either in $[\I+M(\Omega_1)]$ or $[\I+M(\Omega_2)]$, so
$\I_{\Omega_1}$ or $\I_{\Omega_2}$ is $[1]$.
For the second case, we have $\I_{\Omega_1\cap\Omega_2} \neq[1]$ and
$\{\I_{\Omega_1\cap\Omega_2}\} =[1]$. Then there exist a finite
number of proper $\sigma$-binomials $B_1,\ldots,B_s$ and
$\sigma$-monomials $m_1,\ldots,m_s$ in $\F\{\Omega_1\cap\Omega_2\}$
such that $m_iB_i\in \I$ and $\{B_1,\ldots,B_s,
y_i,i\notin\Omega_1\cap\Omega_2\}=[1]$. We thus have $
\{B_1,\ldots,B_s\}=[1]$. Since $m_iB_i\in
\I\cap\F\{\Omega_1\cap\Omega_2\}$, we have $B_i\in\I_{\Omega_1}$ and
$B_i\in\I_{\Omega_2}$ and thus
$\{\I_{\Omega_1}\}=\{\I_{\Omega_2}\}=[1]$.
To prove (3), given $\Omega_1,\Omega_2\in U$ and
$\Omega_1\subset\Omega_2$, we have
$\I_{\Omega_2}=[\I_{\Omega_2}:\m_{\Omega_1} ]$. Set
$R'=\F\{\Omega_1^{\pm}\}\{\{y_i\}_{i\notin \Omega_1}\}$, then
 $$[\I+M(\Omega_1)]R'\cap \F\{\Omega_1^{\pm}\}\subseteq\I_{\Omega_2}R'.$$
Since $\Omega_1\in U$, the $\sigma$-ideal $[\I+M(\Omega_1)]R^{'}$ is
proper. By Lemma \ref{diff-4}, we have $[\I+M(\Omega_1)]R'\cap
\F\{\Omega_1^{\pm}\}=\I R'\cap\F\{\Omega_1^{\pm}\}\subset
\I_{\Omega_2}R'\cap\F\{\Omega_1^{\pm}\}$. So
$\I_{\Omega_1}\cap\F\{\Omega_1\}\subset\I_{\Omega_2}$.

To prove the other driection, let $\I$ be a perfect $\sigma$-ideal
satisfying the three conditions. By the difference Nullstellensatz,
$\I=\cap_{\Omega\in U}\{\I_{\Omega}\}$. By condition (2), $U$ is a
partially ordered set under the inclusion for subsets of
$\{1,\ldots,n\}$. For each $\Omega\in U$, we set
$\J(\Omega)=[\I_{\Omega}\cap \F\{\Omega\}]\F\{\Y\}$ with the
properties that if $\Omega_1\subset \Omega_2$,
$\{\J(\Omega_1)\}\subset\{\J(\Omega_2)\}$. Note that
$[M_{\Omega_1\cap\Omega_2}]\subset[M_{\Omega_1}+M_{\Omega_2}]$. Then
we have
 $$ \I=\cap_{\Omega\in U}\{\I_{\Omega}\}=\cap_{\Omega\in
 U}\{\J(\Omega)+M(\Omega)\}.$$
Now we will prove that
 \begin{equation}\label{eq-pc3}
\cap_{\Omega\in U}\{\J(\Omega)+M(\Omega)\}=\{\cap_{\Omega\in
U}M(\Omega)+\sum_{\Omega\in U}\{\J(\Omega)\cap (\cap_{\Omega_\eta
\nsupseteq \Omega}M(\Omega_\eta))\}\}.  \end{equation}
If $\Omega_2\supseteq\Omega_1$, we have
$\{\J(\Omega_2)+M(\Omega_2)\}\supseteq\{\J(\Omega_2)\}\supseteq\{\J(\Omega_1)\}\supseteq
\{\J(\Omega_1)\cap \cap_{\Omega_\eta \nsupseteq
\Omega_1}M(\Omega_\eta)\}$.
If  $\Omega_2\nsupseteq\Omega_1$, we have
$\{\J(\Omega_2)+M(\Omega_2)\}\supseteq M\{\Omega_2\}\supseteq
\{\J(\Omega_1)\cap \cap_{\Omega_\eta \nsupseteq
\Omega_1}M(\Omega_\eta)\}$. So the left hand side contains the right
hand side of \bref{eq-pc3}.
For the other direction, consider a reflexive prime $\sigma$-ideal
$P\supseteq[\cap_{\Omega\in U}M(\Omega)+\sum_{\Omega\in
U}\{\J(\Omega)\cap \cap_{\Omega_\eta \nsupseteq
\Omega}M(\Omega_\eta)\}]$ and set $V=\{\Omega\in U|M(\Omega)\subset
P\}$. Then $V$ is a finite partially order set and nonempty since
$P\supseteq\cap_{\Omega\in U}M(\Omega)$ and
$\{M_{\Omega_1\cap\Omega_2}\}\subset\{M_{\Omega_1}+M_{\Omega_2}\}$.
Let $\Omega_{0}$ be the smallest element of $V$ such that
$P\supseteq M_{\Omega_{0}}$. At the same time,
$P\supset\J(\Omega_0)\cap \cap_{\Omega_\eta \nsupseteq
\Omega_0}M(\Omega_\eta)$, then $P\supseteq \J(\Omega_0)$. Therefore,
$P\supseteq \J(\Omega_0)+M(\Omega_0)$ and $P$ contains the left hand
side \bref{eq-pc3} and \bref{eq-pc3} is proved. Since
$\cap_{\Omega\in U}M(\Omega)+\sum_{\Omega\in U}\{\J(\Omega)\cap
\cap_{\Omega_\eta \nsupseteq \Omega}M(\Omega_\eta)\}$ is binomial,
the theorem follows from \bref{eq-pc3}.\qedd
%

\section{Toric $\sigma$-ideal and toric $\sigma$-variety}
\label{sec-tv} In this section, we will introduce the concept of
toric $\sigma$-variety and prove some of its basic properties.

\subsection{Toric $\sigma$-variety}
\label{sec-tv1} Let $\F$ be an inversive $\sigma$-field. Following
\cite{wibmer},  we denote the category of $\sigma$-field extensions
of $\F$ by $\mathscr{E}_\F$,  the category of $\E^{n}$ by
$\mathscr{E}_{\F}^{n}$ where $\E \in \mathscr{E}_{\F}$. Let
$(\AH^{*})^{n}$ be the functor from $\mathscr{E}_{\F}$ to
$\mathscr{E}_{\F}^{n}$ satisfying $(\AH^{*})^{n}(\E)=(\E^{*})^{n}$
where $\E \in \mathscr{E}_{\F}$ and
 $\E^{*}=\E\setminus\{0\}$. Let $\AH^{n}$ be the
functor from $\mathscr{E}_{\F}$ to $\mathscr{E}_{\F}^{n}$ satisfying
$(\AH)^{n}(\E)=(\E)^{n}$ where $\E \in \mathscr{E}_{\F}$.
%
%
A {\em $\sigma$-variety} over $\F$ is a functor $\V$ from
$\mathscr{E}_\F$ to the category of sets with the form $\V(P)$ for
$P\subset\F\{\Y\}$ satisfying $\V_{\E}(P)= \{\eta\in\E^n\,|\,
\forall p\in P, p(\eta)=0\}$.

In the rest of this section, let
\begin{equation}\label{eq-A}
 \falpha=\{\balpha_1,\ldots, \balpha_n\}, \hbox{ where }
 \balpha_i \in \Z[x]^m,i=1,\ldots,n,
\end{equation}
and  $\T =(t_1, \ldots, t_m)$ a set of $\sigma$-indeterminates.
We define the following rational $\sigma$-morphism
\begin{equation}\label{eq-torus}
\phi:(\AH^{*})^{m} \longrightarrow (\AH^{*})^{n},  \T \mapsto
\T^{\falpha} = (\T^{\balpha_1},  \ldots, \T^{\balpha_n}).
\end{equation}

Define the functor $T_{\falpha}$ from $\mathscr{E}_{\F}$ to
$\mathscr{E}_{\F}^n$ with $T_{\falpha}(\E)=\Im(\phi_{\E})$ which is
called a {\em quasi $\sigma$-torus} with defining vector $\falpha$.
For each $\E \in \mathscr{E}_{\F}$, $T_{\falpha}(\E)$ has a group
structure with component wise multiplication:
$\T_1^{\falpha}\cdot\T_2^{\falpha} =(\T_1\cdot\T_2)^{\falpha}$ where
$\T_1,\T_2 \in (\E^{*})^{m}$ and $\balpha\in\Zx^m$. If
$\E_{1},\E_{2} \in \mathscr{E}_{\F}$, $\E_{1} \subseteq \E_{2}$,
then $T_{\falpha}(\E_{1})$ is a subgroup of $T_{\falpha}(\E_{2})$.

We now define the toric $\sigma$-variety.
\begin{definition}
A $\sigma$-variety over the $\sigma$-field $\Q$ is called {\em
toric} if it is the Cohn closure of a quasi $\sigma$-torus
$T_{\falpha}$ in $\AH^n$, where  $\falpha$ is given in \bref{eq-A}.
More precisely,  let
 \begin{equation}\label{eq-tv0}
 \I_{\falpha}=\{f\in \Q\{\Y\}\,|\,f(\T^{\balpha_1},  \ldots,
\T^{\balpha_n})=0\}.
 \end{equation}
Then the toric  $\sigma$-variety defined by $\falpha=\{\balpha_{1},
\ldots, \balpha_{n}\}$ is $\XB_{\falpha}=\V(\I_{\falpha})$.
$A=[\balpha_{1}, \ldots, \balpha_{n}]_{m\times n}$ is called the
defining matrix of $\XB_{\falpha}$.
\end{definition}
We make the following reasonable assumption: $A$ does not contain a
zero row, or equivalently, every $t_i$ appears effectively in some
$\T^{\balpha_k}$. Also, no $\balpha_i$ is the zero vector.

\begin{lemma}\label{lm-tv10}
Let $L_{\falpha}$ be the $\Zx$-lattice generated by $\falpha$ given
in \bref{eq-A}. Then $\XB_{\falpha}$ is an irreducible
$\sigma$-variety of dimension $\rk(L_{\falpha})$.
\end{lemma}
\proof It is clear that $\T^{\falpha}$ in \bref{eq-torus} is a
generic point of $\I_{\falpha}$ in \bref{eq-tv0}. Then
$\I_{\falpha}$ is a reflexive prime $\sigma$-ideal. By Theorem 3.20
of \cite{dd-sres,dd-sres1}, $\I_{\falpha}$ is of dimension $\dtrdeg
 \Q\langle\T^{\falpha}\rangle/\Q =
 \rk(A)=\rk(L_{\falpha})$, where $A=[\balpha_{1}, \ldots, \balpha_{n}]_{m\times n}$
is the defining matrix of $\XB_{\falpha}$.\qedd

Let $[y_1-\T^{\balpha_1},\ldots,y_n-\T^{\balpha_n}]$ be the
$\sigma$-ideal generated by $y_i-M_i,i=1,\ldots,n$ in
$\Q\{\Y,\T^{\pm}\}$. Then it is easy to check
 \begin{equation}\label{eq-tv01}
  \I_{\falpha}=[y_1-\T^{\balpha_1},\ldots,y_n-\T^{\balpha_n}]\cap\Q\{\Y\}.
 \end{equation}
Alternatively, let
$[\T^{\balpha_1^+}y_1-\T^{\balpha_1^-},\ldots,\T^{\balpha_n^+}y_n-\T^{\balpha_n^-}]$
be a $\sigma$-ideal in $\Q\{\Y,\T\}$. Then
 \begin{equation}\label{eq-tv02}
  \I_{\falpha}=[\T^{\balpha_1^+}y_1-\T^{\balpha_1^-},\ldots,\T^{\balpha_n^+}y_n-\T^{\balpha_n^-}]:\mb_{\T}\cap\Q\{\Y\}
 \end{equation}
 where  $\mb_{\T}$ is the multiplicative set generated by $t_i^{x^j}$
for $i=1, \ldots, m,  j\in\N$.

\begin{remark}
With formulas \bref{eq-tv01} and \bref{eq-tv02}, a characteristic
set for $\I_{\falpha}$ can be computed with the characteristic set
method in \cite{gao-dcs}. More efficient methods to compute a
characteristic set for $\I_{\falpha}$ will be given in Section
\ref{sec-tv2}.
\end{remark}

\begin{example}\label{ex-tv1}
Let $M_2=\left[\begin{array}{llll}
 2       & x-1      & 0      & 0   \\
 0       & 0          & 2      & x-1       \\
\end{array}\right]$
be the matrix from Example \ref{ex-L1} and $\falpha$ the set of
columns of $M_2$. Let
$\I_{1}=[y_1-t_1^2,y_2-t_1^{x-1},y_3-t_2^2,y_4-t_2^{x-1}]$.
By \bref{eq-tv01}, $\I_{\falpha} = \I_1\cap\Q\{y_1,y_2,y_3,y_4\}$.
With the characteristic set method \cite{gao-dcs}, under the
variable order $y_2 < y_4<y_1 < y_3< t_1 < t_2$, a characteristic
set of $\I_1$ is $y_1y_2^{2}-y_1^x, y_{3}y_{4}^{2}-y_{3}^x,
y_1-t_1^2,t_1y_2-t_1^{x},y_3-t_2^2, t_2y_4-t_2^{x}$.
Then
$$\I_{\falpha} = \sat(y_1y_2^{2}-y_1^x,
y_{3}y_{4}^{2}-y_{3}^x)=[y_1y_2^{2}-y_1^x,
y_{3}y_{4}^{2}-y_{3}^x].$$
\end{example}

The following example shows that some $y_i$ might not appear
effectively in $\I_{\falpha}$.
\begin{example}\label{ex-tv2}
Let $\falpha= \{[1,1]^\tau,[x,x]^\tau,[0,1]^\tau\}$.
By \bref{eq-tv01}, $\I_{\falpha} =
[y_1-t_1t_2,y_2-t_1^{x}t_2^x,y_3-t_2]$ $\cap\Q\{y_1,y_2,y_3\}=[y_1^x-y_2]$
and $y_3$ does not appear in $\I_{\falpha}$.
\end{example}

The following lemma shows that quasi $\sigma$-tori of two sets of
generators of $L_{\falpha}$ are isomorphic.
\begin{lemma}\label{lm-tv11}
Let $\falpha=\{\balpha_{1}, \ldots, \balpha_{n}\}$ and
$\fbeta=\{\bbeta_{1}, \ldots, \bbeta_{s}$\} be two sets of
generators for $L$. Then  $T_{\falpha}$ and $T_{\fbeta}$ are
isomorphic as groups.
\end{lemma}
\proof Let $A_{m\times n}$ and $B_{m\times s}$ be the matrix
representations for the two sets of generators for $L$.
Then there exist matrices $M=(m_{ij})\in\Z[x]^{s \times n}$ and
$N=(n_{ij})\in\Z[x]^{n \times s}$ such that $A=BM,  B=AN$. Hence
$A=ANM,  B=BMN$.
Define two maps
 $\theta_1: \AH^n \Rightarrow \AH^s$ and
 $\theta_2: \AH^s \Rightarrow \AH^n$ by
 $\theta_1(y_1,\ldots,y_n) = (\prod_{i=1}^n y_i^{n_{i1}},\ldots,
 \prod_{i=1}^n y_i^{n_{is}})$ and
 $\theta_2(x_1,\ldots,x_s) = (\prod_{i=1}^s x_i^{m_{i1}},\ldots,
 \prod_{i=1}^s x_i^{m_{in}})$.
Then it is clear that $\theta_1(T_{\falpha})\subset T_{\fbeta}$ and
$\theta_2(T_{\fbeta})\subset T_{\falpha}$. From $A=ANM$ and $B=BMN$,
it can be checked that $\theta_1 \circ \theta_2=id$ on quasi
$\sigma$-torus $T_{\fbeta}$ and $\theta_2 \circ \theta_1=id$ on
quasi $\sigma$-torus $T_{\falpha}$.
It is easy to check $\theta_1$ and $\theta_2$ are group
homomorphisms. Therefore, the quasi tori $T_{\falpha}$ and
$T_{\fbeta}$ are isomorphic as groups.\qedd

In the rest of this section, we will give a description for the
coordinate ring of a toric $\sigma$-variety.
\begin{definition}
$M \subseteq \Z[x]^{m}$ is called an {\em affine $\N[x]$-module} if
there exists a $\fbeta=\{{\bbeta_{1}}, \ldots, {\bbeta_{s}}\}$
$\subset \Z[x]^{m}$ such that $M =\N[x](\fbeta)=\{\sum_{i=1}^s
a_{i}{\bbeta_{i}}\,|\,a_{i}\in \N[x]\}$.
\end{definition}

Given an affine $\N[x]$-module $M \subseteq \Z[x]^{m}$ and a set of
$\sigma$-indeterminates $\T=\{t_1,\ldots,t_m\}$, define the
corresponding {\em affine $\sigma$-algebra}
 $$\F\{M\}=\{\sum_{\bbeta \in M}a_{\bbeta}\T^{\bbeta}\,|\, a_{\bbeta}\in \F,   a_{\bbeta}\ne0 \hbox{ for
finitely many } \bbeta\}.$$
It is easy to check that $\F\{M\}=\F\{N[x]\fbeta\}$ is a
$\sigma$-algebra over $\F$. If $M=\N[x]\fbeta$ with
$\fbeta=\{{\bbeta_{1}}, \ldots, {\bbeta_{s}}\}$, then
$\F\{M\}=\F\{\T^{{\bbeta_{1}}}, \ldots, \T^{{\bbeta_{s}}}\}$.

\begin{example}
As an affine $\N[x]$-module, $\N[x]^{m}$ gives the
$\sigma$-polynomial ring in $\T$: $\F\{\N[x]^{m}\}$ $=\F\{\T\}$.
Let ${\beps}_{1}, \ldots, {\beps}_{m}$ be the standard basis of
dimension $m$. Then $\Z[x]^{m}$ is generated by $A=\{\pm
{\beps}_{1}, \ldots, \pm {\beps}_{m}\}$ as an affine $\N[x]$-module,
and the corresponding affine $\sigma$-algebra is the Laurent
$\sigma$-polynomial ring $\F\{\Z[x]^{m}\}=\F\{\T^{\pm}\}$.
\end{example}

Given an affine $\N[x]$-module $M \subseteq \Z[x]^{n}$, we define
{\em $\Spec^{\sigma}(\F\{M\})$ } to be the set of reflexive prime
$\sigma$-ideals of $\F\{M\}$ \cite{wibmer}.

\begin{theorem}\label{th-tv31}
$X$ is a toric $\sigma$-variety if and only if there exists an
affine $\N[x]$-module $M$ such that $X\cong\Spec^{\sigma}(\Q\{M\})$.
Equivalently, the coordinate ring of $X$ is $\Q\{M\}$.
\end{theorem}
\proof Let $X=\XB_{\falpha}$ and $\IH(X) = \I_{\falpha}$, where
$\falpha$ is defined in \bref{eq-A} and $\I_{\falpha}$ is defined in
\bref{eq-tv0}. Let $M=\N[X]{\falpha}$ be the affine $\N[x]$-module
generated by $\falpha$.
Define the following morphism of $\sigma$-rings
$$
\theta: \Q\{\Y\}\longrightarrow \Q\{M\}, \hbox{ where }
\theta(y_{i})= \T^{\balpha_i},i=1,\ldots,n.
$$
The map $\theta$ is surjective by the definition of $\Q\{M\}$. If
$f\in \ker(\theta)$,  then
 $f(\T^{\balpha_1}, \ldots, \T^{\balpha_i})=0$, which is
equivalent to $f\in \I_\falpha$. Then, $\ker(\theta) =\I_{\falpha}$
and $\Q\{\Y\}/{\I_\falpha} \cong \Q\{M\}$. Therefore
$X\cong\Spec^{\sigma}(\Q\{\Y\}/\I_{\falpha})=\Spec^{\sigma}(\Q\{M\})$.

If $X\cong\Spec^{\sigma}(\Q\{M\})$,  where $M \subseteq \Z^{m}[x]$
is an affine $\N[x]$-module, and $M=\N[x](\balpha_{1}, \ldots,$
$\balpha_{n})$ for $\balpha_{i}\in M$.
Let $\XB_{\falpha}$ be the toric $\sigma$-variety defined by
$\falpha=\{\balpha_{1}, \ldots, \balpha_{n}\}$. Then as we just
proved, the coordinate ring of $X$ is isomorphic to $\Q\{M\}$. Then
$X \cong\XB_{\falpha}$.\qedd

\subsection{Toric $\sigma$-ideal}\label{sec-tv2}
In this section,  we will show that $\sigma$-toric varieties are
defined exactly by toric $\sigma$-ideals which are in a one to one
correspondence with toric $\Zx$-lattices.
%

\begin{definition}
A $\Z[x]$-lattice $L\subseteq \Z[x]^{n}$ is called {\em toric} or
{\em $\Z[x]$-saturated} if for any $p\in\Z[x]$ and $\f \in
\Z[x]^{n}$, $p\f \in L$ implies $\f \in L$.
Let $L$ be a toric $\Z[x]$-lattice in $\Z[x]^{n}$ and $\rho_L$ the
trivial partial character defined on $L$. Then the $\sigma$-ideal
 $$I^{+}(\rho_L)=[\Y^{\f^+}-\Y^{\f^-}\,|\,\f\in L]
    =[\Y^{\a}-\Y^{\b}\,|\,\a-\b\in L, \a, \b \in \N[x]^{n}]$$
 is called a {\em toric $\sigma$-ideal}.
\end{definition}
We list several properties for toric $\sigma$-ideals.
 \begin{itemize}
\item
By Theorem \ref{th-pr1}, a toric $\sigma$-ideal is a reflexive prime
$\sigma$-ideal of dimension $n-\rk(L)$.
\item
By Corollary \ref{lm-nl2}, toric $\sigma$-ideals in $\F\{\Y\}$ are
in a one to one correspondence with toric lattices in $\Zxn$.
\end{itemize}
In Section \ref{sec-alg3}, an algorithm will be given to check
whether a $\Zx$-lattice is toric.
In this section, we will prove the following result which can be
deduced from Lemmas \ref{lm-tvi1} and \ref{lm-tvi4}.
\begin{theorem}\label{th-tvi}
A $\sigma$-variety $V$ over $\Q$ is toric if and only if $\IH(V)$ is
a toric $\sigma$-ideal.
\end{theorem}
%

The following lemma shows that the defining ideal of a toric
$\sigma$-variety is a toric $\sigma$-ideal.
\begin{lemma}\label{lm-tvi1}
Let $\XB_{\falpha}$ be a toric $\sigma$-variety and $A$ the defining
matrix of $\XB_{\falpha}$. Then $\IH(\XB_{\falpha})$ is the toric
$\sigma$-ideal whose support lattice is $\ker(A)$.
\end{lemma}
\proof
%
Let $A=[\balpha_1,\ldots,\balpha_n]_{m\times n}$ be the defining
matrix of the toric $\sigma$-variety and
 $$K_{A}=\ker(A)=\{\f \in \Z[x]^{n}\,|\,A\f=0\}.$$
As a kernel, $K_A$ is clearly a toric $\Zx$-lattice. Then it
suffices to show that $\IH(\XB_{\bf\alpha})=\I_{\falpha} =
\I^{+}(\rho_{K_A})$, where $\rho_{K_A}$ is the trivial partial
character defined over $K_A$ and $\I_{\falpha}$ is defined in
\bref{eq-tv0}.

For $\f \in K_{A}$, we have
 $(\Y^{\f}-1)(\T^{\falpha}) = (\T^{\falpha})^{\f} - 1 =
 \T^{A\f}-1 = 0$.
As a consequence, $(\Y^{\f ^{+}}-\Y^{\f^{-}})(\T^{\falpha})=0$ and
$\Y^{\f ^{+}}-\Y^{\f^{-}} \in \I_{\falpha}$. Since
$\I^{+}(\rho_{K_A})$ is generated by $\Y^{\f ^{+}}-\Y^{\f^{-}}$ for
$\f\in K_A$, we have $ \I^{+}(\rho_{K_A})\subset\I_{\falpha}$.

Consider a well order for elements in $\N^{n}[x]$, which leads to a
well order for $\{\Y^\f, \f\in \N[x]^n\}$ as well as an order for
$\F\{\Y\}$ by comparing the largest $\sigma$-monomial in a
$\sigma$-polynomial.
We will prove $\I_{\falpha} \subset\I^{+}(\rho_{K_A})$. Assume the
contrary, and let $f=\Sigma_{i}a_{i}\Y^{\f_{i}}\in\I_{\falpha}$ be a
minimal element in $\I_{\falpha}\setminus\I^{+}(\rho_{K_A})$ under
the above order. Let $a_0\Y^{\g}$ be the biggest $\sigma$-monomial
in $f$. From $f\in\I_{\falpha}$, we have $f(\T^{\falpha})=0$. Since
$\Y^{\g}(\T^{\falpha})=\T^{A\g}$ is a $\sigma$-monomial about $\T$,
there exists another $\sigma$-monomial $b_0\Y^{\h}$ in $f$ such that
$\Y^{\h}(\T^{\falpha})=\Y^{\g}(\T^{\falpha})$. As a consequence,
$(\Y^{\g}-\Y^{\h})(\T^{\falpha})=\T^{A\h}(\T^{A(\g-\h)}-1)=0$, from
which we deduce $\g-\h\in K_A$ and hence $\Y^{\g}-\Y^{\h} \in
\I_{\falpha}\cap\I^{+}(\rho_{K_A})$.
Then $f-a_{0}(\Y^{\g}-\Y^{\h})\in
\I_{\falpha}\setminus\I^{+}(\rho_{K_A})$, which contradicts to the
minimal property of $f$, since  $f-a_{0}(\Y^{\g}-\Y^{\h})<f$.\qedd

To prove the converse of Lemma \ref{lm-tvi1}, we first introduce a
new concept for $\Zx$-lattices. Let $L\subset \Zxn$ be a $\Zx$-lattice. Define the {\em orthogonal complement} of $L$ to be
 $$ L^C = \{\f\in\Zxn\,|\, \forall \g\in L, \langle \f,\g\rangle=0 \}$$
where $\langle \f,\g\rangle=\f^\tau\cdot \g$ is the dot product of
$\f$ and $\g$. It is easy to show that
\begin{lemma}\label{lm-tvi2}
Let $A_{n\times r}$ be a matrix representation for $L$. Then $L^C
=\ker(A^\tau)= \{\f\in\Zxn\,|\, A^\tau \f=0\}$ and hence $\rk(L^C) =
n - \rk(L)$.
\end{lemma}
Furthermore, we have
\begin{lemma}\label{lm-tvi3}
If $L$ is a toric $\Zx$-lattice, then $L= (L^C)^C$.
\end{lemma}
\proof It is easy to see $L\subset (L^C)^C$. Let $r=\rk(L)$. By
Lemma \ref{lm-tvi2}, $\rk((L^C)^C) = n - (n-r) = r=\rk(L)$.
Let $\K = \Q(x)$. In $\K^n$, the $\Zx$-lattices $L$ and $(L^C)^C$
become vector spaces $\widetilde{L}$ and $\widetilde{(L^C)^C}$ with
$\dim(\widetilde{L})=\rk(L)=\dim(\widetilde{(L^C)^C})$. Since
$\widetilde{L}\subset\widetilde{(L^C)^C}$ and
$\dim(\widetilde{L})=\dim(\widetilde{(L^C)^C})$, we have
$\widetilde{L}=\widetilde{(L^C)^C}$.
Let $\f\in(L^C)^C$. Then $\f\in\widetilde{(L^C)^C}=\widetilde{L}$,
which means that there exists a $p\in\Z[x]$ such that $p\f\in L$.
Since $L$ is toric, we have $\f\in L$ and the lemma is proved.\qedd

From Lemmas \ref{lm-free}, \ref{lm-tvi2}, and \ref{lm-tvi3}, we have
\begin{cor}\label{cor-tvi1}
A toric $\Zx$-lattice is a free $\Zx$-module.
\end{cor}

The following lemma shows that the inverse of Lemma \ref{lm-tvi1} is
also valid.
\begin{lemma}\label{lm-tvi4}
If $\I$ is a toric $\sigma$-ideal in $\Q\{\Y\}$, then $\V(\I)$ is a
toric $\sigma$-variety.
\end{lemma}
\proof Since $\I$ is a toric $\sigma$-ideal, $\I$ is reflexive prime
and $\L(\I)$ is a toric lattice.
Let $\{\f_{1}, \ldots, \f_{s}\}$ be a set of generators of $\L(\I)$,
$F=[\f_{1}, \ldots, \f_{s}]_{n\times s}\in\Zx^{n\times s}$ a matric
representation for $\L(\I)$, and $K_{F}=\{X\in \Z[x]^{n}|F^\tau
X=0\} = \L(\I)^C$. Then by Corollary \ref{cor-tvi1} and Lemma
\ref{lm-tvi2}, $K_{F}$ is a free $\Zx$-module and hence has a basis
$\{\h_{1}, \ldots, \h_{n-r}\}$, where $r=\rk(F)=\rk(\L(\I))$.
Let $H=[\h_{1}, \ldots, \h_{n-r}]_{n\times(n-r)}$ be the matrix with
$\h_i$ as the $i$-th column and $\falpha=\{\balpha_{1}, \ldots,
\balpha_{n}\}$ the rows of $H$.
Consider the toric $\sigma$-variety  $\XB_{\falpha}$ defined by the
following quasi $\sigma$-torus
 $$\phi: (\AH^{\ast})^{n-r}\longrightarrow \AH^{n} $$
where $\phi(\T)=(\T^{\balpha_{1}},\ldots,\T^{\balpha_{n}})$ for
$\T=(t_1,\ldots,t_{n-r})$. To prove the lemma, it is suffices to
show $\XB_\falpha=\V(\I)$ or $\I_\falpha=\I$

By Lemma \ref{lm-tvi1}, $\I_\falpha$ is toric. Since both
$\I_\falpha$ and $\I$ are toric, to prove $\I_\falpha=\I$, we need
only to show $\L(\I_\falpha)=\L(\I)$.
By Lemma \ref{lm-tvi2}, $K_F^C=\ker(H^\tau)$.
Then, by Lemma \ref{lm-tvi1}, $\L(\I_\falpha) = \ker(H^\tau) = K_F^C
=(\L(\I)^C)^C$. By Lemma \ref{lm-tvi3}, $\L(\I_\falpha) = \L(\I)$
and the lemma is proved.\qedd

It should be noted that the proofs of Lemma \ref{lm-tvi1} and  Lemma
\ref{lm-tvi2} also give algorithms to compute the defining ideal for
a toric $\sigma$-variety and the defining matrix for the toric
variety defined by a toric $\sigma$-ideal.
In other words,  a toric $\sigma$-variety is a $\sigma$-variety
which has a Laurent $\sigma$-monomial parametrization, and the
proofs of Lemma \ref{lm-tvi1} and Lemma \ref{lm-tvi2} give
implicitization and parametrization algorithms for these kind of
unirational $\sigma$-varieties \cite{gao-im}.
\begin{example}\label{ex-tvi11}
Continue from Example \ref{ex-tv1}
%
%
Let $\f_1=(1-x,2,0,0)^\tau,\f_2=(0,0,1-x,2)^\tau$. Then
$K_{M_2}=\ker(M_2) = (\f_1,\f_2)\subseteq \Z^{4}[x]$. By Lemma
\ref{lm-tvi1} and Theorem \ref{th-nl3}, with the variable order $y_2
< y_4<y_1 < y_3$, we have
 $\I_{\falpha}=\sat(y_1y_2^{2}-y_1^x, y_{3}y_{4}^{2}-y_{3}^x)=[y_1y_2^{2}-y_1^x, y_{3}y_{4}^{2}-y_{3}^x].$

Conversely, let $B=[\f_1,\f_2]_{4\times 2}$ be the support lattice
of $\IH(\XB_{\falpha})$. Then $M_2^\tau$ is the defining matrix for
$\ker(B^\tau)$. By Lemma \ref{lm-tvi4}, $M_2$ is the defining matrix
for the toric $\sigma$-variety $\XB_{\falpha}$.
\end{example}

For any $\Z[x]$-lattice $L\subseteq \Z[x]^{n}$, the {\em
$\Z[x]$-saturation} of $L$ is defined as
 $$\sat_{\Z[x]}(L)=\{ \f\in \Zxn\, |\, \exists g\in\Z[x], s.t.\, g\f \in
L\}.$$
Similar to Corollary \ref{cor-ld0}, it can be shown that
$\rk(\sat_{\Z[x]}(L)) = \rk(L)$. An algorithm to compute the $\Zx$
saturation is given in Section \ref{sec-alg3}.
\begin{example}\label{ex-ig11}
Continue from Example \ref{ex-ig}. Let $L=([1-x,x-1]^\tau)$. Then
$\sat_{\Z[x]}(L) = ([1,-1]^\tau)$.
The $\sigma$-ideal
$\I=\I^{+}(\rho_L)=\sat(y_1y_2^x-y_1^xy_2)=[y_1^{x^i}
y_2^{x^j}-y_1^{x^j}y_2^{x^i}]$ is reflexive prime but not toric. The
minimal toric $\sigma$-ideal containing $\I$ is $[y_2-y_1]$ with
lattice support $([1,-1]^\tau)$.
This is quite different from the algebraic case \cite{es-bi}, where
all prime binomial ideals generated by binomial of the form
$(\Y^\f-\Y^\g)$ are toric.
\end{example}

%
%
%
%
%
%

\subsection{$\sigma$-torus and toric $\sigma$-variety in terms of group action}
In this section,  we first define the $\sigma$-torus and prove its
properties and then give a abstract description for $\sigma$-toric
varieties in terms of group actions.

Let $T_{\falpha}$ be a quasi $\sigma$-torus and $\XB_{\falpha}$ the
toric $\sigma$-variety defined by ${\falpha}\subset\Zx^m$. In the
algebraic case, $T_{\falpha}$ is a variety, that is, $T_{\falpha}=
\XB_{\falpha}\cap(\CH^*)^m$, where $\CH$ is the field of complex
numbers and $\CH^*=\CH\setminus\{0\}$. The following example shows
that this is not valid in the difference case.
\begin{example}\label{ex-dt1}
In Example \ref{ex-tv1},  $\XB_{\falpha} = \V({\{y_1y_2^{2}-y_1^x,
y_{3}y_{4}^{2}-y_{3}^x\}})$.
Let $P= (-1,1,-1,-1)$. Then $P\in \XB_{\falpha}$. On the other hand,
assume $P\in T_{\falpha}$ which means
$((t_1)^{2},(t_1)^{x-1},(t_2)^{2},(t_2)^{x-1})$ $=(-1,1,-1,-1)$ or the
$\sigma$-equations
$t_1^{2}+1=0,t_1^{x}-t_1=0,t_2^{2}+1=0,t_2^{x}+t_2=0$ have a
solution. From Example \ref{ex-31}, this is impossible. That is,
$T_{\falpha}\varsubsetneq \XB_{\falpha}\cap(\CH^*)^m$.
\end{example}

In order to define the $\sigma$-torus, we need to introduce the
concept of Cohn $*$-closure.  $(\AH^*)^n$ is isomorphic to the
$\sigma$-variety defined by
$\I_0=[y_1z_1-1,\ldots,y_nz_n-1]\subset\F\{\Y,\Z\}$ in $(\AH)^{2n}$.
Furthermore, $\sigma$-varieties in $(\AH^*)^n$ are in a one to one
correspondence with affine $\sigma$-varieties contained in
$\V({\I_0})$ via the map $\theta:(\AH^*)^n\Rightarrow (\AH)^{2n}$
defined by $\theta(a_1,\ldots,a_n)
=(a_1,\ldots,a_n,a_1^{-1},\ldots,a_n^{-1})$.
Let $V\subset(\AH^*)^n$ and $V_1$ the Cohn closure of $\theta(V)$ in
$(\AH)^{2n}$. Then the $\theta^{-1}(V_1)$ is called the {\em Cohn
$*$-closure} of $V$.

Example \ref{ex-dt1} gives the motivation for the following
definition.
\begin{definition}
A $\sigma$-torus is an $\sigma$-variety which is isomorphic to the
Cohn $*$-closure of a quasi $\sigma$-torus in $(\AH^{*})^{n}$.
\end{definition}

\begin{theorem}\label{lm-tv41}
Let $T_{\falpha}$ be the quasi $\sigma$-torus defined in
\bref{eq-torus} by $\falpha$, $T_{\falpha}^*$ the Cohn $*$-closure
of $T_{\falpha}$ in $(\AH^*)^n$, and $\I_{\falpha}$ defined in
\bref{eq-tv0}. Then $T_{\falpha}^*$ is isomorphic to
$\Spec^{\sigma}(\F\{\Y,\Z\}/\widetilde{\I}_{\falpha})$ where
$\Z=\{z_1,\ldots,z_n\}$ is a set of $\sigma$-indeterminates and
$\widetilde{\I}_{\falpha} = [\I_{\falpha},y_1z_1-1,\ldots,y_nz_n-1]$
in $\F\{\Y,\Z\}$.
\end{theorem}
\proof Use the notations just introduced.
Let $\widetilde{T}_{\falpha} = \theta(T_{\falpha})\subset\AH^{2n}$
and $\widetilde{T}_{\falpha}^*$ the Cohn closure of
$\widetilde{T}_{\falpha}$ in $\AH^{2n}$. Then $T_{\falpha}^* =
\theta^{-1}(\widetilde{T}_{\falpha}^*)$ is the Cohn $*$-closure of
$T_{\falpha}$ in $(\AH^*)^n$. Since $\theta$ is clearly an
isomorphism between $\widetilde{T}_{\falpha}^*$ and $T_{\falpha}^*$,
it suffices to show that
$\IH(\widetilde{T}_{\falpha}^*)=\widetilde{\I}_{\falpha}$.

We have  $\IH(\widetilde{T}_{\falpha}^*) = \{p\in\F\{\Y,\Z\}\,|\,
p(\T^{\balpha_1},
\ldots,\T^{\balpha_n},\T^{-\balpha_1},\ldots,\T^{-\balpha_n})=0\}$.
It is clear that $\widetilde{\I}_{\falpha}\subset
\IH(\widetilde{T}_{\falpha}^*)$.
If $f\in \IH(T_{\falpha}^*)$, eliminate $z_1,\ldots,z_n$ from $f$ by
taking the $\sigma$-remainder
$f_1=\prem(f,\{y_1z_1-1,\ldots,y_nz_n-1\})=\prod_{i=1}^n y_i^{t_i} f
- \sum_{i=1}^n\sum_{j=0}^{l_i} g_{ij}(y_iz_i-1)^{x^{s_{ij}}} \in
\F\{\Y\}\cap\IH(T_{\falpha}^*)$, where $t_i\in\N[x], s_{ij}\in\N,
g_{ij}\in\F\{\Y,\Z\}$.
Substituting $y_i$ by $\T^{\balpha_i}$ and $z_i$ by
$\T^{-\balpha_i}$, we have $f_1(\T^{{\falpha_1}}, \ldots,
\T^{{\balpha_n}})=0$, which is equivalent to $f_1\in \I_{\falpha}$.
Then $\prod_{i=1}^n y_i^{t_i} f\in\widetilde{\I}_{\falpha}$ and
hence $\prod_{i=1}^n z_i^{t_i}y_i^{t_i} f = (z_iy_i)^{t_i} f = f
+f_0 \in\widetilde{\I}_{\falpha}$, where $f_0 \in\I_0$. Thus
$f\in\widetilde{\I}_{\falpha}$.\qedd

\begin{cor}\label{cor-tv41}
Use the notations in Lemma \ref{lm-tv41}. $T^*_{\falpha}$ is
isomorphic to $\Spec^{\sigma}(\Q\{\Y^{\pm}\}/\I_{\falpha})$.
\end{cor}
\proof The isomorphism is given by $y_i \Rightarrow y_i$ and $z_i
\Rightarrow y_i^{-1}$ for $i=1,\ldots,n$.\qedd

\begin{cor}\label{cor-tv42}
Let $T_{\falpha}^*$ and $\XB_{\falpha}$ be the $\sigma$-torus and
the toric $\sigma$-variety defined by ${\falpha}$, respectively.
Then $T_{\falpha}^* = \XB_{\falpha}\cap (\AH^*)^n$. Equivalently,
$T_{\falpha}^*$ is a Cohn open set in $\XB_{\falpha}$.
\end{cor}
\proof From Lemma \ref{lm-tv41} and the fact $\I_{\falpha} =
\IH(\XB_{\falpha})$, we have $T_{\falpha}^* = \XB_{\falpha}\cap
(\AH^*)^n$.
Then $\XB_{\falpha}\setminus T_{\falpha}^* = \cup_{i=1}^n
\V({\{\I_{\falpha}, y_i\}})$. Since $\I_{\falpha}$ is a reflexive
and prime $\sigma$-ideal not containing any $y_i$,
$\V({\{\I_{\falpha}, y_i\}})$ is a proper sub-$\sigma$-variety of
$\XB_{\falpha}$, or equivalently, $T_{\falpha}^*$ is an open set of
$\XB_{\falpha}$.\qedd

\begin{cor}\label{cor-tv43}
A $\sigma$-torus $T^*$ is a toric $\sigma$-variety and has a group
structure. As a consequence,  $T^* = t\cdot T^*$ for any $t\in\T^*$.
\end{cor}
\proof From Lemma \ref{lm-tv41}, it suffices to prove the corollary
for $\widetilde{T}_{\falpha}^*=\V({\{\widetilde{\I}_{\falpha}\}})$.
By Theorem \ref{th-tvi}, $\I_\falpha\subset\Q\{\Y\}$ is toric. Then,
it is apparent that $\widetilde{\I}_{\falpha} =
[\I_{\falpha},y_1z_1-1,\ldots,y_nz_n-1]$ is also toric.
To show that $\widetilde{T}_{\falpha}^*$ has a group structure, let
$a, b \in \widetilde{T}_{\falpha}^{*}$. Since $y_iz_i-1\in
\IH(\widetilde{T}_{\falpha}^{*}),i=1,\ldots,n$, $a^{-1}$ is well
defined.
Then for each $\sigma$-binomial
$\Y^{\f}-\Y^{\g}\in\IH(\widetilde{T}_{\falpha}^{*})$, we have
$(\Y^{\f}-\Y^{\g})(a)=0,(\Y^{\f}-\Y^{\g})(b)=0$. Then
$(\Y^{\f}-\Y^{\g})(a b)=0$ and $(\Y^{\f}-\Y^{\g})(a^{-1})=0$. It
implies that the $\sigma$-torus is a group.
Since $\widetilde{T}_{\falpha}^*$ is a group, for $t\in
\widetilde{T}_{\falpha}^*$, we have $t\cdot
\widetilde{T}_{\falpha}^* \subset \widetilde{T}_{\falpha}^*$ and
$t^{-1}\cdot \widetilde{T}_{\falpha}^* \subset
\widetilde{T}_{\falpha}^*$ which implies $t\cdot
\widetilde{T}_{\falpha}^* =\widetilde{T}_{\falpha}^*$. \qedd

Similar to Theorem \ref{th-tv31}, we have the following result.
\begin{theorem}\label{th-tv41}
$T^*$ is a $\sigma$-torus if and only if there exists a $\Zx$-lattice $L$ in $\Zx^m$ such that $T^*\cong\Spec^{\sigma}(\Q\{L\})$,
where
  $$\Q\{L\}=\{\sum_{\bbeta \in L}a_{\bbeta}\T^{\bbeta}\,|\, a_{\bbeta}\in \Q,
   a_{\bbeta}\ne0 \hbox{ for finitely many } \bbeta\}.$$
\end{theorem}
\proof Let $T_{\falpha}$ be the quasi $\sigma$-torus defined by
$\falpha$, $T_{\falpha}^*$ the Cohn $*$-closure of $T_{\falpha}$ in
$(\AH^*)^n$, $L$ the $\Zx$-lattice generated by $\falpha$, and
$\I_{\falpha}$ defined in \bref{eq-tv0}.
For a set of $\sigma$-indeterminates $\Z=\{z_1,\ldots,z_n\}$, define
the following morphism of $\sigma$-rings
$$
\theta: \Q\{\Y,\Z\}\longrightarrow \Q\{L\}, \hbox{ where }
 \theta(y_{i})= \T^{{\balpha_i}},
 \theta(z_{i})= \T^{-{\balpha_i}}, i=1,\ldots,n.
$$
It is easy to check that the map $\theta$ is surjective. By Theorem
\ref{lm-tv41}, $\widetilde{\I}_{\falpha}=\ker(\theta) =
[\I_\alpha,y_1z_1-1,\ldots,y_nz_n-1]$ which is the defining ideal
for $\widetilde{T}_{\falpha}^*$.
Since $\widetilde{T}_{\falpha}^*$ is isomorphic to $T_{\falpha}^*$,
we have $\Q\{\Y,\Z\}/\widetilde{\I}_{\falpha} \simeq \Q\{L\}$ and
$T_{\falpha}^*\cong\Spec^{\sigma}(\Q\{L\})$. We thus proved that a
$\sigma$-torus is isomorphic to $\Spec^{\sigma}(\Q\{L\})$.

If $T^*=\Spec^{\sigma}(\Q\{L\})$,  where $L=\Z[x](\balpha_{1},
\ldots, \balpha_{n})$ is a $\Zx$-lattice for $\balpha_{i}\in \Zx^m$.
Let $T_{\falpha}^*$ be the $\sigma$-torus defined by
$\falpha=\{\balpha_{1}, \ldots, \balpha_{n}\}$. Then as we just
proved, the coordinate ring of $T^*_{\falpha}$ is isomorphic to
$\Q\{L\}$. Then $T^*$ is isomorphic to $T^*_{\falpha}$.\qedd
%

As a consequence, we can prove that the tori defined by two sets of
generators of the same $\Zx$-lattice are isomorphic.
\begin{cor}\label{cor-tv25}
Let $\falpha=\{\balpha_{1}, \ldots, \balpha_{n}\}$ and
${\fbeta}=\{{\bbeta}_{1}, \ldots, {\bbeta}_{s}$\} be two sets of
generators for the same $\Zx$-lattice. Then $T_{\falpha}^*$ and
$T_{{\fbeta}}^*$ are isomorphic as $\sigma$-varieties and groups.
\end{cor}
\proof By Theorem \ref{th-tv41}, it suffices to show that
$\Q\{\Z[x]{\falpha}\}$ is isomorphic to $\Q\{\Z[x]{\fbeta}\}$ as
$\sigma$-rings.
Use notations introduced in the proof of Lemma \ref{lm-tv11}.
Define
 $$\phi:\Q\{\Z[x]{\falpha}\}\Rightarrow\Q\{\Z[x]{\fbeta}\}
 \hbox{ and }
  \theta:\Q\{\Z[x]{\fbeta}\}\Rightarrow \Q\{\Z[x]{\falpha}\}$$
by setting $\phi(\T^{\falpha}_i)=\prod_{k=1}^s
\T^{m_{ki}{\fbeta}_k}$ and $\theta(\T^{\bbeta}_j)=\prod_{k=1}^n
\T^{n_{ki}{\falpha}_k}$.
It is clear that $\phi$ and $\theta$ are $\sigma$-morphisms and
group homomorphisms. Since $A=ANM$ and $B=BMN$, it can be checked
that $\phi\circ\theta = id$ and $\theta\circ\phi = id$. Then,
$\Q\{\Z[x]{\falpha}\}$ is isomorphic to $\Q\{\Z[x]{\fbeta}\}$ and
hence $T_{\falpha}^*$ and $T_{{\fbeta}}^*$ are isomorphic as
$\sigma$-varieties. \qedd

Since $T_{\falpha}^*$ is Cohn open in $\XB_{\falpha}$, as a
consequence of Corollary \ref{cor-tv25}, we have
\begin{cor}\label{cor-tv26}
Let $\falpha=\{\balpha_{1}, \ldots, \balpha_{n}\}$ and
${\fbeta}=\{{\bbeta}_{1}, \ldots, {\bbeta}_{s}\}$ be two sets of
generators for the same $\Zx$-lattice. Then $\XB_{\falpha}$ and
$\XB_{\fbeta}$ are birationally equivalent  $\sigma$-varieties.
\end{cor}

\begin{example}\label{ex-tv21}
Let $N_1=\left[\begin{array}{lll}
 x       & 1      & x   \\
 x-1     & 1      & x +1 \\
\end{array}\right]$,
 $N_2=\left[\begin{array}{lll}
  0      & 1      &  0  \\
 -1      & 1      &  1 \\
\end{array}\right]$, and
 $N_3=\left[\begin{array}{lll}
  0      & 1      &  0  \\
  0      & 0      &  1 \\
\end{array}\right]$.
Then $(N_1)=(N_2)=(N_3)$ and the corresponding toric $\sigma$-ideals
are $\I_1=[y_1y_3-y_2^{2x}]$, $\I_2=[y_1y_3-1]$, and $\I_3=[y_1-1]$,
respectively.
Then, $\V({\I_i}),i=1,2,3$ are isomorphic as tori in $(\AH^*)^3$ and
birationally equivalent as $\sigma$-varieties in $\AH^3$.
\end{example}

An algebraic torus is isomorphic to $(\CH^*)^m$ for some $m\in\N$
\cite{cox-2010}. The following example shows that this is not valid
in the  difference case.
\begin{example}\label{ex-ext1}
Let $\balpha_1 = (2)$, $\balpha_2 = (x)$, and $\balpha =
\{\balpha_1,\balpha_2\}$. Then $T^*_{\falpha} =
\V({\{y_1^x-y_2^2\}})$. We claim that $T^*_{\falpha}$ is not
isomorphic $\AH^*$. By Theorem \ref{th-tv41}, we need to show
$\E_1=\Q\{t,t^{-1}\}$ is not isomorphic to
$\E_2=\Q\{s^2,s^{-2},s^x,s^{-x}\}$, where $t$ and $s$ are
$\sigma$-indeterminates.
Suppose the contrary, there is an isomorphism $\theta:
\E_1\Rightarrow\E_2$ and $\theta(t)=p(s)\in \E_2$.
Then there exists a $q(z)\in\Q\{z\}$ such that $s^2=q(p(s))$ which
is possible only if $q=z,p=s^2$.
Since $s^x\in\E_2$, there exists an $r(z)\in\Q\{z\}$ such that
$s^x=r(s^2)$ which is impossible.
%
\end{example}


Let $T^{*} \subseteq (\AH^*)^{n}$ be a $\sigma$-torus and $ X
\subseteq \AH^{n}$  a $\sigma$-variety.  $T^{*}$ is said to {\em act
on $X$} if there exists a $\varphi \in \Q\{y_{1}, \ldots, y_{n},
z_{1}, \ldots, z_{n}\}$ such that for all $\E \in  \mathscr{E}_\F$,
$$
\varphi:T^{*}(\E)\times X(\E)\longrightarrow X(\E)
$$
is a $\sigma$-morphism as well as a group action of $T^{*}(\E)$ on
$X(\E)$, that is, $\varphi({\bm{1}}, \h) = \h$ and
$\varphi({\beps}_1\cdot {\beps}_2, \h) = \varphi({\beps}_1,
\varphi({\beps}_2, \h))$.

%
%
%
%
The following result gives a new characterization of toric
$\sigma$-varieties.
\begin{theorem}\label{th-toric}
A $\sigma$-variety $X$ is toric if and only if $X$ contains a
$\sigma$-torus $T^{*}$ as an open subset and with a group action of
$T^{*}$ on $X$ extending the natural group action of $T^{*}$ on
itself.
\end{theorem}
\proof $``\Rightarrow"$
For $\falpha$ given in \bref{eq-A}, let $X=\XB_{\falpha}$ be the
toric $\sigma$-variety defined by $\falpha$ and $T_{\falpha}^*$ the
$\sigma$-torus defined by $\falpha$.
By Corollary \ref{cor-tv42}, $T_{\falpha}^*$ is open in $X$.
The group action of $T_{\falpha}^{*}$ on itself can be extended to
$\AH^{n}$: $T_{\falpha}^{*}\times \AH^{n}\longrightarrow \AH^{n}$ by
$(a_1,...,a_n)\cdot (t_1,\ldots,t_n)=(a_1t_1,...,a_nt_n)$.
By Corollary \ref{cor-tv43}, for $t\in T_{\falpha}^{*}$,
$T_{\falpha}^{*}=t\cdot T_{\falpha}^{*} \subseteq t\cdot X$ and
$t\cdot X$ is an irreducible $\sigma$-variety whose defining
$\sigma$-ideal is $\I_t=\{f\in\Q\{\Y\}\,|\,f(t\Y)\in \IH(X)\}$.
Since $X$ is the Cohn closure of $T_{\falpha}\subset T_{\falpha}^*$
in $\AH^{n}$, we have $X\subseteq t\cdot X$. Set $t:=t^{-1}$, we
have $X=t\cdot X$. So $T_{\falpha}^{*}\times
T_{\falpha}^{*}\longrightarrow T_{\falpha}^{*}$ can be extended to
$T_{\falpha}^{*}\times X \longrightarrow X$.

$``\Leftarrow"$ Let $X$ be a $\sigma$-variety containing a
$\sigma$-torus $T^*$ which is isomorphic to $T^*_\balpha$ for
$\falpha$ defined in \bref{eq-A}. Then, we have the following
commutative diagram:
 \begin{equation}\label{eq-cg}
\xymatrix{
    T_{\falpha}^{*}\times T_{\falpha}^{*} \ar[r]^-{\phi} \ar[d]^-{(id,\omega)} & T_{\falpha}^{*} \ar[d]^-{\omega} \\
    T_{\falpha}^{*}\times X \ar[r]^-{\widetilde{\phi}}         & X}
 \end{equation}
where $\phi$ is the group operation of $T_{\falpha}^*$,
$\widetilde{\phi}$ is the extension of $\phi$ to $T^*\times X$, and
$\omega$ is the inclusion map of $T_{\falpha}^*$ into $\X$.

Since $T_{\falpha}^{*}$ is a $\sigma$-torus, by Theorem
\ref{th-tv41}, there exists a $\Z[x]$-lattice $L=(\balpha)$ in
$\Zx^m$ such that $T_{\falpha}^{*}\cong\Spec^{\sigma}(\Q\{L\})$. Let
$\Theta(X)$ be the set of all the $\sigma$-polynomial functions on
$X$ and $\Theta(X)\cong \Q\{X\}$.
From \bref{eq-cg}, we obtain the following commutative diagram
$$
\xymatrix{
   \Theta(X) \ar[r]^-{\widetilde{\Phi}} \ar[d]^-{\Omega} & \Q\{L\}\bigotimes_{\Q} \Theta(X)\ar[d]^-{(id,\Omega)}\\
    \Q\{L\}   \ar[r]^-{\Phi}     &  \Q\{L\}\bigotimes_{\Q}\Q\{L\}}
$$
Since $T_{\falpha}^{*}$ is an open set in $X$, ${\Omega}:
\Theta(X)\Rightarrow\Q\{L\}$ is the inclusion map of $\Theta(X)$
into $\Q\{L\}$. Therefore, each element $p \in \Theta(X)$ can be
written as $p=\sum_{\bbeta\in L}a_{\bbeta}\T^\bbeta$, where
$a_\bbeta$ is zero except a finite number of $\bbeta$.
From the morphism $\Phi$ in the diagram, we have $\Phi(p) =
\sum_{\bbeta\in L}a_{\bbeta} \T^\bbeta \bigotimes \T^\bbeta \in
\Q\{L\} \bigotimes_{\Q}\Q\{L\}$.
From the injective morphism $(id,\Omega)$ in the diagram,  we have
$\overline{p}=(id,\Omega)(\widetilde{\Phi}(p)) =\Phi(p) =
\sum_{\bbeta\in L}a_{\bbeta} \T^\bbeta \bigotimes \T^\bbeta \in
\Q\{L\}\bigotimes_{\Q} \Theta(X) $.
This is possible only if $\T^\bbeta\in \Theta(X)$ for each
$a_\bbeta\T^\bbeta\in L$. For otherwise, $\overline{p}$ should
contain a term of the form $\T^\bbeta\bigotimes (a_1
\T^{\bbeta_1}+a_2\T^{\bbeta_2})$ where $a_i\in \Q$,
$\bbeta_1\ne\bbeta_2$, and $\bbeta_1\not\in\Theta(X)$ and
$\bbeta_2\not\in\Theta(X)$. Thus, $\overline{p}$ contains the terms
$a_1\T^\bbeta\bigotimes \T^{\bbeta_1}+a_2 \T^\bbeta\bigotimes
\T^{\bbeta_2}$, which implies $\bbeta=\bbeta_1=\bbeta_2$, a
contradiction.
Therefore, there exists an affine $\N[x]$-module $S\subseteq L $
such that $\Theta(X)=\bigoplus_{\bbeta\in S} b_\bbeta \T^\bbeta$.
Since $\Zxn$ is Notherian,
%
there exists a finite set $\bgamma \subseteq L$ such that
$S=\N[x]\bgamma$ is an affine $\N[x]$-module and $\Theta(X) \cong
\Q\{S\}$. Thus $X=\Spec^{\sigma}(\Q\{S\})$.\qedd

Following Theorems \ref{th-tv31}, \ref{th-tvi},  and \ref{th-toric},
we have Theorem \ref{th-m3}.

\subsection{$\sigma$-Chow form and order of toric $\sigma$-variety}
In this section, we show that the $\sigma$-Chow form
\cite{dd-chowform,gao} of a toric $\sigma$-variety $\XB_{\falpha}$
is the sparse $\sigma$-resultant \cite{dd-sres} with support
${\falpha}$.
As a consequence, we can give a bound for the order of
$\XB_{\falpha}$.

Let ${\falpha}=\{\balpha_1,\ldots, \balpha_n\}$ be a subset of
$\Z[x]^m$ and $\XB_{\falpha}$ the toric $\sigma$-variety defined by
${\falpha}$.
In order to establish a connection between the $\sigma$-Chow form of
$\XB_{\falpha}$ and the $\sigma$-sparse resultant with support
${\falpha}$, we assume that ${\falpha}$ is {\em Laurent
transformally essential} \cite{dd-sres}, that is
 $$\rk(M_L)=m$$
where $M_L$ is the matrix with $\balpha_i$ as the $i$-th column.
Let $\T =(t_1, \ldots, t_m)$ be a set of $\sigma$-indeterminates.
Here, ${\falpha}$ is Laurent transformally essential means that
there exist indices $k_1,\ldots,k_m\in\{1,\ldots,n\}$ such that
$\T^{\balpha_{k_1}},\ldots,\T^{\balpha_{k_m}}$ are transformally
independent over $\Q$ \cite{dd-sres}.

Let
$\mathcal{A}=\{M_1=\T^{\balpha_{1}},\ldots,M_n=\Y^{\balpha_{n}}\}$
and
 \begin{equation}\label{eq-P}\P_i=u_{i0}+u_{i1}M_1+\cdots+u_{in}M_n\,(i=0,\ldots,m)\end{equation}
$m+1$ generic Laurent $\sigma$-polynomials with the same support
$\falpha$. Denote ${\bu_i}=(u_{i0},\ldots,u_{in}),$ $i=0,\ldots,m$.
Since $\A$ is Laurent transformally essential, the $\sigma$-sparse
resultant of $\P_0, \P_1,\ldots,$ $\P_m$ exists \cite{dd-sres},
which is denoted by $R_{\falpha}\in\Q\{\bu_0,\ldots,\bu_m\}$.

By Lemma \ref{lm-tv10},   $\XB_{\falpha}\subset \AH^n$ is an
irreducible $\sigma$-variety of dimension $\rk(M_L)=m$. Then, the
$\sigma$-Chow form of $\XB_{\falpha}$, denoted by
$C_{\falpha}\in\Q\{\bu_0,\ldots,\bu_m\}$, can be obtained by
intersecting $\XB_{\falpha}$ with the following generic
$\sigma$-hyperplanes \cite{dd-chowform}
 $$\L_i=u_{i0}+u_{i1}y_1+\cdots+u_{in}y_n\,(i=0,\ldots,m).$$
We have
\begin{theorem}\label{th-toric-chow}
Up to a sign, the sparse $\sigma$-resultant $R_{\falpha}$ of
$\P_i\,(i=0,\ldots,m)$ is the same as the $\sigma$-Chow form
$C_{\falpha}$ of $\XB_{\falpha}$ w.r.t. the generic hyperplanes
$\L_i\,(i=0,\ldots,m)$.
\end{theorem}
\proof Let $[\P_0,\ldots,\P_m]$ be the $\sigma$-ideal  in
$\Q\{\bu_0,\ldots,\bu_m,\T^{\pm}\}$. From \cite{dd-sres},
 $$[\P_0,\P_1,\ldots,\P_m]\cap\Q\{{\bu_0},\ldots,{\bu_m}\}=\sat(R_{\falpha},R_1,\ldots,R_l)$$
is a reflexive prime $\sigma$-ideal of codimension one in
$\Q\{\bu_0,\ldots,\bu_m\}$. Let $\I_{\falpha} = \IH(\XB_\balpha)$
and $[\I_{\falpha},\L_0,\L_1,\ldots,\L_m]$ be generated in
$\Q\{\Y,\bu_0,\ldots,\bu_m\}$. From \cite{dd-chowform},
 $$[\I_{\falpha},\L_0,\L_1,\ldots,\L_m]\cap
\Q\{{\bu_0},\ldots,{\bu_m}\}=\sat(C_{\falpha},C_1,\ldots,C_t)$$ is a
reflexive prime $\sigma$-ideal of codimension one in
$\Q\{\bu_0,\ldots,\bu_m\}$.
Then, it suffices to show
 $$[\P_0,\P_1,\ldots,\P_m]\cap
 \Q\{{\bu_0},\ldots,{\bu_m}\}=[\I_{\falpha},\L_0,,\L_1,\ldots,\L_m]\cap
 \Q\{{\bu_0},\ldots,{\bu_m}\}.$$
Let $\I_\T=[y_1-M_1,\ldots,y_n-M_n]$ be the $\sigma$-ideal generated
by $y_i-M_i,i=1,\ldots,n$ in $\Q\{\Y,\T^{\pm}\}$. By \bref{eq-tv01},
$\I_{\falpha}=\I_{\T}\cap\Q\{\Y\}$.
Then, $[\I_{\falpha},\L_0,\ldots,\L_m]\cap
\Q\{{\bu_0},\ldots,{\bu_m}\}=[y_1-M_1,\ldots,y_n-M_n,\L_0,\ldots,\L_m]\cap
\Q\{{\bu_0},\ldots,{\bu_m}\}=[y_1-M_1,\ldots,y_n-M_n,\P_0,\ldots,\P_m]\cap
\Q\{{\bu_0},\ldots,{\bu_m}\}$.
Since $\P_i\in\Q\{\bu_0,\ldots,\bu_n,\T^{\pm}\}$ does not contain
any $y_i^{x^j}$ in $\Theta(\Y)$, we have
$[y_1-M_1,\ldots,y_n-M_n,\P_0,\ldots,\P_m]\cap
\Q\{{\bu_0},\ldots,{\bu_m}\} = [\P_0,\ldots,\P_m]\cap
\Q\{{\bu_0},\ldots,{\bu_m}\}$, and the theorem is proved.\qedd

To give a bound for the order of $\XB_{\falpha}$, we need to
introduce the concept of Jacobi number.
Let $M=({m}_{ij})$ be an $m\times m$ matrix with elements either in
$\N$ or $-\infty$. A {\em diagonal sum} of $M$ is any sum
${m}_{1\sigma(1)}+{m}_{2\sigma(2)}+\cdots+{m}_{m\sigma(m)}$ with
$\sigma$ a permutation of $1,\ldots,m$. The {\em Jacobi number} of
$M$ is the maximal diagonal sum of $M$, denoted by $\Jac(M)$
\cite{dd-sres}.

Let $\falpha=\{{\falpha}_1,\ldots, {\falpha}_n\}\subset\Zx^{m}$ and
$A=(a_{ij})_{m\times n}$ the matrix with $\balpha_i$ as the $i$-th
column.
For each $i\in\{1,\ldots,m\}$, let $o_{i} = \max_{k=1}^n
\deg(a_{ik},x)$ and assume that $\deg(0,x)=-\infty$. Since $A$ does
not contain zero rows, no $a_{ij}$ is $-\infty$.
For a $p(x)\in\Zx$, let $\underline{\deg}(p,x) = \min\{k\in\N\,|\,
s.t.\, \coeff(p,x^k)\ne0 \}$ and $\underline{\deg}(0,x)=0$.
For each $i\in\{1,\ldots,m\}$, let $\underline{o}_i=\min_{k=1}^{n}
\underline{\deg}(a_{ik},x)$ and $\underline{o}=\sum_{i=1}^m
\underline{o}_i$. Then we have
\begin{theorem}\label{th-tv-ord}
Use the notations just introduced. Let $\XB_{\falpha}$ be the toric
variety defined by $\falpha$. Then $\ord(\XB_{\falpha}) \le
\sum_{i=1}^m (o_i-\underline{o}_i)$.
\end{theorem}
\proof Use the notations in Theorem \ref{th-toric-chow}. Since
$\P_i$ in \bref{eq-P} have the same support for all $i$,
$\ord(R_{\falpha},\bu_i)$ are the same for all $i$.
The {\em order matrix} for $\P_i$ given in \bref{eq-P} is
$O=(\ord(\P_i,t_j))_{(m+1)\times m}=(o_{ij})_{(m+1)\times m}$, where
$o_{ij} = o_j$. That is, all rows of $O$ are the same.
Let $\overline{O}$ be obtained from $O$ by deleting the any row of
$O$. Then $J=\Jac(\overline{O})=\sum_{i=1}^m o_i$.
By Theorem 4.17 of \cite{dd-sres}, $\ord(R_{\falpha},\bu_i) \le J-
\underline{o} = \sum_{i=1}^m (o_i-\underline{o}_i)$.
By Theorem 6.12 of \cite{dd-chowform}, $\ord(\XB_{\falpha}) =
\ord(C_{\falpha},\bu_i)$ for each $i=0,\ldots,m$.
By Theorem \ref{th-toric-chow}, $C_{\falpha} = R_{\falpha}$. Then
the theorem is proved.\qedd

\section{Algorithms for $\Zx$-lattice and binomial $\sigma$-ideal}\label{sec-alg}
In this section, we give algorithms to decide whether a given $\Zx$-lattice $L$ is $x$-, $\Z$-, $M$-, or $\Zx$-saturated, and in the
negative case to compute the $x$-, $\Z$-, $M$-, or $\Zx$-saturation
of $L$.
Based on these algorithms, we  give algorithms to compute the
reflexive closure, the well-mixed closure, the perfect closure of a
finitely generated Laurent binomial $\sigma$-ideal and an algorithm to
decompose a finitely generated perfect $\sigma$-ideal as the
intersection of reflexive prime $\sigma$-ideals.

In this section, the base files $\F$ is assumed to be inversive and
algebraically closed, since such conditions are required in Theorems \ref{th-pr1}, \ref{th-l2}, \ref{th-wm}, \ref{th-nl5}.

The following well-known algorithms will be used.
\begin{itemize}
\item Let $\fb$ be a finite set of elements in $\Zxn$.
We need to compute the reduced Gr\"obner basis of the $\Zx$-module
$(\fb)$ \cite[p. 197]{cox-1998}.

\item
Let $\D$ be $\Z$, $\Q[x]$, $\Z_p[x]$, or $\Q[x]/(q(x))$, where
$q(x)$ is an irreducible polynomial in $\Q[x]$. Then $\D$ is either
a PID or a field. For a finite set $S\subset \D^n$ and a matrix
$M\in\D^{n\times s}$, we need to compute the Hermite normal form of
the $\D$-module generated by $S$ and a basis for the $\D$-module:
$\ker(M) = \{X\in\D^s\,|\, MX=0\}$ \cite[p.68, p.74]{cohen}.
\end{itemize}

\subsection{$x$-saturation of $\Zx$-lattice}
\label{sec-alg1}

In this section, we give algorithms to check whether a $\Zx$-lattice
$L$ is $x$-saturated and in the negative case to compute the
$x$-saturation of $L$.

Let $\f_1,\ldots,\f_s\in\Zxn$ and $L=(\f_1,\ldots,\f_s)$. If $L$ is
not $x$-saturated, then there exist $g_i\in\Zx$ such that
 $\sum_{i=1}^s g_i \f_i = x \h$
and $\h\not\in L$. Setting $g_i(x) = g_i(0) + x\widetilde{g}_i(x)$
and $\widetilde{\h} = \h - \sum_{i=1}^s \widetilde{g}_i(x)\f_i$, we
have
  \begin{equation}\label{eq-x11}
 \sum_{i=1}^s g_i(0) \f_i = x \widetilde{\h}
 \end{equation}
where $\widetilde{\h}\not\in L$. Set $x=0$ in the above equation, we
have
  $$\sum_{i=1}^s g_i(0) \f_i(0) = 0,$$
that is, $G=(g_1(0),\ldots,g_s(0))^\tau$ is in the kernel of the
matrix $F=[\f_1(0),\ldots,\f_s(0)]\in\Z^{n\times s}$, which can be
obtained by existing algorithms \cite[page 74]{cohen}. From $G$ and
\bref{eq-x11}, we may compute $\widetilde{\h}$. This observation
leads to the following algorithm.

\begin{algorithm}[H]\label{alg-satX}
  \caption{\bf --- XFACTOR$(\f_1,\ldots,\f_s)$} \smallskip
  \Inp{A generalized Hermite normal form $\{\f_1,\ldots,\f_s\}\subset\Zxn$.}\\
  \Outp{$\emptyset$, if $L = (\f_1,\ldots,\f_s)$ is $x$-saturated;
    otherwise, a finite set
    $\{(\h_i,\e_i)\,|\, i=1,\ldots,r\}$ such that
    $\e_i=(e_{i1},\ldots,e_{is})^\tau\in\Z^s$, $\h_i\not\in L$, and
     $x \h_i = \sum_{l=1}^s e_{il} \f_l \in L$,
    $i=1,\ldots,r$.}\medskip

  \noindent
  1. Set $F=[\f_1(0),\ldots,\f_s(0)]\in\Z^{n\times s}$.\\
  2. Compute a basis $E\subset\Z^s$ of the $\Z$-module $\ker(F)$ with the algorithm in
  \cite[page 74]{cohen}.\\
  3. Set $H = \emptyset$.\\
  4. While $E \not= \emptyset$\\
   \SPC 4.1. Let $\e=(e_1,\ldots,e_s)^\tau\in E$ and $E = E\setminus \{\g\}$.\\
   \SPC 4.2. Let $\h=(e_1\f_1+\cdots+e_s\f_s)/x$.\\
   \SPC 4.3. If $\grem(\h,\{\f_1,\ldots,\f_s\})\ne0$,
   then add $(\h,\e)$ to $H$.\\
  5. Return $H$.
\smallskip
\end{algorithm}

We now give the algorithm to compute the $x$-saturation of a
$\Zx$-lattice.
\begin{algorithm}[H]\label{alg-satX1}
  \caption{\bf --- SATX$(\f_1,\ldots,\f_s)$} \smallskip
  \Inp{A finite set $\fb=\{\f_1,\ldots,\f_s\}\subset \Zxn$.}\\
  \Outp{A set of generators of $\sat_{x}(\f_1,\ldots,\f_s)$ .}\medskip

  \noindent
  1. Compute the generalized Hermite normal form $\gb$ of $\fb$. \\
  2. Set $H=$ {\bf XFACTOR}$(\gb)$.\\
  3. If $H=\emptyset$, then output $\gb$;
      otherwise set $\fb = \gb\cup \{\h\,|\, (\h,\f)\in H\}$ and goto step 1.  \\
\smallskip
\end{algorithm}

\begin{example}
Let $\C$ be the following generalized Hermite normal form.
\[ \C =[\f_1,\f_2,\f_3]
 = \left[  \begin{array}{llllll}
    -x+2               & 1      & 1\\
    3x+2               & 1      & 2x+1 \\
    0                  & 2x     & x^2 \\
\end{array} \right]. \]
In {\bf XFACTOR}$(\C)$, the kernel of the following matrix
\[
 [\f_1(0), \f_2(0), \f_3(0)] = \left[  \begin{array}{llllll}
  2               & 1      & 1\\
  2               & 1      & 1 \\
  0               & 0     & 0 \\
\end{array} \right] \]
is generated by $\e_1 = (0, -1,1)^\tau$ and $\e_2 = (1,-2,0)^\tau$.
In step 4.2 of {\bf XFACTOR}, we have  $ \C \e_1 = -\f_2+\f_3 = (0,
2x, x^2-2x)^\tau = x(0,2,x-2)^\tau$. One can check that
$(0,2,x-2)^\tau\not\in (\C)$.
In {\bf SATX}, computing the generalized Hermite normal form of
$\C\cup \{ (0,2,x-2)^\tau \} $, we have
\[ \C_1
 = \left[  \begin{array}{llllll}
    -x+2               & 1      & 0\\
    3x+2               & -3      & 2 \\
    0                  & 4     & x-2 \\
\end{array} \right]. \]
{\bf XFACTOR}$(\C_1)$ returns $\emptyset$. So, $(\C_1)$ is
 the $x$-saturation of $(C)$.
%
\end{example}

\begin{theorem}
Algorithms~{\bf SATX} and {\bf XFACTOR} are correct.
\end{theorem}
\proof From the output of Algorithm~{\bf XFACTOR}, in step 3 of {\bf
SATX}, we have $(\gb)\subsetneq(\gb\cup \{\h\,|\, (\h,\f)\in H\})\subseteq\sat_x(\fb)$.
Since $\Z[x]^n$ is a Noetherian $\Zx$-module, {\bf SATX} will
terminate and return the $x$-closure of $(\fb)$. So, it suffices to
show the correctness of Algorithm {\bf XFACTOR}.

We first explain step 4.2 of Algorithm {\bf XFACTOR}. Since
$\e\in\ker(F)$, $\h(0)=[\f_1(0),\ldots,$ $\f_s(0)]\e =
[0,\ldots,0]^\tau$. Therefore, $x$ is a factor of $e_1\f_1+\cdots
+e_s\f_s$ and thus $\h=(e_1\f_1+\cdots +e_s\f_s)/x\in\Zxn$.

To prove the correctness of Algorithm {\bf XFACTOR}, it suffices to
show that $L=\sat_x (L)$ if and only if for each $\e\in E$,
$e_1\f_1+\cdots +e_s\f_s = x\h$ implies $\h\in L$.

Let $E=\{\e_1,\ldots,\e_k\}$ where $\e_i\in\Z^s$. If $L=\sat_x (L)$,
then it is clear that $(\f_1,\ldots, \f_s)\e_i = x\h_i$ implies
$\h_i\in L$. To prove the other direction, let $[\f_1,\ldots,
\f_s]\e_i = x\h_i$ for $1\le i\le k$, where $\h_i\in L$.
Let $x\f\in L$. Then $x\f = \sum_{i=1}^s c_i(x)\f_i$, where
$c_i(x)\in\Z[x]$. If for each $i$, $x | c_i(x)$, then we have $\f =
\sum_{i=1}^s (c_i(x)/x) \f_i \in L$, and the lemma is proved.
Otherwise, set $x=0$ in $x\f = \sum_{i=1}^s c_i(x)\f_i$, we obtain
$\sum_{i=1}^s c_i(0)\f_i(0) = 0$. Hence
$Q=[c_1(0),\ldots,c_s(0)]^\tau \in \ker(F)$ and hence there exist
$a_i\in\Z, i=1,\ldots,k$ such that $Q = \sum_{i=1}^k a_i\e_i$. Then,
\[
 \begin{array}{llll}
 [\f_1,\ldots, \f_s]Q & =& \sum_{i=1}^k a_i [\f_1,\ldots, \f_s]\e_i
  =  \sum_{i=1}^k a_i x\h_i = x\widetilde{\h}, \\
\end{array}\]
where $\widetilde{\h}=\sum_{i=1}^k a_i \h_i\in L$. Then,
\[
 \begin{array}{llll}
 x\f & =& \sum_{i=1}^s c_i(x)\f_i \\
 &=& \sum_{i=1}^s c_i(0)\f_i + \sum_{i=1}^s x\overline{c}_i(x)\f_i \\
 &=& [\f_1,\ldots, \f_s]Q + x\sum_{i=1}^s \overline{c}_i(x)\f_i \\
 &=& x\widetilde{\h} + x\sum_{i=1}^s \overline{c}_i(x)\f_i,
\end{array}\]
where $\overline{c}_i(x) = (c_i(x)-c_i(0))/x \in \Z[x]$. Hence,  $\f
=\widetilde{\h} + \sum_{i=1}^s \overline{c}_i(x)\f_i \in L$ and the
lemma is proved. \qedd



\subsection{$\Z$-saturation of $\Zx$-lattice}
\label{sec-alg2}

The key idea to compute $\sat_\Z(L)$ for a $\Zx$-lattice
$L\in\Z[x]^n$ is as follows. Let $\fb=\{\f_1,\ldots,\f_s\}$.
Then $(\fb)$ is not $\Z$-saturated if and only if a linear
combination of $\f_i$ contains a nontrivial prime factor in $\Z$,
that is, $\sum_i g_i\f_i = p \f$, where $p$ is a prime number and
$\f\not\in (\fb)$.
%
%
%
Furthermore, $\sum_i g_i\f_i = p \f$ with $g_i\ne 0~\mod ~p$ is
valid if and only if $\f_1,\ldots,\f_s$ are linear dependent over
$\Z_p[x]$. The fact that $\Z_p[x]$ is a PID allows us to compute
such linear relations using methods of Hermite normal
forms~\cite{cohen}. The following algorithm is based on this
observation.

\begin{algorithm}[H]\label{alg-satZ1}
  \caption{\bf --- ZFACTOR} \smallskip
  \Inp{A generalized Hermite normal form $\C=\{\c_1,\ldots,\c_s\}\subset \Zxn$ given in  \bref{ghf}.}\\
  \Outp{$\emptyset$, if $L = (\C)$ is  $\Z$-saturated;
    otherwise, a finite set $\{(\h_i,k_i,\e_i)\,|\,i=1,\ldots,r\}$, such that
    $\h_i\in\Zxn$, $k_i\in\N$, $\e_i=(e_{i1},\ldots,e_{i,s})^\tau\in\Zx^s$,
    $\h_i\not\in L$ and $k_i \h_i = \sum_{l=1}^s e_{il} \c_l \in L$ for $i=1,\ldots,r$.}\medskip

  \noindent
  1. Read the numbers $t,r_i,l_i, c_{r_i,1,0},i=1,\ldots,t$ from \bref{ghf}.\\
  2. Set $ q = \prod_{i=1}^t c_{r_i,1,0}\in\N $.\\
  3. For any prime factor $p$ of $q$ do\\
    \SPC 3.1.    Set $F=[ \c_{r_1,l_1},\c_{r_2,l_2},\ldots,\c_{r_t,l_t}]\in \Z_p[x]^{n\times t}$.\\
    \SPC 3.2.  Compute a basis $G\subset\Z_p[x]^s$ of the $\Z_p[x]$-module $\ker(F)$
              with method  in \cite[page 74]{cohen}.\\
    \SPC 3.3.  If $G\ne\emptyset$, for each $\g=[g_1,\ldots,g_t]^\tau\in G$,
              $\sum_{i=1}^t g_i\c_{r_i,l_i} = p\h$ in $\Zxn$. \\
    \SPC\SPC  Return the set of $(\h,p,\e)$ where $\e\in\Zx^s$ and the place of $\e$ corresponding to $\c_{r_i,l_i}$
    \SPC\SPC  is $g_i$ and other places are zero.  \\
   \SPC 3.4.  Compute the Hermite normal form
     $\B=\{\b_1,\ldots,\b_t\}$ of $\{ \c_{r_1,l_1},\ldots,\c_{r_t,l_t} \}$ in $\Z_p[x]^n$.\\
   %
   \SPC 3.5.  Let $\C_{-} = \{\f_1,\ldots,\f_l\}$ be given in \bref{eq-inf} and
              $\widetilde{\f}_i =\grem(\f_i,\B)= \f_i +
                \sum_{k=1}^t a_{i,k} \b_k$,\\
    \SPC\SPC           in $\Z_p[x]^n$,   where
                $a_{i,k}\in\Z_p[x]$.\\
   \SPC 3.6.  If $\widetilde{\f}_i=0$ for some $i$, then
             $\f_i + \sum_{k=1}^t a_{i,k} \b_k = p_i\h_i$ in $\Zxn$.\\
    \SPC\SPC   Return the set of $(\h_i,p_i,\e_i)$ where $\e_i$ is a vector in $\Zx^s$
    such that\\
   \SPC\SPC  $(\c_1,\ldots,\c_s)\e_i=\f_i + \sum_{k=1}^t a_{i,k} \b_k= p_i\h_i$.  \\
   \SPC 3.7. Set $E=[\widetilde{\f}_1,\ldots,\widetilde{\f}_l]\in\Z_p[x]^{n\times l}$.\\
   \SPC 3.8. Compute a basis $D$ of $\{X\in\Z_p^l\,|\, EX=0\}$ as a vector space over $\Z_p$. \\
   \SPC 3.9.  If $D\ne\emptyset$, for each $\b=[b_1,\ldots,b_l]^\tau\in D$,
              $\sum_{i=1}^l b_i\widetilde{\f}_{_i} = p\h$ in $\Zxn$.\\
   \SPC\SPC            Return the set of $(\h,p,\e)$ where $\e_i$ is a vector in $\Zx^s$
    such that\\
    \SPC\SPC $(\c_1,\ldots,\c_s)\e=\sum_{i=1}^l b_i\widetilde{\f}_{_i} = p\h$.  \\
  4. Return $\emptyset$.
\smallskip
\end{algorithm}

\begin{remark}
In steps 3.6 and 3.9, we need to compute $\e_i$.
Since $\B=\{\b_1,\ldots,\b_t\}$ is the Hermite normal form of $\c=\{
\c_{r_1,l_1},\ldots,\c_{r_t,l_t} \}$ in $\Z_p[x]^n$, there exists an
invertible matrix $M_{t\times t}$ such that $[\b_1,\ldots,\b_t] =
[\c_{r_1,l_1},\ldots,\c_{r_t,l_t}] M$.
In Step 3.6, $\e_i$ can be obtained from the relation $\f_i +
\sum_{k=1}^t a_{i,k} \b_k = p_i\h_i$ and the relation
$[\b_1,\ldots,\b_t] = [\c_{r_1,l_1},\ldots,\c_{r_t,l_t}] M$. Step
3.9 can be treated similarly.
\end{remark}

\begin{remark}
In step 3.8, we need to compute a basis for the vector space
$\{X\in\Z_p^l\,|\, EX=0\}$ over $\Z_p$. We will show how to do this.
A matrix $F\in\Z_p[x]^{m\times s}$ is said to be in standard form if
$F$ has the structure in \bref{ghf} and $\deg(c_{r_i,k_1},x) <
\deg(c_{r_i,k_2},x)$ for $i=1,\ldots,t$ and $k_1 < k_2$.

The matrix $E\in\Z_p[x]^{n\times l}$ can be transformed into
standard form using the following operations: (1) exchange two
columns and (2) add the multiplication of a column by an element
from $\Z_p$ to another column. Equivalently, there exists an
inversive matrix $U\in\Z_p^{l\times l}$ such that $E\cdot U = S$ is
in standard form. Suppose that the first $k$ columns of $S$ are zero
vectors. Then the first $k$ columns of $U$ constitute a basis for
$\ker(E)$.
This can be proved similarly to that of the algorithm to compute a
basis for the kernel of a matrix over a PID \cite[page 74]{cohen}.
\end{remark}


We now give the algorithm to compute the $\Z$-saturation.
\begin{algorithm}[H]\label{alg-satZ}
  \caption{\bf --- SATZ($\f_0,\ldots,\f_s$)} \smallskip
  \Inp{A set of vectors $\fb=\{\f_0,\ldots,\f_s\}\subset\Zxn$.}\\
  \Outp{A generalized Hermite normal form $\gb$ such that $(\gb)=\sat_{\Z} (\fb)$.}\medskip

  \noindent
  1. Compute a generalized Hermite normal form $\gb$ of $\fb$.\\
  2. Set $S=${\bf ZFACTOR}$(\gb)$.\\
  3. If $S=\emptyset$, return $\gb$;
   otherwise set $\fb = \gb\cup \{\h\,|\, (\h,k,\f)\in H\}$ and goto step 1.
\smallskip
\end{algorithm}

\begin{example}
Let $\C$ be the following generalized Hermite normal form:
 \[ \C = \left[
\begin{array}{llllll}
    x^2+2x-2               & x+2      & 1 \\
    0                      &  4       & 2x \\
\end{array} \right]. \]
Then, $t=2, r_1 = 1, l_1 = 1, r_2 = 2, l_2 = 2, q = 4$, $\c_{1,1} =
[x^2+2x-2,0]^\tau$, $\c_{2,1} = [x+2,4]^\tau$, $\c_{2,2} =
[1,2x]^\tau$.

Apply algorithm {\bf ZFACTOR} to $\C$. We have $p=2$. In steps 3.1
and 3.2, $F =\left[\begin{array}{llllll}
    x^2       & 1\\
    0         & 0\\
\end{array}\right]$ and $\ker(F)$ is generated by
$G=\{[-1,x^2]^\tau\}$.
In step 3.3, $x^2\c_{22}-\c_{11}= 2(1-x,x^3)^\tau$ and return
$(1-x,x^3)^\tau$.

In Algorithm {\bf SATZ},  $(1-x,x^3)^\tau$ is added into $\C$ and
\[ \C_1 = \left[  \begin{array}{llllll}
    x^2+2x-2               & x+2      & 1  & 1-x\\
    0                      &  4       & 2x & x^3 \\
\end{array} \right], \]
which is also a generalized Hermite normal form.

Applying Algorithm {\bf ZFACTOR} to $\C_1$. We have $p=2$ and $t=2$.
In steps 3.1-3.3, $G=\emptyset$. In step 3.4,
 $\B=\left[\begin{array}{llllll}
    x^2       & 1-x\\
    0         & x^3\\
\end{array}\right]$.
In step 3.5, $\C_{-} =\left[  \begin{array}{lll}
     x +2  & 1 & x \\
        4    & 2x & 2x^2  \\
\end{array} \right]$ and $\widetilde{\f}_i\ne0$ for all $i$.
In step 3.7, $E=
 \left[  \begin{array}{llllll}
       x   & 1 & x \\
       0   & 0 & 0 \\
\end{array} \right]. $
In Step 3.8, $D=\{\b\}$, where $\b=[1,0,-1]^\tau$.
In Step 3.9,  $(x+2,4)^\tau - (x,2x^2)^\tau = 2(x+1,x^2+2)^\tau$.
Add $(x+1,x^2+2)^\tau$ is added into $\C_1$ and
compute the generalized Hermite normal form, we have
\[ \C_2 = \left[  \begin{array}{llllll}
    x^2+2x-2               & x+2      & 1  & x+1   \\
    0                      &  4       & 2x & x^2+2 \\
\end{array} \right]. \]
Apply Algorithm {\bf ZFACTOR} again, it is shown that $\C_2$ is
$\Z$-saturated.
\end{example}

We will prove the correctness of the algorithm. We denote by
$\sat_p(L)$ the set $\{\f\in\Zxn\, |\, p\f\in L \}$ where $p\in\Z$ a
prime number.
An infinite set $S$ is said to be {\em linear independent} over a
ring $R$ if any finite set of $S$ is linear independent over $R$,
that is $\sum_{i=1}^k a_i \g_i=0$ for $a_i\in R$ and $\g_i \in S$
implies $a_i=0, i=1,\ldots, k$.

\begin{lemma}\label{lm-satz1}
Let  $\C=\{ \c_1,\ldots, \c_s\}$ be a generalized Hermite normal
form and $L=(\C)$. Then $\sat_p(L) = L$ if and only if $\C_{\infty}$
is linear independent over $\Z_p$, where $\C_{\infty}$
is defined in \bref{eq-inf}.
\end{lemma}
\proof  $``\Rightarrow"$ Assume the contrary, that is, $\C_{\infty}
= \{\h_1,\h_2,\ldots\}$ defined in \bref{eq-inf} is linear dependent
over $\Z_p$. Then there exist $a_i\in\Z_p$ not all zero, such that
$\sum_{i=1}^r a_i\h_i = 0$ in $\Z_p[x]^n$ and hence $\sum_{i=1}^r
a_i\h_i = p\g$ in $\Z[x]^n$.
By Lemma \ref{lm-ext1}, $\h_i$ are linear independent over $\Z_p$
and hence $\g\ne{\bf{0}}$.
Since $\sat_p(L) = L$, we have $\g\in L$.
By Lemma \ref{lm-ext2}, there exist $b_i\in\Z$ such that $\g =
\sum_{i=1}^r b_i\h_i$. Hence $\sum_{i=1}^r(a_i-pb_i)\h_i = 0$ in
$\Zxn$. By Lemma \ref{lm-ext1}, $a_i = pb_i$ and hence $a_i=0$ in
$\Z_p[x]$, a contradiction.

$``\Leftarrow"$ Assume the contrary, that is, there exists a
$\g\in\Z[x]^n$, such that $\g\not\in L$ and $p\g \in L$. By Lemma
\ref{lm-ext2}, $p\g = \sum_{i=1}^r a_i\h_i$, where $a_i\in\Z$. $p$
cannot be  a factor of all $a_i$. Otherwise,  $\g = \sum_{i=1}^r
\frac{a_i}{p}\h_i\in L$. Then some of $a_i$ is not zero in $\Z_p$,
which means $\sum_{i=1}^r a_i\h_i=0$ is nontrivial linear relation
among $\C_i$ over $\Z_p$, a contradiction. \qedd

From the $``\Rightarrow"$  part of the above proof, we have
\begin{cor}\label{cor-satz1}
Let  $\sum_{i=1}^r a_i \h_i=0$ be a nontrivial linear relation among
$\h_i$ in $\Z_p[x]^n$, where $a_i\in \Z_p$. Then, in $\Z[x]^n$,
$\sum_{i=1}^r a_i \h_i=p\h$ and $\h\not\in (\C)$.
\end{cor}

\begin{lemma}\label{lm-satz11}
Let  $\C=\{ \c_1,\ldots, \c_s\}$ be a generalized Hermite normal
form and $L=(\C)$. Then $\sat_p(L) = L$ if and only if $\C_{\infty}$
is linear independent over $\Z_p$ for the prime factors of $q$
defined in step 2 of Algorithm {\bf ZFACTOR}.
\end{lemma}
\proof By Definition \ref{def-ghf}, the leading monomial of
$x^k\c_{r_i,j}\in\C_{\infty}$ is of the form
$c_{r_i,j,0}x^{k+d_{r_i,j}}{\beps}_{r_i}$ and $c_{r_i,l_i,0} |\ldots
| c_{r_i,2,0} $ $| c_{r_i,1,0}$.
If  $p$ is coprime with $\prod_{i=1}^t c_{r_i,1,0}$, then
$c_{r_i,j,0}\ne0\,\, \mod \,\, p$ for $1\le j\le l_i$. Therefore,
the leading monomials of the elements of $\C_{\infty}$ are linear
independent over $\Z_p$, and hence $\C_{\infty}$ is linear
independent over $\Z_p$. Therefore, it suffices to consider prime
factors of $\prod_{i=1}^t c_{r_i,1,0}$.\qedd

To check whether $\C_{\infty}$ is linear independent over $\Z_p$, we
first consider a subset of $\C_{\infty}$ is linear independent in the following lemma.
\begin{lemma}\label{lm-satz2}
Let  $\C$ be the generalized Hermite normal form given in
\bref{ghf}. Then $\C^{+}$ defined in \bref{eq-inf} is linear
independent over $\Z_p$ if and only if
$\{\c_{r_1,l_1},\c_{r_2,l_2},\ldots,\c_{r_t,l_t}\}$ are linear
independent over $\Z_p[x]$.
\end{lemma}
\proof This is obvious since $\sum_{i} \sum_j a_{i,j} x^j
\c_{r_i,l_i} =\sum_{i} p_i \c_{r_i,l_i}$, where $a_{i,j}\in\Z$ and
$p_i = \sum_j a_{i,j} x^j$.\qedd

\begin{lemma}\label{lm-satz3}
Let  $\B$ be an Hermite normal form  in $\Z_p[x]^n$ and
$\gb=\{\g_1,\ldots,\g_r\}\subset\Z_p[x]^n$. Then $\gb\cup
\B_{\infty}$ is linear dependent over $\Z_p$ if and only if
\begin{itemize}
\item either $\widetilde{\g}_i=\grem(\g_i,\B)=0$ in $\Z_p[x]^n$ for some $i$, or

\item the residue set $\{\grem(\g_i,\B)\, |\, i=1,\ldots,r\}$ are linear
dependent over $\Z_p$.
\end{itemize}
\end{lemma}
\proof We may assume that $\grem(\g_i,\B)=0$ does not happen, since
it gives a nontrivial linear relation of $\gb\cup \B_{\infty}$.
By Lemma \ref{lm-ext2}, $\widetilde{\g}_i =\g_i$ $\mod \,
\B_{\infty}$.
$\gb\cup \B_{\infty}$ is linear dependent over $\Z_p$ if and only if
there exist $a_i\in\Z_p$ not all zero such that $\sum_i a_i\g_i=0$
\mod \, $\B_{\infty}$ over $\Z_p$, which is valid if and only if
$\sum_i a_i \widetilde{\g_i}=0$ $\mod \,$ $\B_{\infty}$. Since
$\widetilde{\g}_i$ are G-reduced with respect to $\B$, $\sum_i a_i
\widetilde{\g}_i=0$ $\mod \,$ $\B_{\infty}$ if and only if $\sum_i
a_i \widetilde{\g}_i=0$, that is $\widetilde{\g}_i$ are linear
dependent over $\Z_p$.\qedd

\begin{theorem}
Algorithm~{\bf SATZ$(\f_0,\ldots,\f_s)$} is correct.
\end{theorem}
\proof Since $\Zxn$ is Notherian, it suffices to show that Algorithm
{\bf ZFACTOR} is correct.
Let $C = \{\c_1,\ldots,\c_s\}$. By Lemma~\ref{lm-satz1}, to check
whether $\sat_\Z(\c_1,\ldots,\c_s)$ is $\Z$-saturated, we need only
to check for any $p$ prime, $\C_{\infty}$ is linear independent on
$\Z_p$.

By Lemma \ref{lm-satz11}, we need only to consider prime factors of
$\prod_{i=1}^t c_{r_i,1,0}$ in step 3 of the algorithm.

In steps 3.1 and 3.2, we check whether $\C^{+}$ in \bref{eq-inf} is
linear independent over $\Z_p$. By Lemma \ref{lm-satz2}, we need
only to consider whether $\C_1 = \{\c_{r_1,l_1},
\c_{r_2,l_2},\ldots, \c_{r_t,l_t}\}$ is linear independent over
$\Z_p[x]$. It is clear that $\C_1$ is linear independent over
$\Z_p[x]$ if and only if $G=\emptyset$, where $\G$ is given in step
3.2.

In step 3.3, we handle the case that $\C_1$ is linear dependent over
$\Z_p$. If $G\ne\emptyset$ for any $\g=[g_1,\ldots,g_t]^\tau\in G$,
$\sum_{i=1}^t g_i\c_{r_i,l_i} =0$ in $\Z_p[x]$. Hence $\sum_{i=1}^t
g_i\c_{r_i,l_i}=p\h$ where $\h\in\Zxn$. By Corollary
\ref{cor-satz1}, $\h\not\in L$. The correctness of Algorithm {\bf
ZFACTOR} is proved in this case.

In steps 3.4 - 3.10, we handle the case where $\C^{+}$ is linear
independent over $\Z_p$.
In step 3.4, we compute the Hermite normal form of $\C_1$ in
$\Z_p[x]^n$, which is possible because $\Z_p[x]^n$ is a PID
\cite{cohen}. Furthermore, we have \cite{cohen}
 $$[\c_{r_1,l_1},\ldots,\c_{r_t,l_t}] N =[\b_1,\ldots,\b_t]$$
where $\{\b_1,\ldots,\b_t\}$ is an Hermite normal form and $N$ is an
inversive matrix in $\Z_p[x]^{t\times t}$.
Then $\C_{\infty}=\C_{-} \cup\C^{+}$ is linear independent over
$\Z_p$ if and only if
\begin{equation}\label{eq-satzp1}
 \widetilde{\C} = \C_{-} \cup \B_{\infty} = \C_{-} \cup
\cup_{j=0}^{\infty}\{x^j\b_1,\ldots,x^j\b_t\}
 \hbox{ is linear independent over }\Z_p.
\end{equation}
By Lemma \ref{lm-satz3}, property \bref{eq-satzp1} is valid if and
only if $\grem(\c,\B)\ne0$ for all $\c\in\C_{-}$ and the residue set
$\widetilde{\C}_{-}$ is linear independent over $\Z_p$, which are
considered in step 3.7 and steps 3.8-3.10, respectively.
Then we either prove $L$ is $\Z$-saturated or find a nontrivial
linear relation for elements in $\C_{\infty}$ over $\Z_p$. By
Corollary \ref{cor-satz1}, such a relation leads to an
$\h\in\sat_{\Z}(L)\setminus L$. The correctness of the algorithm is
proved.\qedd

As a direct consequence of Lemma \ref{lm-wmc} and Algorithm~{\bf ZFACTOR},
we have the algorithm to compute the M-saturation.
\begin{algorithm}[H]\label{alg-satM}
  \caption{\bf --- SATM($\f_0,\ldots,\f_s$)} \smallskip
  \Inp{A set of vectors $\fb=\{\f_0,\ldots,\f_s\}\subset\Zxn$.}\\
  \Outp{A generalized Hermite normal form $\gb$ such that $\sat_M(\fb)=(\gb)$.}\medskip

  \noindent
  1. Using Algorithm~{\bf ZFACTOR},
  we can compute $m_i\in\N$ and $\g_i\in\Zxn,i=1,\ldots,s$\\
  \SPC such that
     $\sat_\Z (\fb) = (\g_1,\ldots,\g_s)$  and $m_i\g_i\in (\fb)$.\\
  2. Let $S=\emptyset$ and for $i=1,\ldots,s$, if $m_i\ne1$ then $S=S\cup\{(x-o_{m_i})\g_i\}$.\\
  3. Compute the generalized Hermite normal form $\gb$ of $\fb\cup S$ and return $\gb$.\\
  \smallskip
\end{algorithm}

Notices that if $m_i=1$ then $o_{m_i}=0$ and $\g_i\in(\fb)$.
The numbers $m_i$ need not to be unique for the following reasons.
Suppose $m_i=n_ik$ and $n_i\g_i\in (\fb)$. Then by Corollary \ref{cor-per2},
$o_{m_i} = o_{n_i} + c n_i$ and hence $(x-o_{n_i})\g_i = (x-o_{m_i})\g_i + cn_i\g_i\in (\fb)$.
That is, we can replace $m_i$ by its factor $n_i$.

\begin{lemma}\label{lm-satm1}
Algorithm {\bf{SATM}} is correct.
\end{lemma}
\proof
Let $L_1=(\fb)$ and $L_2=(\fb, (x-o_{m_1})\g_1,\ldots,(x-o_{m_s})\g_s)$.
We claim that $\sat_\Z(L_1) = \sat_\Z(L_2)$.
Since $L_1\subset L_2$, $\sat_\Z(L_1) \subset \sat_\Z(L_2)$.
Since $\sat_\Z (L_1) = (\g_1,\ldots,\g_s)$, we have $L_2\subset\sat_\Z (L_1)$
and hence $\sat_\Z(L_2)\subset\sat_\Z (L_1)$. The claim is proved.
Then $\sat_\Z(L_2)=(\g_1,\ldots,\g_s)$ and $m_i\g_i\in L_1\subset L_2$.
Since $(x-o_{m_i})\g_i\in L_2,i=1,\ldots,s$, $L_2$ is M-saturated by Lemma \ref{lm-wmc}.\qedd

\subsection{$\Zx$-saturation of $\Zx$-lattice}
\label{sec-alg3}
The following algorithm checks whether a generalized Hermite normal
form $\C$ is $\Zx$-saturated or toric, and if not, it will return a
set of elements in $\sat_{\Zx}(\C)\setminus (\C)$.
\begin{algorithm}[H]\label{alg-D-satZX}
  \caption{\bf --- ZXFACTOR$(\c_1,\ldots,\c_s)$} \smallskip
  \Inp{A generalized Hermite normal form $\C=\{\c_1,\ldots,\c_s\}\subset \Zxn$ given in \bref{ghf}.}\\
  \Outp{$\emptyset$, if $L = (\C)$ is  $\Z[x]$-saturated;
      otherwise, a finite set $\{\h_1,\ldots,\h_r\}\subset\Zxn$
    such that $\h_i\not\in L$ and $\h_i\in \sat_{\Z[x]}( L)$, $i=1,\ldots,r$.}\medskip

  \noindent
  1. Let $S=${\bf{ZFACTOR}}$(\C)$. If $S\ne\emptyset$ return $S$. \\
  2. For any prime factor $p(x)$ of $\prod_{i=1}^t c_{r_i,1}\in\Zx\setminus\Z$.\\
   \SPC 2.1. Set $M=[\c_{r_1,1},\ldots,\c_{r_t,1}]\in \Z[x]^{n\times t}$,
   where $\c_{r_i,1}$ can be found in \bref{ghf}.\\
   \SPC 2.2. Compute a finite basis $B=\{\b_1,\ldots,\b_l\}$ of $\ker(M)=\{X\in\Q[x]^t \,|\,MX=0\}$\\
    \SPC\SPC as a vector space in $(\Q[x]/(p(x)))^t$.\\
   \SPC 2.3. If $B\ne\emptyset$,\\
   \SPC\SPC 2.3.1.   For each $\b_i$, let $M\b_i= p(x)
   \frac{\g_i}{m_i}$,
     where $\g_i\in\Zxn$ and $m_i\in\Z$.\\
   \SPC\SPC 2.3.2. Return $\{\g_1,\ldots,\g_l\}$ \\
  3. Return $\emptyset$.
\smallskip
\end{algorithm}

\begin{algorithm}[H]\label{alg-satZX}
  \caption{\bf --- SATZX$(\f_1,\ldots,\f_s)$} \smallskip
  \Inp{A finite set $\fb=\{\f_1,\ldots,\f_s\}\subset\Zxn$.}\\
  \Outp{A basis of $\sat_{\Z[x]}(\f_1,\ldots,\f_s)$ .}\medskip

  \noindent
  1. Compute a generalized Hermite normal form $\gb$ of $\fb$.\\
  2. Set $S=${\bf ZXFACTOR}$(\gb)$.\\
  3. If $S=\emptyset$, return $\gb$; otherwise set $\fb=\gb\cup S$ and go to step 1.
\smallskip
\end{algorithm}

\begin{example}
Let \[ \C
 = \left[  \begin{array}{llllll}
    x              & x^2+1     \\
    2x^2+1            & 0    \\
    0                 & 4x^2+2 \\
\end{array} \right]. \]
Apply Algorithm~{\bf{ZXFACTOR}} to $\C$. In step, 1, $S=\emptyset$
and $\C$ is $\Z$-saturated.
In step 2, the only irreducible factor of $\prod_{i=1}^t
c_{r_i,1}\in\Zx$ is $p(x)=2x^2+1$.
In step 2.1, $M=\C$ and in step 2.2, $B=\{[-1,2x]^\tau\}$. In step
2.3.1, $M\cdot [-1,2x]^\tau =  2x\c_{2,1}-\c_{1,1} =
p(x)[x,-1,4x]^\tau = 0~\mod ~p(x)$ and $\{[x,-1,4x]^\tau\}$ is
returned.

In Algorithm {\bf{SATZX}}, $\h=[x,-1,4x]^\tau$ is added into $ \C $
and the generalized Hermite normal form of $\C\cup\{\h\}$ is
 \[ \C_1
 = \left[  \begin{array}{llllll}
    x         & 1    \\
    2x^2+1    & x   \\
    0          & 2 \\
\end{array} \right]. \]
Apply Algorithm~{\bf{ZXFACTOR}} to $\C_1$,  one can check that
$\C_1$ is $\Z[x]$-saturated.
\end{example}

In the rest of this section, we will prove the correctness of the
algorithm.
Denote $L_{\Q}$ to be the $\Q[x]$-module generated by $L$ in
$\Q[x]^n$. Similar to the definition of $\sat_{\Z[x]}(L)$, we can
define $\sat_{\Q[x]}(L_{\Q})$. $L_{\Q}$ is called {\em $\Q[x]$-saturated}
if $\sat_{\Q[x]}(L_{\Q})=L_{\Q}$. The following lemma
gives a criterion for whether $L$ is $\Zx$-saturated.
\begin{lemma}\label{lm-satzx1}
A $\Zx$-lattice $L$ is $\Zx$-saturated if and only if $\sat_\Z(L) =
L$ and $\sat_{\Q[x]} (L_{\Q}) = L_{\Q}$.
\end{lemma}
\proof
%
$``\Rightarrow"$ If $L=(\f_1,\ldots, \f_s)$ is $\Zx$-saturated, then
$\sat_\Z(L) = L$. If $\sat_{\Q[x]} (L_{\Q}) \ne L_{\Q}$, then there
exists an $h(x)\in\Q[x]$ and a $\g\in\Q[x]^n$, such that $h(x)\g\in
L_\Q$ but $\g\not\in L_\Q$. From $h(x)\g\in L_\Q$, we have $h(x)\g =
\sum_{i=1}^s q_i(x)\f_i$ where $q_i(x)\in\Q[x]$. By clearing the
denominators of the above equation, there exist $m_1, m_2\in \Z$
such that $m_1h(x)\in\Zx$, $m_2\g\in\Zxn$, and $m_1h(x)\cdot
m_2\g\in L$. Since $L$ is $\Zx$-saturated, $m_2\g\in L$, which
contradicts to $\g\not\in L_\Q$.

$``\Leftarrow"$ For any $h(x)\in\Zx$ and $\g\in\Zxn$,  if $h(x)\g\in
L$, we have $h(x)\g\in L_\Q$, and hence $\g\in L_\Q$ since
$\sat_{\Q[x]} (L_{\Q}) = L_{\Q}$. From $\g\in L_\Q$, there exists an
$m\in \Z$ such that $m\g\in L$ which implies $\g\in L$ since $L$ is
$\Z$-saturated. \qedd

The following lemma shows that a generalized Hermite normal form
becomes an Hermite normal form in $\Q[x]^n$.
\begin{lemma}\label{lm-satzx2}
Let $\C$ be the generalized Hermite normal form given in \bref{ghf}.
Then $(\C)=(\c_{r_1,1},\ldots,$ $\c_{r_t,1})$ as $\Q[x]$-modules in
$\Q[x]^n$ and $[\c_{r_1,1},\ldots,\c_{r_t,1}]$ is an Hermite normal
form.
\end{lemma}
\proof It is clear that $[\c_{r_1,1},\ldots,\c_{r_t,1}]$ is an
Hermite normal form.
We will prove $(\C)=(\c_{r_1,1},\ldots,$ $\c_{r_t,1})$ by induction.
By 3) of Definition \ref{def-ghf}, $S(\c_{r_1,1},\c_{r_1,2}) = x^u
\c_{r_1,1} - a\c_{r_1,2}$ ($u\in\N$ and $a\in\Z$) can be reduced to
zero by $\c_{r_1,1}$, which means $\c_{r_1,2} = q(x)\c_{r_1,1}$
where $q(x)\in \Q[x]$. Hence, $(\c_{r_1,1},\c_{r_1,2})=
(\c_{r_1,1})$ as $\Q[x]$-modules.
Suppose for $k<l_1$, $(\c_{r_1,1},\ldots,\c_{r_1,k})= (\c_{r_1,1})$
as $\Q[x]$-modules. We will show that
$(\c_{r_1,1},\ldots,\c_{r_1,k+1})= (\c_{r_1,1})$ as $\Q[x]$-modules.
Indeed, by 3) of Definition \ref{def-ghf},
$S(\c_{r_1,1},\c_{r_1,k+1}) = x^v \c_{r_1,1} - b\c_{r_1,k+1}$
($v\in\N$ and $b\in\Z$) can be reduced to zero by
$\c_{r_1,1},\ldots,\c_{r_1,k}$ and hence, $\c_{r_1,k+1} \in
(\c_{r_1,1})$ as $\Q[x]$-modules. Then we have
$(\c_{r_1,1},\ldots,\c_{r_1,l_1})= (\c_{r_1,1})$ as $\Q[x]$-modules.
For the rest of the polynomials in $\C$, the proof is similar. \qedd

The following lemma gives a criterion for a $\Q[x]$-module to be
$\Q[x]$-saturated.
\begin{lemma}\label{lm-satzx3}
Let $\C$ be the generalized Hermite normal form given in \bref{ghf}
and $L=(\C)$. Then $\sat_{\Q[x]}(L_{\Q})=L_{\Q}$ if and only if
$\C_1=\{\c_{r_1,1},\ldots,\c_{r_t,1}\}$ is linear independent over
$\K_{p(x)}=\Q[x]/(p(x))$ for any irreducible polynomial
$p(x)\in\Z[x]$.
\end{lemma}
\proof  $``\Rightarrow"$ Assume the contrary, that is, $\C_1$ are
linear dependent over $\K_{p(x)}$ for some $p(x)$. Then there exist
$g_i\in\Q[x]$  not all zero in $\K_{p(x)}$, such that $\sum_{i=1}^t
g_i\c_{r_i,1} = 0$ in $\K_{p(x)}^n$ and hence $\sum_{i=1}^t g_i
\c_{r_i,1} = p(x)\g$ in $\Q[x]^n$.
Since $\C_1$ is in triangular form and is clearly linear independent
in $\Q[x]^n$, we have $\g\ne{\bf{0}}$.
Since $\sat_{\Q[x]}(L_{\Q})=L_{\Q}$, we have $\g\in L_{\Q}$.
Then, there exist $f_i\in\Q[x]$ such that $\g = \sum_{i=1}^t
f_i\c_{r_i,1}$. Hence $\sum_{i=1}^t(g_i-pf_i)\c_{r_i,1} = 0$ in
$\Q[x]^n$. Since $\C_1$ is linear independent in $\Q[x]^n$, $g_i =
pf_i$ and hence $g_i=0$ in $\K_{p(x)}$, a contradiction.

$``\Leftarrow"$ Assume the contrary, that is, there exists a
$\g\in\Q[x]^n$, such that $\g\not\in L_{\Q}$ and $p(x)\g \in L_{\Q}$
for an irreducible polynomial $p(x)\in\Z[x]$. Then, by Lemma
\ref{lm-satzx2} we have $p\g = \sum_{i=1}^t f_i\c_{r_i,1}$, where
$f_i\in\Q[x]$. $p$ cannot be a factor of all $f_i$. Otherwise, $\g =
\sum_{i=1}^t \frac{f_i}{p}\c_{r_i,1}\in L_{\Q}$. Then some of $f_i$
is not zero in $\K_{p(x)}$, which means $\sum_{i=1}^t
f_i\c_{r_i,1}=0$ is a nontrivial linear relation among $\C_1$ over
$\K_{p(x)}$, a contradiction. \qedd

From the $``\Rightarrow"$  part of the above proof, we have
\begin{cor}\label{cor-satzx1}
Let $\C$ be the generalized Hermite normal form given in \bref{ghf}
and $\sum_{i=1}^t f_i \c_{r_i,1}$ $=0$ a nontrivial linear relation
among $\c_{r_i,1}$ in $(\Q[x]/(p(x)))^n$, where $p(x)$ is an
irreducible polynomial in $\Zx$ and  $f_i\in \Q[x]$. Then, in
$\Q[x]^n$, $\sum_{i=1}^r f_i \c_{r_i,1}=p(x)\g$ and $\g\not\in (\C)$
as a $\Q[x]$-module.
\end{cor}

\begin{theorem}
Algorithm~{\bf SATZX$(\f_1,\ldots,\f_s)$} is correct.
\end{theorem}
\proof Since $\Zxn$ is Notherian, we need only to show the
correctness of Algorithm {\bf ZXFACTOR}.
In step 1, $L$ is $\Z$ saturated if and only if $S\ne\emptyset$.
In step 2, we claim that $L$ is $\Q[x]$-saturated if and only if
$B=\emptyset$ and if $B\ne\emptyset$ then $\g_i$ in step 2.3.1 is
not in $L$. In step 3, $L$ is both $\Z$ and $\Q[x]$ saturated. By
Lemma \ref{lm-satzx1}, $L$ is $\Z[x]$-saturated and the algorithm is
correct. So, it suffices to prove the claim about step 2.

By Lemma~\ref{lm-satzx3}, to check whether
$\sat_{\Q[x]}(L_{\Q})=L_{\Q}$, we  need only to check whether for
any irreducible polynomial $p(x)\in\Zx$,
$\C_{1}=\{\c_{r_1,1},\ldots,\c_{r_t,1}\}$  is linear independent
over $\K_{p(x)} = \Q[x]/(p(x))$.
If  $p(x)$ is not a prime factor of  $\prod_{i=1}^t c_{r_i,1}$, then
the leading monomials of $\c_{r_i,1}, i=1,\ldots,t$ are nonzero.
Since $\C_1$ is an Hermite normal form, $\C_1$ is linear independent
over $\K_{p(x)}$. Hence, we need only to consider prime factors of
$\prod_{i=1}^t c_{r_i,1}$ in step 2 of the algorithm.
In step 2.3, it is clear that if $\C_1=\emptyset$ then $\C_1$ is
linear independent over $\K_{p(x)}$.
For $\b_i\in B$, since $M\b_i=0$ over $\K_{p(x)}$, $M\b_i = p(x)
\h_i$ where $\h_i\in\Q[x]^t$. Hence $\h_i = \frac{\g_i}{m_i}$ for
$\g_i\in\Zx^t$ and $m_i\in\Z$.
By Corollary \ref{cor-satzx1}, $\g_i\not\in L$.
Therefore, step 2 returns a set of nontrivial factors of $L$ if $L$
is not $\Z[x]$-saturated. The claim about step 2 is proved. \qedd

\subsection{Algorithms for binomial $\sigma$-ideals}
\label{sec-alg4}
In this section, we will present several algorithms for Laurent
binomial and binomial $\sigma$-ideals, and in particular a
decomposition algorithm for binomial $\sigma$-ideals.
%
%
We first give an algorithm to compute the reflexive closure for a
Laurent binomial $\sigma$-ideal.
\begin{algorithm}[H]\label{alg-ref}
  \caption{\bf --- REFLEXIVE} \smallskip
  \Inp{$P$: a finite set of Laurent $\sigma$-binomials in $\F\{\Y^{\pm}\}$, where $\F$ is inversive.}\\
  \Outp{$\A$: a regular and coherent Laurent binomial $\sigma$-chain such that
    $[\A]$ is the reflexive closure of $[P]$.}
 \medskip
 \noindent
 1. Check whether $[P]=[1]$  with Lemma \ref{lm-bi2}. If $[P]=[1]$, return $1$.
 Otherwise, set $F=P$.\\
 2. Let $F=\{\Y^{\f_1}-c_1,\ldots,\Y^{\f_r}-c_r\}$ and $\fb=\{\f_1,\ldots,\f_r\}$.\\
 3. Compute the generalized Hermite normal form $\gb$ of $\fb$.\\
 4. Let $\gb=\{\g_1,\ldots\g_s\}$ and $\g_i = \sum_{k=1}^r a_{i,k}
 \f_k$, where $a_{i,k}\in \Zx$, $i=1,\ldots,s$.\\
 5. Let $G=\{g_1,\ldots,g_s\}$, where $g_i = \Y^{\g_i}-d_i$ and $d_i =
 \prod_{k=1}^r c_k^{a_{i,k}}$, $i=1,\ldots,s$.\\
 6. $H=${\bf{XFACTOR}}$(\gb)$.\\
 7. If $H=\emptyset$, return $G$.\\
 8. Let $H= \{(\h_i,\e_i)\,|\, i=1,\ldots,r\}$ and
 $\e_i=(e_{i1},\ldots,e_{is})$.\\
 9. Let $F:=G\cup\{\Y^{\h_i}- \sigma^{-1}(\prod_{j=1}^s d_j^{e_{i,j}}),i=1,\ldots,r\}$,
                   and go to step 2.\\
\end{algorithm}
\begin{theorem}\label{th-algref}
Algorithm  {\bf REFLEXIVE} is correct.
\end{theorem}
\proof The algorithm basically follows the proof of Theorem
\ref{lm-ld0}. In steps 2-5, $(\fb)$ and $(\gb)$ are the support
lattices of $[F]$ and $[G]$, respectively. By Lemma \ref{lm-bi3},
$[F]=[G]$. By Lemma \ref{lm-t4}, $G$ is a regular and coherent
$\sigma$-chain. In step 7, if $H=\emptyset$, then $(\gb)$ is
$x$-saturated, and by Theorem \ref{th-pr1}, $[G]$ is reflexive and
the theorem is proved.
Otherwise, we execute steps 8 and 9. Let $\I_1=[F]$, $L_1=(\fb)$,
$\I_2 = [G\cup\{\Y^{\h_i}- \sigma^{-k_i}(\prod_{j=1}^s
d_j^{b_{i,j}}),i=1,\ldots,t\}]$, and $L_2 =\L(\I_2)$.  Then, we have
$\I_1\varsubsetneq\I_2\subset \I_x$ and $L_1\varsubsetneq L_2\subset
L_x$, where $L_x=\sat_x(L_1)$ and $\I_x$ is the reflexive closure of
$\I_1$.
Similar to the proof of Theorem \ref{lm-ld0}, the algorithm will
terminate and output the reflexive closure of $[P]$.\qedd

\begin{remark}
Following Algorithm 5, we can give an algorithm to check whether a Laurent
binomial $\sigma$-ideal is well-mixed or perfect, and in the negative case to
compute the well-mixed or perfect closure of the $\sigma$-ideal. The details are
omitted.
\end{remark}

We give a decomposition algorithm for Laurent binomial
$\sigma$-ideals.
\begin{algorithm}[H]\label{alg-ldec}
  \caption{\bf --- DECLAURENT} \smallskip
  \Inp{$P$: a finite set of Laurent $\sigma$-binomials in $\F\{\Y^{\pm}\}$.}\\
  \Outp{$\emptyset$, if $\{P\}=[1]$ or
  regular and coherent Laurent binomial
  $\sigma$-chains $\C_1,\ldots,\C_t$ in $\F\{\Y^{\pm}\}$
  such that $[\C_i]$ are Laurent reflexive prime $\sigma$-ideals and
 $\{P\} = \cap_{i=1}^t [\C_i]$.}

 \medskip
 \noindent
 1. Let $F=${\bf{REFLEXIVE}}$(P)$. If $F=1$, return $\emptyset$. \\
 2. Set $\RB=\emptyset$ and $\FB=\{F\}$.\\
 3. While $\FB\ne\emptyset$.\\
 \SPC 3.1. Let $F=\{\Y^{\f_1}-c_1,\ldots,\Y^{\f_r}-c_r\}\in\FB$,
  $\fb=\{\f_1,\ldots,\f_r\}$,
  $\FB=\FB\setminus\{F\}$.\\
 \SPC 3.2. Check whether $[F]$ is proper with Lemma \ref{lm-bi2}. If
  $[F]=[1]$, go to step 3.\\
 \SPC 3.3. Compute the generalized Hermite normal form $\gb=\{\g_1,\ldots\g_s\}$ of
 $\fb$,\\
 \SPC \SPC where $\g_i = \sum_{k=1}^r a_{i,k} \f_k$ for $a_{i,k}\in \Zx$, $i=1,\ldots,s$.\\
 \SPC 3.4. Let $G=\{g_1,\ldots,g_s\}$, where $g_i = \Y^{\g_i}-d_i$ and $d_i =
 \prod_{k=1}^r c_k^{a_{i,k}}$, $i=1,\ldots,s$.\\
 \SPC 3.5. $H=${\bf{ZFACTOR}}$(\gb)$.\\
 \SPC 3.6. If $H=\emptyset$, add $G$ to $\RB$.\\
 \SPC 3.7. Let $H=\{(\h_i,k_i,\e_i)\,|\, i=1,\ldots,r\}$ and $\e_i=(e_{i1},\ldots,e_{is})$.\\
 \SPC 3.8. For $i=1,\ldots,t$,
  let $r_{i,1},\ldots,r_{i,k_i}$ be the $k_i$-th roots of  $\prod_{j=1}^s d_j^{e_{i,j}}$.\\
 \SPC 3.9. For $l_1=1,\ldots,k_1$,\ldots,$l_t=1,\ldots,k_t$,
            add $G\cup\{\Y^{\h_1}- r_{1,l_1},\ldots,\Y^{\h_t}- r_{t,l_t}\}$ to $\FB$.\\
 4.  Return $\RB$.
\end{algorithm}

\begin{theorem}
Algorithm {\bf DECLAURENT} is correct.
\end{theorem}
\proof The algorithm basically follows the proof of Theorem
\ref{th-l2}. The proof is similar to that of Theorem
\ref{th-algref}. \qedd

In the rest of this section, we give a decomposition algorithm for
binomial $\sigma$-ideals and hence a proof for Theorem \ref{th-m4}.
Before giving the main algorithm, we give a sub-algorithm {\bf
DECMONO} which treats the $\sigma$-monomials. Basically, it gives
the following decomposition
$$ \V(\prod_{i=1}^n y_i^{f_i}) =
 \V(y_1)\cup\V(y_2/y_1)\cup\cdots\cup\V(y_n/\{y_1,\ldots,y_{n-1}\})$$
 where $0\ne f_i\in\N[x]$ and $\V(y_c/S)$ is the set of zeros of $y_c=0$ not
 vanishing any of the variables in $S$. The correctness of the
 algorithm comes directly from the above formula.
\begin{algorithm}[H]\label{alg-decmono}
  \caption{\bf --- DECMONO} \smallskip
  \Inp{$(\Y_0,B,\Y_1)$: $\Y_0,\Y_1$ are disjoint subsets of $\Y$ and
   $B$ a finite set of $\sigma$-binomials in $\F\{\Y\}$.}\\
  \Outp{$(\Y_{0i},B_i,\Y_{1i})$: $\Y_{0i},\Y_{1i}$ are disjoint subsets of $\Y$,
   $B_i$ contains no $\sigma$-monomials,
   and  $\V(\Y_0\cup B/\Y_1) = \cup_{i=1}^r \V(\Y_{0i}\cup B_i/\Y_{1i})$.}

 \medskip
  \noindent
 1. Set $\RB=\emptyset$ and $\FB=\{(\Y_0,B,\Y_1)\}$.\\
 2. While $\FB\ne\emptyset$.\\
 \SPC 2.1. Let $C=(\Y_0,B,\Y_1)\in\FB$, $\FB=\FB\setminus\{C\}$.\\
 \SPC 2.2. For all $y_c\in\Y_0$, let $B_1=B_{y_c=0}$ (replace $y_c$ by $0$) and delete $0$ from $B_1$.\\
 \SPC 2.3. If $B_1$ contains no $\sigma$-monomials, add $(\Y_0,B_1,\Y_1)$ to $\RB$ and goto step 2.\\
 \SPC 2.4. Let $M=\prod_{i=1}^k y_{c_i}^{f_i}\in B_1$, where $0\ne f_i\in\N[x]$.
           $B_1=B_1\setminus\{M\}$.\\
 \SPC 2.5. Let $\Y_2:=\{y_{c_1},\ldots,y_{c_k}\}\setminus\Y_1$.
 If $\Y_2=\emptyset$, go to step 2; else let $\Y_2=\{y_{t_1},\ldots,y_{t_s}\}$.\\
 \SPC 2.6. For $i=1,\ldots,s$,
    add $(\Y_0\cup\{y_{t_i}\},B_1,\Y_1\cup\{y_{t_1},\ldots,y_{t_{i-1}}\})$ to $\FB$.\\
  3. Return $\RB$.
\smallskip
\end{algorithm}

%
%

We now give the main algorithm. The algorithm basically follows the
proof of Theorem \ref{th-nl5}. The main modification is that instead
of the perfect ideal decomposition
$$\{F\}= (\{F\}:\mb)\bigcap \cap_{i=1}^n \{F,y_i\},$$
we use the following zero decomposition
$$\V(F)= \V(\{F\}:\mb)\bigcup \cup_{i=1}^n
  \V(F\cup\{y_i\}/\{y_1,\ldots,y_{i-1}\}).$$
The purpose of using the later decomposition is that many redundant
components can be easily removed by the following criterion
$\V(F/D)=\emptyset$ if $F\cap D\ne\emptyset$, which is done in step
2.5 of Algorithm {\bf DECMONO}.
\begin{algorithm}[H]\label{alg-dec}
  \caption{\bf --- DECBINOMIAL} \smallskip
  \Inp{$F$: a finite set of $\sigma$-binomials in $\F\{\Y\}$.}\\
  \Outp{$\emptyset$, if $\{F\}=[1]$
  or $(\C_1,\Y_1),\ldots,(\C_r,\Y_r)$,
  where $\Y_i\subset\Y$ and
  $\C_i$ are regular and coherent $\sigma$-chains containing $\sigma$-binomials of variables in
  $\Y\setminus\Y_i$
  such that
  $\sat(\C_i)$ are reflexive prime $\sigma$-ideals
 and
 $\{F\} = \cap_{i=1}^r \sat(\C_i)$.}

 \medskip
  \noindent
 1. Set $\RB=\emptyset$ and $\FB=${\bf DECMONO}$(\emptyset,F,\emptyset)$.\\
 2. While $\FB\ne\emptyset$.\\
 \SPC 2.1. Let $C=(\Y_0,B,\Y_1)\in\FB$, $\FB=\FB\setminus\{C\}$.\\
 \SPC 2.2. If $B=\emptyset$, add $(\Y_0,\Y_1)$ to $\RB$.\\
 \SPC 2.3. Let $E=$ {\bf DECLAUENT}$(B)$ in $\F\{\Z^{\pm}\}$, where $\Z=\Y\setminus\Y_0$ and $m=|\Z|$. \\
 \SPC 2.4. If $E=\emptyset$ goto step2.\\
 \SPC 2.5. Let $E=\{E_1,\ldots,E_l\}$ and $E_{l}=$
    $\{\Z^{\f_{l,1}}-c_{l,1},\ldots, \Z^{\f_{l,s_l}}-c_{l,s_l}\}$,
    where $\f_{l,k}\in\Zx^{m}$.\\
 \SPC 2.7. Add
         $(\{\Y_0,     \Z^{\f_{l,1}^+}-c_{l,1}\Z^{\f_{l,1}^-},
              \ldots, \Z^{\f_{l,s_l}^+}-c_{l,s_l}\Z^{\f_{l,s_l}^-}\},\Y_1)$
              to $\RB$, $l=1,\ldots,k$.\\
 \SPC 2.8. Let $\Z = \{y_{c_1},\ldots,y_{c_s}\}$. For
 $i=1,\ldots,s$, do \\
 \SPC\SPC $\FB=\FB\cup$ {\bf DECMONO}$(\Y_0\cup\{y_{c_i}\},B,\Y_1\cup\{y_{c_1},\ldots,y_{c_{i-1}}\})$.\\
  3. Return $\RB$.
\smallskip
\end{algorithm}

\begin{theorem}\label{th-decbi}
Algorithm {\bf DECBINOMIAL} is correct.
\end{theorem}
\proof In step 1, $\V(F)$ is decomposed as $\V(F) = \cup_{i=1}^r
\V(\Y_{0i}\cup B_i/\Y_{1i})$ and $\FB=\{(\Y_{0i},B_i,$
$\Y_{1i});i=1,\ldots,r\}$.
In step 2, we will treat the components of $\FB$ one by one. In step
2.1, the component $(\Y_0,B,\Y_1)$ is taken from $\FB$.
In steps 2.3-2.4, $\{B\}$ is decomposed as
 $$\{B\}=\cap^{k}_{l=1}[E_l]$$
 in $\F\{\Z^{\pm}\}$, where $E_l$ are regular and coherent $\sigma$-chains and $[E_l]$
 are reflexive prime ideals.
By \bref{eq-nbi1} and Corollary \ref{cor-ascr1}, we have
 \begin{equation}\label{eq-ma01}
 \{B\}:\mb=\{B\}\cap\F\{\Z\}= \cap^{k}_{l=1} ([E_l]\F\{\Z^{\pm}\})\cap \F\{\Y\}
=\cap^{k}_{l=1}\sat(E^+_l),\end{equation}
 where   $E^+_l=\{\Z^{\f_{l,1}^+}-c_{l,1}\Z^{\f_{l,1}^-},
              \ldots, \Z^{\f_{l,s_l}^+}-c_{l,s_l}\Z^{\f_{l,s_l}^-}\}$,
              $l=1,\ldots,k$.
Since $E_{l}$ is regular and coherent, by Lemma \ref{lm-asc2},
$E^+_{l}$ is also regular and coherent.

Since $B\subset \F\{\Z\}$, we have the following zero decomposition
$$\V(\Y_0\cup B/\Y_1)= \V(\Y_0\cup(\{B\}:\mb)/\Y_1)
  \cup_{i=1}^s \V(\Y_0\cup B\cup\{y_{c_i}\}/\{y_{c_1},\ldots,y_{c_{i-1}}\}\cup\Y_1),$$
 where $\V(\Y_0\cup B\cup\{y_{c_i}\}/\{y_{c_1},\ldots,y_{c_{i-1}}\}\cup\Y_1)$ is
 further simplified with algorithm {\bf DECMONO} in step 2.8.
From \bref{eq-ma01},
 $$\V(\Y_0\cup\{B\}:\mb/\Y_1)=\cup^{k}_{l=1}\V(\sat(\Y_0,E^+_l)/\Y_1)
  =\cup^{k}_{l=1}\V([\Y_0,\sat(E^+_l)]/\Y_1),$$
where $\{\Y_0,E^+_l\}$ is a regular and coherent $\sigma$-chain
since $E^+_l$ does not contain variables in $\Y_0$. The above
formula explains why $(\{\Y_0,E^+_l\},\Y_1)$ is added to $\RB$ in
steps 2.5-2.7.

Let the algorithm returns $\RB=\{(\C_i,\Y_i);i=1,\ldots,m\}$. From
the above proof, we have $\V(F)=\cup^{k}_{l=1}\V(\sat(\C_i)/\Y_i)$.
Since $\Y_i\cap \C_i=\emptyset$ and $\sat(\C_i)$ is a prime
$\sigma$-ideal, the Cohn closure of $\V(\sat(\C_i)/\Y_i)$ is
$\V(\sat(\C_i))$ and hence
 $$\V(F)=\cup^{k}_{l=1}\V(\sat(\C_i)/\Y_i)=\cup^{k}_{l=1}\V(\sat(\C_i)).$$
By the difference Hilbert  Nullstellensatz,
 $$\{F\}=\cap^{k}_{l=1}\{\sat(\C_i)\} =\cap^{k}_{l=1}\sat(\C_i).$$
The algorithm terminates, since after each execution of step 2, in
the new components $(\Y_{0l},B_l,\Y_{1l})$ added to $\FB$ in step
2.8, $B_l$ contains at least one less variables than $B$. \qedd
%
%

The following example shows that even in the binomial case, a proper
irreducible $\sigma$-chain is not necessarily strong irreducible
\cite{gao-dcs}.
\begin{example}\label{ex-dec2}
Let $\A=\{A_1,A_2,A_3\}$ be a proper irreducible $\sigma$-chian,
where $A_1 = y_1^{x^k} - y_1^2$, $A_2 = y_2^{x^k} - y_2^2$, $A_3 =
y_1y_3^{3x} - y_2^x$, and $k\ge2$.
Then $A_3^{x^k}=y_1^2y_3^{3(x+k)}- y_2^{2x}\,\mod [A_1,A_2]$. Thus
$A_3^{x^k}-(y_1y_3^{3x} +
y_2^x)A_3=y_1^2(y_3^{3(x+k)}-y_3^{3x})\,\mod [A_1,A_2]$ is reducible
and $\A$ is not strong irreducible.
\end{example}

\begin{remark}
It is still open to give a minimal decomposition
for finitely generated binomial $\sigma$-ideals.
Or equivalently, the Ritt problem \cite[p 191]{kol}
is open even for binomial $\sigma$-chains.
\end{remark}

\section{Conclusion}
In this paper, we initiate the study of binomial $\sigma$-ideals and
toric $\sigma$-varieties. Two basic tools used are the $\Zx$-lattice
and the characteristic set instead of the $\Z$-lattice and the
Gr\"obner basis used in the algebraic case.  Since $\Zx$ is not a PID, a
matrix with entries in $\Zx$ does not have a Hermite normal form.
As an alternative, we introduce the concept of generalized Hermite
normal form which is equivalent to a reduced Gr\"obner basis for
$\Zx$-lattices. It is shown that a set of Laurent $\sigma$-binomials
is a regular and coherent $\sigma$-chain if and only if their
supports form a generalized Hermite normal form.
%

For Laurent binomial $\sigma$-ideals, three main results are proved.
Canonical representations for proper Laurent binomial
$\sigma$-ideals are given in terms of Gr\"obner basis of $\Zx$-lattices, regular and coherent $\sigma$-chains in $\F\{\Y^{\pm}\}$,
and partial characters over $\Zxn$.
We also give criteria for a
Laurent binomial $\sigma$-ideal to be reflexive, well-mixed, perfect, and prime
in terms of its support lattice. It is also shown that the reflexive, well-mixed,
and perfect closure of a Laurent binomial $\sigma$-ideal is still
binomial.
Finally, it is shown that a perfect Laurent $\sigma$-ideal
can be written as the intersection of Laurent reflexive prime
binomial $\sigma$-ideals with the same support lattice. Most of these
results are also extended to the $\sigma$-binomial case.

A toric $\sigma$-variety is defined as the Cohn closure of the image
of Laurent $\sigma$-monomial maps. Three characterizing properties
of toric $\sigma$-varieties are proved in terms of its coordinate
ring, its defining ideals, and group actions. In particular, a
$\sigma$-variety is toric if and only if its defining ideal is
a toric $\sigma$-ideal, meaning a binomial $\sigma$-ideal whose
support lattice is $\Zx$-saturated. 

Finally, algorithms are given for all the main results in the paper,
that is, to decide whether a finitely generated Laurent binomial
$\sigma$-ideal is reflexive, well-mixed, perfect, prime, or toric, and to decompose a
finitely generated perfect binomial $\sigma$-ideal as intersection
of reflexive prime binomial $\sigma$-ideals.


\begin{thebibliography}{99}

\bibitem{bouziane} 
 D. Bouziane, A. Kandri Rody, H. Ma$\hat{a}$rouf.
 Unmixed-dimensional Decomposition of a Finitely Generated Perfect Differential Ideal.
 {\em Journal of Symbolic Computation}, 31(6), 631-649, 2001.

\bibitem{cohen}
 H. Cohen.
 {\em A Course in Computational Algebraic Number Theory.}
 Springer-Verlag, Berlin, 1993.

\bibitem{cohn}
 R.~M. Cohn.
 {\em Difference Algebra}.
 Interscience Publishers,  New York,  1965.

\bibitem{cox-1998}
 D. Cox,  J. Little, D. O'Shea.
 {\em Using Algebraic Geometry},  Springer-Verlag,  New York,  1998.

\bibitem{cox-2010}
 D. Cox,  J. Little,  H. Schenck.
 {\em Toric Varieties},  Springer-Verlag,  New York,  2010.

\bibitem{dema}
 M. Demazure.
 Sous-groupes Alg\'ebriques de Rang Maximum du Groupe de Cremona,
 {\em Ann. Sci. \'Ecole Norm. Sup.}, {3}, 507-588, 1970.

\bibitem{es-bi}
 D. Eisenbud and B. Sturmfels.
 Binomial Ideals.
 {\em Duke Math. J.},  84(1), 1-45, 1996.

\bibitem{fulton}
 W. Fulton.
 {\em Introduction to Toric Varieties},
  Princeton Univ. Press, Princeton, USA, 1993.

\bibitem{gao}
 X.~S. Gao,  W. Li,  C.~M. Yuan.
 Intersection Theory in Differential Algebraic Geometry:
 Generic Intersections  and the Differential Chow Form.
 {\em Trans. of Amer. Math. Soc.},  365(9),  4575-4632,  2013.

\bibitem{gao-dcs}
 X.~S. Gao,  Y. Luo,  C.~M. Yuan.
 A Characteristic Set Method for Ordinary Difference Polynomial Systems.
 {\em Journal of Symbolic Computation},  44(3),  242-260,  2009.

\bibitem{gao-dcs1}
 X.~S. Gao, C. Yuan, and G. Zhang.
 Ritt-Wu's Characteristic Set Method
 for Ordinary Difference Polynomial Systems with Arbitrary Ordering.
 {\em Acta Mathematica Scientia}, 29(3,4), 1063-1080, 2009.

\bibitem{gao-im}
 X.~S. Gao and S.~C. Chou.
 Implicitization of Rational Parametric Equations.
 {\it Journal of Symbolic Computation}, 14, 459-470, 1992.

\bibitem{gelfand}
 I.~M. Gelfand,  M. Kapranov, A. Zelevinsky.
 {\em Discriminants, Resultants and Multidimensional Determinants}.
 Boston, Birkh\"auser, 1994.


\bibitem{GCG1992}
K.~O. Geddes, S.~R. Czapor, G. Labahn.
 {\em Algorithms for Computer Algebra}.
 Kluwer Academic Publishers, Boston, MA, 1992.


\bibitem{Hrushovski1}
 E. Hrushovski.
 The Elementary Theory of the Frobenius Automorphisms.
 Available from http://www.ma.huji.ac.il/\~\,ehud/, July, 2012.

\bibitem{Kei}
 K.~I. Iima and Y. Yoshino.
 Gr\"{o}bner Bases for the Polynomial Ring with Infinite Variables and Their
 Applications,
 {\em Communications in Algebra},  37:10,  3424-3437, 2009.

\bibitem{Kandri}
 A. Kandri Rody and  D. Kapur.
Computing a Gr\"{o}bner Basis of a Polynomial Ideal over a Euclidean Domain.
 {\em Journal of Symbolic Computation},  6(1), 37-57, 1988.

\bibitem{mum}
 G. Kempf, F. Knudsen, D. Mumford, B. Saint-Donat.
 {\em Toroidal Embeddings I},
 Lecture Notes in Math. 339, Springer-Verlag, Berlin, 1973.

\bibitem{kol}
 E. R. Kolchin.  {\em Differential Algebra and Algebraic Groups}.
 Academic Press,  New York and London, 1973.

\bibitem{levin}
 A. Levin.
 {\em Difference Algebra}.
 Springer-Verlag,  New Work, 2008.


\bibitem{dd-sres}
 W. Li, C.~M. Yuan, X.~S. Gao.
 Sparse Difference Resultant.
 {\em Proc. ISSAC 2013}, 275-282, ACM Press, New York, 2013.

\bibitem{dd-sres1}
 W. Li, C.~M. Yuan, X.~S. Gao.
 Sparse Difference Resultant.
 {\em Journal of Symbolic Computation},
 68, 169-203, 2015.

\bibitem{d-sres}
 W. Li, C.~M. Yuan, X.~S. Gao.
 Sparse Differential Resultant for Laurent Differential Polynomials.
 {\em Foundations of Computational Mathematics},
 15(2), 451-517, 2015.

\bibitem{dd-chowform}
 W. Li and Y.~H. Li. Difference Chow Form. 
 {\em Journal of Algebra}, 428, 67-90, 2015.

\bibitem{med1}
 A.  Medvedev and T. Scanlon.
 Invariant Varieties for Polynomial Dynamical Systems.
 {\em  Annals of Mathematics}, 179(1), 181-177, 2014.

\bibitem{miller1}
 E. Miller and B. Sturmfels.
 {\em Combinatorial Commutative Algebra}.
 Springer-Verlag, New York, 2005.

\bibitem{oda}
 K. Miyake and T. Oda.
 Almost Homogeneous Algebraic Varieties under Algebraic Torus Action,
 in {\em Manifolds-Tokyo}, (A. Hattori, ed.), 373-381,
 Univ. Tokyo Press, Tokyo, 1975.

\bibitem{oda-book}
 T. Oda.
 {\em Convex Bodies and Algebraic Geometry}.
 Springer, New York, 1988.

\bibitem{pachter1}
 L. Pachter and B. Sturmfels (eds).
 {\em Algebraic Statistics for Computational Biology}.
 Cambridge University Press, 2005.



\bibitem{ritt-dd1}
 J.~F. Ritt and J.~L. Doob.
 Systems of Algebraic  Difference Equations.
 {\em American Journal of Mathematics}, {55}(1), 505-514, 1933.

\bibitem{sata}
 I. Satake.
 On the Arithmetic of Tube Domains (blowing-up of the point at infinity),
 {\em Bull. Amer. Math. Soc.} {79}, 1076-1094, 1973.


\bibitem{put-dd}
 M. van der Put and M.~F. Singer.
 {\em Galois Theory of Difference Equations},
 Springer, Berlin and New York, 1997.

\bibitem{wibmer}
 M. Wibmer.
 Algebraic Difference Equations. Preprint, 2013.

\end{thebibliography}
\end{document}